\documentclass[11pt,a4paper]{article}
\usepackage{amsmath}
\usepackage{amssymb}
\usepackage{color}

\usepackage{graphicx}
\usepackage{amsfonts,amsmath, amssymb}

\numberwithin{equation}{section}

\usepackage{graphics}
\usepackage{epsfig}
\newtheorem{claim}{\bf \t}[part]

\newtheorem{theorem}{Theorem}[section]

\newtheorem{lemma}[theorem]{Lemma}
\newtheorem{proposition}[theorem]{Proposition}
\newtheorem{remark}[theorem]{Remark}

\def\v{\varepsilon}

\def\t{\theta}
\def\k{\kappa}
\def\n{\nu}
\def\m{\mu}
\def\a{\alpha}
\def\b{\beta}
\def\g{\gamma}
\def\d{\delta}
\def\l{\lambda}

\def\r{\rho}

\def\Om{\Omega}

\def\i{\infty}
\def\f{\frac}
\def\MZ{\mathcal{Z}}
\def\MH{\mathcal{H}}

\begin{document}

\title{Uniform regularity and vanishing dissipation limit for  the full compressible Navier-Stokes system  in  3-D bounded domain}

\author{ {Yong Wang$^{\S}$\footnote{ \indent
			Email addresses:  yongwang@amss.ac.cn(Yong Wang)
		} 
	}
\\
  \ \\
	{\small \it $^\S$Institute of Applied Mathematics, AMSS, CAS, Beijing 100190, China}
}

\date{ }
\maketitle


\begin{abstract}
In the present paper, we study the uniform regularity and vanishing dissipation limit for the full compressible Navier-Stokes system whose viscosity and heat conductivity  are allowed to vanish at different order. The problem is studied in a 3-D bounded domain with Navier-slip type boundary conditions \eqref{1.9}. It is shown that there exists a unique strong solution to the full compressible Navier-Stokes system with the boundary conditions \eqref{1.9} in  a finite time interval   which is independent of the viscosity and heat conductivity. The solution is uniform bounded in $W^{1,\infty}$ and a conormal Sobolev space. Based on such uniform estimates, we prove the convergence of the solutions of the full compressible  Navier-Stokes to the corresponding solutions of the full compressible Euler system in $L^\infty(0,T;L^2)$,$L^\infty(0,T;H^1)$ and $L^\infty([0,T]\times\Omega)$ with a rate of convergence.

\


\

AMS: 35Q35, 35B65, 76N10
\end{abstract}


\section{Introduction and Main Results}
The motion of a compressible viscous, heat conductive, ideal polytropic fluid is governed by the following full compressible Navier-Stokes equations(FCNS)
{\small\begin{equation}\label{1.1}
\begin{cases}
\r_t^\v+\mbox{div}(\r^\v u^\v)=0,\\
(\r^\v u^\v)_t+\mbox{div}(\r^\v u^\v \otimes u^\v)+\nabla p^\v=\mu\v \Delta u^\v+(\mu+\l)\v \nabla div u^\v,\\
(\r^\v E^\v)_t+\mbox{div}(\r^\v u^\v E^\v+p^\v u^\v)=\k(\v) \Delta\t^\v+\mbox{div}(\vec{\tau}^\v u^\v),
\end{cases}
x\in \Omega,~ t>0
\end{equation}}
where $\Omega$ is a bounded smooth domain of $\mathbb{R}^3$. Here $\r^\v, u^\v$ and $E^\v$ represent density, velocity and total energy, respectively. The pressure function $p^\v$, total energy $E^\v$ are  given by
{\small\begin{equation*}
p^\v=R\r^\v\t^\v, ~~E^\v=c_{v}(\t^\v+\f12|u^\v|^2),
\end{equation*}}
where $\t^\v$  is temperature and $c_v$  is a positive constant. For the simplicity of presentation, we normalize $c_v$ to be 1. The tensor $\vec{\tau}^\v$ is represented
{\small\begin{equation*}
 \vec{\tau}^\v=\l\v\mbox{div}u^\v I+2\mu\v Su^\v,~~\mbox{with}~Su^\v=\f12(\nabla u^\v+(\nabla u^\v)^T).
\end{equation*}}
Here  $\mu, \l$ are given constants satisfying the following physical restriction
{\small\begin{equation}
\mu>0,~~2\m+3\l>0,
\end{equation}}
and the parameter $\v>0$ is the inverse of the Reynolds number. $\k(\v)>0$ is the heat conductivity  which is  assumed to depend on $\v$.

We impose the full compressible Navier-Stokes equations with the following Navier-slip type boundary conditions
{\small
\begin{equation}\label{1.9}
		u^\v\cdot n=0,~~
		((Su^\v) n)_\tau=(Au^\v)_\tau,~~\mbox{and}~~
		n\cdot\nabla\t^\v=\n\t^\v,~~
		\mbox{on}~\partial\Omega.
	\end{equation}
}
where $n$ is the outward  unit normal to $\partial\Omega$, $u_\tau$ represents the tangential part of $u$, $A$ is a smooth symmetric  matrix and $\nu$ is a given constant.  For smooth solutions, it is noticed that
{\small\begin{equation}\nonumber
\left(2S(v)n-(\nabla\times v)\times n\right)_\tau=-(2S(n)v)_\tau,
\end{equation}}
see \cite{Xiao-Xin-3} for details. Then the boundary condition \eqref{1.9} can  be rewritten in the form of the vorticity as
{\small\begin{equation}\label{1.9-1}
u^\v\cdot n=0,~~
n\times \omega^\v=[Bu^\v]_\tau,~~\mbox{and}~~
n\cdot\nabla\t^\v=\n\t^\v,~~
\mbox{on}~\partial\Omega.
\end{equation}}
where $\omega^\v=\nabla\times u^\v$ is the vorticity   and  $B=2(A-S(n))$ is a symmetric matrix. Actually,  it turns out  that the form \eqref{1.9-1} will be  more convenient than \eqref{1.9} in the energy estimates.

We are interested in the existence of strong solutions of \eqref{1.1} with uniform bounds on an interval of time independent of the viscosity and heat conductivity, and the vanishing dissipation limit to the full compressible Euler flow  as $\v$ and $\k(\v)$ vanish, i.e,
{\small\begin{equation}\label{1.7}
\begin{cases}
\r_t+\mbox{div}(\r u)=0,\\
(\r u)_t+\mbox{div}(\r u \otimes u)+\nabla p=0,\\
[\r(\t+\f12|u|^2)]_t+\mbox{div}[\r u(\t+\f12|u|^2)+pu]=0,
\end{cases}
\mbox{as}~\v,\k\rightarrow 0+,
\end{equation}}
with slip boundary condition
\begin{equation}\label{1.7-2}
u\cdot n|_{\partial\Omega}=0.
\end{equation}
There has lots of literatures on the inviscid limit for {\it incompressible} Navier-Stokes equations. The inviscid limit of Cauchy problem has been studied by many authors,  see for instances \cite{Constantin,Constantin-1,Kato-1,Masmoudi}.
However, in the presence of a physical boundary, the problems become challenging due to the appearance of  boundary layers. As illustrated   by Prandtl's theory,  the inviscid limit of the incompressible Navier-Stokes with no-slip boundary condition to the incompressible Euler flows with slip boundary condition $\eqref{1.7-2}$ is a very difficult problem.  Sammartino-Caflisch \cite{S-C-1,S-C-2}  proved the  convergence of the incompressible Navier-Stokes flows to the Euler flows away from the boundary and to the prandtl flows near the boundary at the inviscid limit for the analytic initial data. Recently, Maekawa\cite{Maekawa} proved such limit when the initial vorticity is located away from the boundary  in the 2-D half plane. On the other hand, for the incompressible Navier-Stokes system with Navier-slip boundary condition
\eqref{1.9}(without the heat flux part), lots of important progress has been made on this problem. The uniform $H^3$ bound and the inviscid limit to Euler flow was proved by Xiao-Xin\cite{Xiao-Xin-1} for flat boundaries which was generalized to $W^{k,p}$ in \cite{Beiro-1,Beiro-2} by Veiga-Crispo soon later.  However, one can not obtain such results for general curved boundaries since boundary layers may appear due to non-trivial curvature as pointed out by Iftimie and Sueur\cite{Iftimie}, where the inviscid limit  was also obtained in $L^\infty(0,T,; L^2)$ by a careful construction of boundary layer expansions. In order  to investigate precisely the asymptotic structure and get the convergence in stronger norms such as $L^\infty(0,T; H^s)(s>0)$, stronger  estimates  are needed. Recently, Masmoudi-Rousset \cite{Masmoudi-R} established a conormal uniform estimates for 3-D   domains with the Naiver-slip boundary condition, which implies the uniform boundedness of Lipschitz-norm for the velocity field. This allows  to obtain the inviscid limit in $L^\infty$-norm by a compactness argument. Based on the uniform estimates in \cite{Masmoudi-R}, better convergence with rates have been obtained in \cite{Gie-Kelliher} and \cite{Xiao-Xin-2}. In particular, Xiao-Xin  \cite{Xiao-Xin-2} has proved the convergence in $L^\infty(0,T; H^1)$ with a rate.

For the {\it isentropic  compressible} Navier-Stokes equations,  Xin-Yanagisawa \cite{Xin-Y} studied the vanishing viscosity limit of the linearized problem  with the no-slip boundary condition in 2-D half plane. For the Navier-slip boundary condition case, Wang-Williams \cite{Wang-Williams} constructed a boundary layer solution of the compressible Navier-Stokes equations  in 2-D half plane. The layers constructed in \cite{Wang-Williams} are of width $O(\sqrt\v)$ as the Prandtl boundary layer, but the amplitude are of $O(\sqrt\v)$ which is similar to the one \cite{Iftimie} for the incompressible case. It is also shown \cite{Wang-Williams} that  the boundary layers for the density is weaker than the one for the velocity.  So, in general, it is impossible to obtain the $H^3$ or $W^{2,p}(p>3)$ estimates for the compressible Navier-Stokes system with the Navier-slip boundary condition. Recently, Paddick \cite{Paddick} obtained an uniform conormal Sobolev estimates for  the  isentropic compressible Navier-Stokes system in the 3-D half-space.  Wang-Xin-Yong \cite{Wang-Xin-Yong} also obtained an uniform regularity for isentropic compressible Navier-Stokes equations with Navier-slip boundary conditions in 3-D domain with curvature, especially, the inviscid limit was also obtained with rate of convergence in $L^\infty([0,T]\times\Omega)$ and $L^\infty([0,T];H^1)$.  The fact that the boundary layer for density is weaker than the one for velocity fields was also shown in \cite{Wang-Xin-Yong}.

For the {\it full  compressible} Navier-Stokes equations, the study is quite limit. Under the assumption  that  the viscosity and heat conductivity converge to zero at the same order,  Ding-Jiang \cite{Ding-Jiang} studied the zero viscosity and heat conductivity limit for the linearized compressible Navier-Stokes-Fourier equations with no-slip boundary condition  in the half plane.

However, there is no uniform regularity and vanishing dissipation limit results for the {\it full compressible} Navier-stokes equations \eqref{1.1} with Navier-slip type boundary conditions \eqref{1.9} in a bounded domain. The aim of this paper is to investigate the uniform regularity for the solutions of the full compressible Navier-Stokes system \eqref{1.1}  even if  the viscosity and heat conductivity  converge to zero at  different order.  Compared to the isentropic case\cite{Wang-Xin-Yong,Paddick}, it is difficult to obtain the Lipschitz estimates for the solutions of \eqref{1.1} due to the appearance of temperature and the strongly coupled system of $(\r^\v,u^\v,\t^\v)$.   On the other hand, the amplitude and width of  boundary layers for velocity and  temperature may be not at the same order if the viscous and heat conductivity vanish at different order. Then  the interaction of the two different amplitude boundary layers may  arise difficulties  in the analysis of obtaining uniform regularity, especially, in the Lipschitz estimates.  To overcome these difficulties, some new ideas and observations are need.  It is also very important to study the vanishing dissipation limit. Especially,  we shall investigate how the rate of convergence is influenced by the thermal boundary layers.

Before stating our main results, we first explain the notations and conventions used throughout this paper.
Similar to \cite{Masmoudi-R,Wang-Xin-Yong}, one assumes that the bounded domain $\Omega\subset\mathbb{R}^3$  has   a covering such that
{\small\begin{equation}\label{2.0}
	\Omega\subset\Omega_0\cup_{k=1}^n\Omega_k,
\end{equation}}
where $\overline{\Omega}_0\subset\Omega$ and in each $\Omega_k$ there exists a function $\psi_k$ such that
{\small\begin{equation}
	 \Omega\cap\Omega_k=\{x=(x_1,x_2,x_3)~|~x_3>\psi_k(x_1,x_2)\}\cap\Omega_k
	 ~~\mbox{and}~~\partial\Omega\cap\Omega_k=\{x_3=\psi_k(x_1,x_2)\}\cap\Omega_k.\nonumber
\end{equation}}
$\Omega$ is said to be  $\mathcal{C}^m$ if the functions $\psi_k$ are $\mathcal{C}^m$-function. To define the Sobolev conormal spaces, one considers $(Z_k)_{1\leq k\leq N}$ a finite set of generators of vector fields that are tangent to $\partial\Omega$ and set
{\small\begin{equation*}\label{1.3}
	H^m_{co}=\Big\{f\in L^2(\Omega)~|~Z^If\in L^2(\Omega), ~~\mbox{for}~|I|\leq m  \Big\},
\end{equation*}}
where $I=(k_1,\cdots, k_m)$. The following notations will be used:
{\small\begin{equation*}\label{1.4}
	\|u\|^2_{m}=\|u\|^2_{H^m_{co}}=\sum_{j=1}^3\sum_{|I|\leq m}\|Z^Iu_j\|^2_{L^2},~~\|u\|^2_{m,\infty}=\sum_{|I|\leq m}\|Z^Iu\|^2_{L^\infty}~~\mbox{and}~~\|\nabla Z^m u\|^2=\sum_{|I|=m}\|\nabla Z^Iu\|^2_{L^2}.
\end{equation*}}
Noting that by using the covering of $\Om$, one can always assume that each vector field $(\rho,u,\t)$ is supported in one of the $\Om_i$, moreover, in $\Om_0$ the norm $\|\cdot\|_m$ yields a control of the standard $H^m$ norm, whereas if $\Om_i\cap \partial\Om\neq \O$, there is no control of the normal derivatives.

Denote by  $C_k$   a positive constant independence of $\v,\k\in(0,1] $ which depends only on the $\mathcal{C}^k$-norm of the functions $\psi_j,j=1\cdots,n$. Since $\partial\Omega$ is given locally by $x_3=\psi(x_1,x_2)$(we omit the subscript $j$ for notational convenience), it is convenient to use the coordinates:
{\small\begin{equation*}\label{3.5}
	\Psi:~(y,z)\longmapsto (y,\psi(y)+z)=x.
\end{equation*}}
A local basis is thus given by the vector fields $(e_{y^1},e_{y^2},e_z)$ where $e_{y^1}=(1,0,\partial_1\psi)^t$, $e_{y^2}=(0,1,\partial_2\psi)^t$ and $e_{z}=(0,0,-1)^t$. On the boundary $e_{y^1}$ and $e_{y^2}$ are tangent to $\partial\Omega$, and in general, $e_z$ is not a normal vector field. By using this parametrization, one can take as suitable vector fields compactly supported in $\Omega_j$ in the definition of the $\|\cdot\|_m$ norms:
{\small\begin{equation}\label{3.6}
	 Z_i=\partial_{y^i}=\partial_i+\partial_i\psi\partial_z,~i=1,2, ~~Z_3=\varphi(z)\partial_z,
\end{equation}}
where $\varphi(z)=\f{z}{1+z}$ is smooth, supported in $\mathbb{R}_+$ with the property $\varphi(0)=0$,
$\varphi'(0)>0$, $\varphi(z)>0$ for $z>0$. It is easy to check that
{\small\begin{equation*}\label{3.7}
	Z_kZ_j=Z_jZ_k,~~j,~k=1,2,3,
\end{equation*}}
and
{\small\begin{equation*}\label{3.8}
	 \partial_zZ_i=Z_i\partial_z,~i=1,2,~~\mbox{and}~~\partial_zZ_3\neq Z_3\partial_z.
\end{equation*}}

In this paper, we shall still denote by $\partial_j,~j=1,2,3$ or $\nabla$ the derivatives in the physical space. The coordinates of a vector field $u$ in the basis  $(e_{y^1},e_{y^2},e_z)$ will be denoted by $u^i$, thus
{\small\begin{equation}\label{3.9}
	u=u^1e_{y^1}+u^2e_{y^2}+u^3e_{z}.
\end{equation}}
We shall denote by $u_j$ the coordinates in the standard basis of $\mathbb{R}^3$, i.e, $u=u_1e_1+u_2e_2+u_3e_3$. Denote by $n$ the unit outward normal in the physical space which is given locally by
{\small\begin{equation}\label{3.10}
	n(x)\equiv n(\Psi(y,z))=\f{1}{\sqrt{1+|\nabla\psi(y)|^2}}\left(\begin{array}{cccc} &\partial_1\psi(y)\\&\partial_2\psi(y)\\&-1\end{array}\right)\doteq\f{-N(y)}{\sqrt{1+|\nabla\psi(y)|^2}},
\end{equation}}
and by $\Pi$ the orthogonal projection
{\small\begin{equation}\label{3.10-1}
	\Pi u\equiv\Pi(\Psi(y,z))u=u-[u\cdot n(\Psi(y,z))]n(\Psi(y,z)),
\end{equation}}
which gives the orthogonal projection onto the tangent space of the boundary. Note that $n$ and $\Pi$ are defined in the whole $\Omega_k$ and do not depend on $z$.

For later use and notational convenience, set
{\small\begin{equation}\label{2.1}
	 \mathcal{Z}^\a=\partial_t^{\a_0}Z^{\a_1}=\partial_t^{\a_0}Z_1^{\a_{11}}Z_2^{\a_{12}}Z_3^{\a_{13}}.
\end{equation}}
where $\a, \a_0, \a_1$ are the differential  multi-indices with $\a\doteq(\a_0,\a_1),~\a_1=(a_{11},\a_{12},\a_{13})$,
and we also use the following notations
{\small\begin{equation}\label{2.3}
	\|f(t)\|^2_{\mathcal{H}^m}=\sum_{|\a|\leq m}\|\mathcal{Z}^\a f(t)\|^2_{L^2_x},~~
	\|f(t)\|_{\mathcal{H}^{k,\infty}}=\sum_{|\a|\leq k}\|\mathcal{Z}^\a f(t)\|_{L^\infty_x},
\end{equation}}
for smooth space-time function $f(x,t)$.

Firstly, we consider the uniform regularity of the solutions of full compressible Navier-Stokes system \eqref{1.1} with the Navier-slip type boundary conditions \eqref{1.9}. Since the viscous and thermal boundary layers  may appear in the presence of physical boundaries, so one needs to design a suitable functional space. Here  the functional space $X_m^\v(T)$ for   functions $(\r,u,\t)=(\r,u,\t)(x,t)$  is defined as
{\small\begin{eqnarray}
X^\v_m(T)=\Big\{(\r,u,\t)\in L^\infty([0,T],L^2);~~
\mbox{esssup}_{0\leq t\leq T}\|(\r,u,\t)(t)\|_{X_m^\v}<+\infty\Big\},
\end{eqnarray}}
where the norm $\|(\cdot,\cdot,\cdot)\|_{X_m^\v}$ is given by
{\small\begin{align}\label{1.11-1}
&\|(\r,u,\t)(t)\|^2_{X^\v_m}=\|(\r,u,\t)(t)\|^2_{\mathcal{H}^m}+\|\nabla u(t)\|^2_{\mathcal{H}^{m-1}}+\sum_{k=0}^{m-2}\|\partial_t^k \nabla(\r,\t)(t)\|^2_{m-1-k}+\v\|\nabla\partial_t^{m-1}\r(t)\|^2\nonumber\\
&~~~~~~~~~~~~~~~~~~~~~~~~~~~~~~+\v\|\nabla\partial_t^{m-2}\mbox{div}u(t)\|^2+\k(\v) \|\Delta\partial_t^{m-2}\t(t)\|^2+\|\nabla(\r,u,\t)(t)\|^2_{\mathcal{H}^{1,\infty}}\nonumber\\[1mm]
&~~~~~~~~~~~~~~~~~~~~~~~~~~~~~~~~~~~~~~~~+\v\|\nabla^2u(t)\|^2_{L^\infty}
+\v\|\nabla(\r\t)(t)\|^2_{\mathcal{H}^{2,\infty}}.
\end{align}}
We remark that the term $\v\|\nabla(\r\t)(t)\|^2_{\mathcal{H}^{2,\infty}}$ included in \eqref{1.11-1} is important for us to obtain the Lipschitz estimates  even though such term is slightly strange. And we will explain the reason of including such term after Theorem \ref{thm1.1} below.

In the present paper,  we supplement the full compressible Navier-Stokes equations \eqref{1.1} with the initial data
{\small\begin{equation}\label{1.11}
(\r^\v,u^\v,\t^\v)(x,0)=(\r_0^\v,u_0^\v,\t^\v)(x),
\end{equation}}
such that
{\small\begin{align}\label{1.13}
\sup_{0< \v\leq 1}\|(\r_0^\v,u_0^\v,\t^\v_0)\|^2_{X^\v_m}
\leq \tilde{C}_0,~~~0<\hat{C}_0^{-1}\leq \r_0^\v,~ \t_0^\v\leq \hat{C}_0<\infty
\end{align}}
where $\hat{C}_0>0$, $\tilde{C}_0>0$ are positive constants independent of $\v\in(0,1]$, and  the time derivatives of initial data  in \eqref{1.13} are defined through the full compressible Navier-Stokes system \eqref{1.1}. Thus,   the initial data  $(\r^\v_0,u^\v_0,\t^\v_0)$  is assumed to have a higher space regularity and compatibilities. Notice that the {\it a priori} estimates in Theorem \ref{thm3.1} below is obtained in the case that the approximate solution is sufficient smooth up to the boundary, therefore, in order to obtain a selfcontained result, one  needs  to assume that the approximate initial data  satisfies the boundary compatibility conditions, i.e. \eqref{1.9}(or equivalent to \eqref{1.9-1}).  For the initial data $(\r_0^\v,u_0^\v,\t_0^\v)$  satisfying \eqref{1.13}, it is  not clear if there exists an approximate sequence  $(\r_0^{\v,\d},u_0^{\v,\d},,\t^{\v,\d}_0)$($\d$ being a regularization parameter), which satisfy the boundary compatibilities  and $\|(\r_0^{\v,\d}-\r_0^\v,u_0^{\v,\d}-u_0^\v,\t_0^{\v,\d}-\t_0^\v)\|_{X^\v_m}\rightarrow0$ as $\d\rightarrow0$. Therefore, we set
{\small\begin{eqnarray}\label{Initial}
& X_{NS,ap}^{\v,m}=\Big\{(\r,u,\t)\in H^{3m}(\Omega)~\Big| \partial_t^k\r, \partial_t^k\t,~\partial_t^ku,k=1,\cdots,m~\mbox{are defined through the system} \eqref{1.1} \nonumber\\
&~~~~~~~~~~~~ \mbox{and}~\partial_t^ku, \partial_t^k\t,k=0,\cdots,m-1~ \mbox{satisfy   the boundary compatibility conditions}\Big\},
\end{eqnarray}}
and
{\small\begin{equation}\label{Initial space}
 X^{\v,m}_{NS}=\mbox{The closure of}~ X^{\v,m}_{NS,ap}~ \mbox{in the norm }~ \|(\cdot,\cdot,\cdot)\|_{X^\v_m}.
\end{equation}}

If the heat conductivity $\k(\v)$ decays too fast as $\v\rightarrow0+$, then the possible interaction between the viscous boundary layers and thermal  boundary layers is strong and it is  hard to get the uniform regularity. So, in order to control the possible interaction between the viscous boundary layers and the thermal boundary layers, throughout this paper, we assume that the heat conductivity is a continuous function of $\v$ and  satisfies
{\small \begin{align}\label{1.8}
 \v^4\leq C\k(\v)<\infty,~~\mbox{for}~~\v\in (0,1],
 \end{align}}
where the $C>0$ is some positive constant.  Then  our uniform regularity result is as follows:
\begin{theorem}[Uniform Regularity]\label{thm1.1}
Let $m$ be an integer satisfying $m\geq 6$, $\k(\v)$ satisfies \eqref{1.8},  $\Omega$ be a $\mathcal{C}^{m+2}$ domain and $A\in\mathcal{C}^{m+1}(\partial\Omega)$. Consider the initial data $(\r^\v_0, u^\v_0,\t^\v_0)\in X^{\v,m}_{NS}$ given in \eqref{1.11} and satisfying \eqref{1.13}. Then there exists a time $T_0>0$ and $\tilde{C}_1>0$  independent of  $\v\in(0,1]$, such that there exists a unique solution $(\r^\v, u^\v,\t^\v)$ of \eqref{1.1}, \eqref{1.9} and \eqref{1.11}  on $[0,T_0]$  and satisfies the estimates:
{\small\begin{equation}\label{1.19}
(2C_0)^{-1}\leq \r^\v(t), \t^\v(t)\leq 2C_0~~\forall t\in [0,T_0],
\end{equation}}
and
{\small\begin{align}\label{1.18}
&\sup_{0\leq t\leq T_0}\|(\r^\v,u^\v,\t^\v)(t)\|^2_{X^\v_m}+\int_0^{T_0}\|\nabla\partial_t^{m-1}(\r^\v,\t^\v)(t)\|^2+\|\nabla (\r^\v\t^\v)(t)\|^2_{\MH^{2,\infty}}dt+\v\int_0^{T_0}\|\nabla u^\v(t)\|^2_{\mathcal{H}^{m}}dt\nonumber\\
&~~~+\k(\v)\int_0^{T_0}\|\nabla \t^\v(t)\|^2_{\mathcal{H}^{m}}dt+\sum_{k=0}^{m-2}\int_0^{T_0}\v\|\nabla^2\partial_t^{k}u^\v(t)\|^2_{m-1-k}+\k(\v)\|\Delta\partial_t^{k}\t^\v(t)\|^2_{m-1-k}dt
\nonumber\\
&~~~+\v^2\int_0^{T_0}\|\nabla^2\partial_t^{m-1}u^\v(t)\|^2dt+\k(\v)^2\int_0^{T_0}\|\Delta\partial_t^{m-1}\t^\v(t)\|^2+\|\nabla\MZ^{m-2}\Delta\t^\v(t)\|^2dt\leq \tilde{C}_1<\infty,
\end{align}}
where $\tilde{C}_1$ depends only  on $\hat{C}_0,~\tilde{C}_0$ and $C_{m+2}$.
\end{theorem}
\begin{remark}
The novelty of this work is that we allow the viscous and heat conductivity to vanish at  different order.  And it is noted that there are many functions satisfy the condition \eqref{1.8}. For example, it is easy to see that \eqref{1.8} holds provided
$\k(\v)=\v^b$ where $b$ is constant such that $0\leq b\leq 4$.
\end{remark}
\begin{remark}
In order to obtain the  Lipschitz estimates included in \eqref{1.11-1}, one needs to use the pointwise estimates because the boundary layers prevent to obtain uniform estimate in $H^3(\Omega)$(or $W^{2,p}, p>0$). So, one has to deal with the possible interaction of the viscous and thermal  boundary layers in the pointwise estimates. Indeed,  the restriction \eqref{1.8} is used to control such possible interaction, see Lemma \ref{lem3.17}-Lemma \ref{lem3.19} for details.
\end{remark}
\begin{remark}
For the solution $(\r^\v,u^\v,\t^\v)(t)$ of \eqref{1.1},\eqref{1.9-1},\eqref{1.11}, the boundary conditions \eqref{1.9}(or equivalently \eqref{1.9-1}) are satisfied in the trace sense for every fixed $\v\in(0,T_0]$ and $t\in(0,T_0]$.
\end{remark}

We now describe the main ideas of the proof of Theorem \ref{thm1.1}. And it turns out that it suffices to establish the estimates \eqref{1.18}. It is noted that there are two parts included in \eqref{1.18}, i.e., the conormal energy estimates part and the pointwise estimates part.  Firstly, by complicated  conormal energy estimates, one can obtain
{\small\begin{align}\label{1.11-10}
	&\|(\r^\v,u^\v,\t^\v)(t)\|^2_{\mathcal{H}^m}+\sum_{k=0}^{m-2}\|\partial_t^k \nabla(\r^\v,u^\v,\t^\v)(t)\|^2_{m-1-k}+\|\partial_t^{m-1}\omega^\v(t)\|^2+\v\|\nabla\partial_t^{m-1}\r^\v(t)\|^2+\v\int_0^t\|\nabla u^\v(\tau)\|^2_{\mathcal{H}^{m}}\nonumber\\
	 &+\k\int_0^t\|\nabla\t^\v(\tau)\|^2_{\mathcal{H}^{m}}+\sum_{k=0}^{m-2}\int_0^t\|\partial_t^{k}(\sqrt{\v}\nabla^2u^\v,\sqrt{\k}\Delta\t^\v)(\tau)\|^2_{m-1-k}+\int_0^t\v^2\|\nabla^2\partial_t^{m-1}u^\v(\tau)\|^2,
	\end{align}}
at the cost of
{\small\begin{equation}\label{1.11-11}
\int_0^t\|\nabla\MZ^{m-1}\mbox{div}u^\v\|^2+\k^2\|\nabla\MZ^{m-2}\Delta\t^\v\|^2d\tau~\mbox{and}~\int_0^t\|\partial_t^{m-1}\nabla(\r^\v,\t^\v)\|^2d\tau,
\end{equation}}
see Lemma \ref{lem2.2}-Lemma \ref{lem4.3} and Lemma \ref{lem5.1} below. By using the structure of mass equation  and energy equation, respectively, one can bound the first part of \eqref{1.11-11} at the cost of $\int_0^t\|\partial_t^{m-1}\nabla(\r^\v,\t^\v)\|^2d\tau$, see Lemma \ref{lem4.5} below. So, it suffices to bound the second part of \eqref{1.11-11}. Considering $\int_0^t\int\partial_t^{m-2}\eqref{2.5}_3\cdot\nabla\partial_t^{m-2}\mbox{div}u^\v+\partial_t^{m-1}\nabla\eqref{2.5}_2\cdot
\f{\nabla\partial_t^{m-1}\t^\v}{\t^\v}dxd\tau$, one can control the second part of \eqref{1.11-11} and $\v\|\nabla\partial_t^{m-2}\mbox{div}u^\v(t)\|^2+\k(\v) \|\Delta\partial_t^{m-2}\t^\v(t)\|^2$, see Lemma \ref{lem4.8} below. Therefore, combining the above estimates, one can obtain the conormal energy estimates of \eqref{1.18} except $\|\nabla\partial_t^{m-1}u^\v\|^2$, see also \eqref{5.39} below.

Next, we try to establish the pointwise estimates part. However, it is difficult to obtain such estimates because the equations of $\r^\v,u^\v,\t^\v$ are strongly coupled and the viscosity and  heat conductivity are not at the same order. Actually, if one estimate $\|\nabla\r^\v\|_{\MH^{1,\infty}}$ directly, then one has  to deal with the high order derivative term $\int_0^t\|\nabla\t^\v\|^2_{\MH^{2,\infty}}d\tau$, however it is hard to control this term in our functional space. Instead, we try to control the pointwise estimates of $\nabla(\r^\v\t^\v), \nabla\t^\v$ and $\nabla u^\v$. This is key to overcome the difficulty. Indeed, we can obtain(see Lemma \ref{lem3.18} below)
{\small\begin{equation}
\|\nabla(\r^\v\t^\v)(t)\|^2_{\MH^{1,\infty}}+\v\|\nabla(\r^\v\t^\v)(t)\|^2_{\MH^{2,\infty}}+\int_0^t\|\nabla(\r^\v\t^\v)(\tau)\|^2_{\MH^{2,\infty}}\leq \v\int_0^t\|\nabla\t^\v(\tau)\|^2_{\MH^{3,\infty}}d\tau+\cdots.
\end{equation}}
where $\cdots$ means terms can be controlled.
Since the strength and width of the thermal boundary layers is connected with $\k$, then the first term on the RHS of above is actually the interaction of viscous   and  thermal boundary layers. If the heat conductivity $\k$ vanishes too fast, it will be very hard to control  such interaction term. To overcome the difficulty, we assume that the decay rate of $\k$ satisfies \eqref{1.8}. Then   the interaction term can be controlled as following:
{\small\begin{equation}
\v\int_0^t\|\nabla\t^\v(\tau)\|^2_{\MH^{3,\infty}}d\tau\leq \v^4\int_0^t\|\Delta\t^\v(\tau)\|^2_{\MH^{4}}d\tau+\cdots\leq \k\int_0^t\|\Delta\t^\v(\tau)\|^2_{\MH^{4}}d\tau+\cdots,
\end{equation}}
where the last term of above has already been controlled in the conormal energy estimates part. It is worth to point out that the above interaction estimate will be employed repeatedly throughout the pointwise estimates part. On the other hand, to control the pointwise estimate of  $\nabla\t^\v$, the most difficult part is to deal with the term $p^\v\nabla\mbox{div}u^\v$ which comes from the term $p^\v\mbox{div}u^\v$ on the LHS of energy equation(see also \eqref{6.105} below). Actually, if  $p^\v\nabla\mbox{div}u^\v$ is regarded as a source term,  it will be very hard to control $\int_0^t\|p^\v\nabla\mbox{div}u^\v(\tau)\|^2_{\MH^{1,\infty}}d\tau$ because the derivative is too high.
We remark that such difficulty does not arise in the isentropic case \cite{Wang-Xin-Yong}. So, to overcome the difficulty, some new idea is needed. Fortunately, we find that the term  $p^\v\nabla\mbox{div}u^\v$ can be represented as following
{\small\begin{align}\label{6.107-1}
p^\v\nabla\mbox{div}u^\v=R\r^\v[\nabla\t^\v_t+(u^\v\cdot\nabla)\nabla\t^\v]-R[\nabla(\r^\v\t^\v)_t+(u^\v\cdot\nabla)\nabla(\r^\v\t^\v)]+\mbox{lower order terms}.
	\end{align}}
It is noted that the first part(i.e., the hardest part) on the RHS of \eqref{6.107-1} can be absorbed into the main part of equation(see \eqref{6.109} below), while the second part on the RHS of \eqref{6.107-1} can be regarded as a source term because the term $\int_0^t\|\nabla(\r^\v\t^\v)(\tau)\|^2_{\MH^{2,\infty}}$ has already been controlled above. This observation is key to close the pointwise estimates. And it is also one of the main reason to include $\v\|\nabla(\r^\v\t^\v)\|^2_{\MH^{2,\infty}}$ in our functional space. Based on  the above observation, one can obtain the control of $\|\nabla\t^\v\|^2_{\MH^{1,\infty}}$. Later, by similar arguments as \cite{Wang-Xin-Yong}, one can control $\|\nabla u^\v\|^2_{\MH^{1,\infty}}+\v\|\nabla^2u^\v\|^2_{L^\infty}$. Finally, in order  to estimate $\|\nabla\partial_t^{m-1}u^\v\|^2$, we still need to obtain the uniform bound of $\|\partial_t^{m-1}\mbox{div}u^\v\|^2$ which is hard to get by the conormal energy estimate because  some  boundary terms are hard to control. By employing the mass equation, it is found that $\|\partial_t^{m-1}\mbox{div}u^\v\|^2$ can be controlled by the conormal energy estimates  and pointwise estimates obtained above. Therefore, combining all the above estimates, one proves \eqref{1.18}.

\

Based on the uniform estimates  in Theorem 1.1, using a similar arguments as \cite{Wang-Xin-Yong,Masmoudi-R}, one can  justify the  vanishing dissipation limit of solutions of full compressible Navier-Stokes system \eqref{1.1} to the solutions of the full  Euler equations \eqref{1.7} in $L^\infty$-norm by a strong compactness argument, but without convergence rate. In the present paper, we are interested in the vanishing dissipation limit with rates of convergence.

We supplement the full  Euler equations \eqref{1.7} and the full compressible Navier-Stokes equations \eqref{1.1} with the same initial data $(\r_0,u_0,\t_0)$ satisfying
{\small\begin{equation}\label{1.20}
0<\hat{C}_0^{-1}\leq\r_0,\t_0\leq \hat{C}_0~~\mbox{and}~~ (\r_0,u_0,\t_0)\in H^3\cap X^{\v,m}_{NS}~\mbox{with}~m\geq6.
\end{equation}}
It is well known that there exists a unique smooth solution $(\r,u,\t)\in H^3$ for the problem \eqref{1.7}, \eqref{1.7-2} with initial data  $(\r_0,u_0,\t_0)$ at least locally in time $[0,T_1]$ where $T_1>0$ depends only on $\|(\r_0,u_0,\t_0)\|_{H^3}$. On the other hand, it follows from Theorem \ref{thm1.1} that there exists a time $T_0$ and $\tilde{C}_1$ independent of $\v\in(0,1]$, such that there exists a unique solution $(\r^\v,u^\v,\t^\v)$ of \eqref{1.1}, \eqref{1.9} with the initial data $(\r_0,u_0,\t_0)$ and satisfies   \eqref{1.19}, \eqref{1.18}.

Then we justify the vanishing dissipation limit as follows:
\begin{theorem}[Vanishing Dissipation Limit]\label{thm1.2}
Based on the above preparations, under the assumptions of Theorem \ref{thm1.1} and $\k(\v)\rightarrow0+$ as $\v\rightarrow0+$,  there exists $T_2=\min\{T_0,T_1\}>0$, which is independent of $\v>0$, such that
{\small\begin{align}
\|(\r^\v-\r,u^\v-u,\t^\v-\t)(t)\|^2_{L^2}&+\int_0^t\v\|(u^\v-u)(\tau)\|^2_{H^1}+\k(\v)\|(\t^\v-\t)(\tau)\|^2_{H^1}d\tau\nonumber\\
&\leq C\max\{\v^{\f32},\k(\v)^{\f32}\},~t\in[0,T_2],\label{1.8-0}\\
\|(\r^\v-\r,u^\v-u,\t^\v-\t)(t)\|^2_{H^1}&+\int_0^t\v\|(u^\v-u)(\tau)\|^2_{H^2}+\k(\v)\|(\t^\v-\t)(\tau)\|^2_{H^2}d\tau\nonumber\\
&\leq C\max\{\v^{\f12},\k(\v)^{\f13}\},~t\in[0,T_2],\label{1.9-0}	
\end{align}}
and
{\small\begin{align}
\|(\r^\v-\r,u^\v-u)\|_{L^\infty(\Omega\times[0,T_2])}&\leq \|(\r^\v-\r,u^\v-u)\|^{\f25}_{L^2}\cdot \|(\r^\v-\r,u^\v-u)\|^{\f35}_{W^{1,\infty}}\nonumber\\
&\leq C\max\{\v,\k(\v)\}^{\f3{10}},\label{1.10-2}
\end{align}}
where $C$ depends only on the norm $\|(\r_0,u_0,\t)\|_{H^3}+\|(\r_0,u_0,\t_0)\|_{X^{\v}_{m}}$.
\end{theorem}

\begin{remark}
It is easy to see that  $k(\v)=\v^a$ with $0< a \leq 4$ satisfy the condition of Theorem \ref{thm1.2}.
\end{remark}
\begin{remark}
Compared to the isentropic case \cite{Wang-Xin-Yong}, one can see that  the convergence rates of vanishing dissipation limit are influenced by the decay rate of heat conductivity. In particular,  for the case $k(\v)=\v$, Theorem \ref{thm1.2} implies that the convergence rate in $L^\infty(0,T_2;H^1)$ is $\v^{\f13}$ which is slower than the isentropic case \cite{Wang-Xin-Yong} whose corresponding rate is $\v^{\f12}$. This is mainly due to the influence of thermal boundary layers, see Lemma \ref{lem14.2} below. If one can prove $\k(\v)^{\f32}\int_0^t\|\nabla\Delta\t^\v\|^2d\tau$ is uniformly bounded, then the convergence rate of \eqref{1.9-0} could be improved to be $\max(\v^{\f12},\k(\v)^{\f12})$, however it is very hard to obtain such uniform estimate in our framework.
\end{remark}
\begin{remark}
By the same arguments as Theorem \ref{thm1.2}, one can also prove the dissipation limit of full compressible Navier-Stokes system to the following system:
{\small\begin{equation}\label{1.7-1}
	\begin{cases}
	\r_t+\mbox{div}(\r u)=0,\\
	(\r u)_t+\mbox{div}(\r u \otimes u)+\nabla p=0,\\
	[\r(\t+\f12|u|^2)]_t+\mbox{div}[\r u(\t+\f12|u|^2)+pu]=\k_0\Delta\t,
	\end{cases}
	\mbox{if}~\k(\v)\rightarrow \k_0>0~\mbox{as}~\v\rightarrow0+.
	\end{equation}}
with boundary conditions
$$u\cdot n|_{\partial\Omega}=0,~n\cdot\nabla\t|_{\partial\Omega}=\nu\t.$$
\end{remark}

The rest of the paper is organized as follows: In section 2, we collect some elementary facts and inequalities that will be used later.   We prove  the a priori estimates Theorem \ref{thm3.1} in section 3 which is the main part of this paper. By using the a priori estimates, we  prove Theorem \ref{thm1.1}  in section 4.  Section 5 is devoted to the proof of  Theorem \ref{thm1.2}.

\

\noindent{\bf Notations:} Throughout this paper, the positive
generic constants that are independent of $\v$ are denoted by
$c,C$(may depend on $\mu,\l$).    $\|\cdot\|$  denotes the standard
$L^2(\Om;dx)$ norm, and $\|\cdot\|_{H^m}~(m=1,2,3,\cdots)$
denotes the Sobolev $H^m(\Om;dx)$ norm. The notation $|\cdot|_{H^m}$ will be used for the standard Sobolev norm of functions defined on $\partial\Omega$. Note that this norm involves only tangential  derivatives.  $P(\cdot)$ denotes a polynomial function.

\section{Preliminaries }
The following lemma \cite{Xiao-Xin-1,Teman} allows one to control the $H^m(\Om)$-norm of a vector valued function $u$ by its $H^{m-1}(\Om)$-norm of $\nabla\times u$ and $\mbox{div} u$, together with the $H^{m-\f12}(\partial\Om)$-norm of $u\cdot n$.

\begin{proposition}\label{prop3.1}
	Let $m\in \mathbb{N}_+$ be an integer. Let $u\in H^m$ be a vector-valued function. Then, there exists a constant $C>0$ independent $u$, such that
{\small	\begin{equation}\label{3.1}
	\|u\|_{H^m}\leq C\left(\|\nabla \times u\|_{H^{m-1}}+\|\mbox{div}u\|_{H^{m-1}}
	+\|u\|_{H^{m-1}}+|u\cdot n|_{H^{m-\f12}(\partial\Om)}\right).
	\end{equation}}
and
{\small\begin{equation}\label{3.1-1}
\|u\|_{H^m}\leq C\left(\|\nabla \times u\|_{H^{m-1}}+\|\mbox{div}u\|_{H^{m-1}}
+\|u\|_{H^{m-1}}+|n\times u|_{H^{m-\f12}(\partial\Om)}\right).
\end{equation}}
\end{proposition}

In this paper, we shall use repeatedly the  Gagliardo-Nirenbirg-Morser type inequality, whose proof can be find in \cite{Gues}. First, define the space
{\small\begin{equation}\label{2.2}
\mathcal{W}^m(\Om\times[0,T])=\{f(x,t)\in L^2(\Omega\times[0,T]) ~|~ \mathcal{Z}^\a f \in L^2(\Omega\times[0,T]) , ~~|\a|\leq m~\}.
\end{equation}}
Then, the Gagliardo-Nirenbirg-Morser type inequality is as follows:
\begin{proposition}\label{prop3.2}
	For $ u,v\in L^\infty(\Om\times[0,T])\cap \mathcal{W}^m(\Om\times[0,T])$ with  $m\in \mathbb{N}_+$ be an integer. It holds that
{\small	\begin{equation}\label{3.2}
	 \int_0^t\|(\mathcal{Z}^{\b}u\mathcal{Z}^{\g}v)(\tau)\|^2d\tau
	\lesssim \|u\|^2_{L^\infty_{t,x}}\int_0^t\|v(\tau)\|^2_{\mathcal{H}^m}d\tau
	 +\|v\|^2_{L^\infty_{t,x}}\int_0^t\|u(\tau)\|^2_{\mathcal{H}^m}d\tau,~~|\b|+|\g|=m.
	\end{equation}}
	
\end{proposition}

We also need the following anisotropic Sobolev embedding and trace estimates:
\begin{proposition}\label{prop3.3}
	Let  $m_1\geq 0,~m_2\geq 0$ be integers, $f\in H^{m_1}_{co}(\Om)\cap  H^{m_2}_{co}(\Om)$  and   $\nabla f\in H^{m_2}_{co}(\Om)$.\\
	1) The following anisotropic Sobolev embedding holds:
{\small	\begin{equation}\label{3.3}
	\|f\|^2_{L^\infty}\leq C \Big(\|\nabla f\|_{H^{m_2}_{co}}+ \|f\|_{H^{m_2}_{co}}\Big)\cdot\|f\|_{H^{m_1}_{co}},
	\end{equation}}
	provided $m_1+m_2\geq 3$.\\[2mm]
	2) The following trace estimate holds:
{\small	\begin{equation}\label{3.4}
	|f|^2_{H^{s}(\partial\Om)}\leq C\Big(\|\nabla f\|_{H^{m_2}_{co}}+ \|f\|_{H^{m_2}_{co}}\Big)\cdot\|f\|_{H^{m_1}_{co}}.
	\end{equation}}
	provided $m_1+m_2\geq 2s\geq 0$.
\end{proposition}

\noindent\textbf{Proof}. The proof is just a using of  the covering $\Om\subset \Om_0\cup_{k=1}^n\Om_k$ and Proposition 2.2 in \cite{Masmoudi-R-1}, the details are   omitted here.  $\hfill\Box$



\section{A priori Estimates}

The aim of this section is to prove the following {\it a priori} estimates, which is a crucial step to prove Theorem \ref{thm1.1}. For notation convenience, we drop the superscript $\v$ throughout this section.

\begin{theorem}\label{thm3.1}
Let $m$ be an integer satisfying $m\geq 6$, $\k(\v)$ satisfies \eqref{1.8},  $\Omega$ be a $\mathcal{C}^{m+2}$ domain and $A\in\mathcal{C}^{m+1}(\partial\Omega)$. For very sufficiently smooth solution defined on $[0,T]$ of \eqref{1.1} and \eqref{1.9} in $\Omega$, then we have
{\small\begin{equation}\label{3.0-3}
|\r(x,0)|\exp(-\int_0^t\|\mbox{div}u\|_{L^\infty}d\tau)\leq\r(x,t)\leq |\r(x,0)|\exp(\int_0^t\|\mbox{div}u\|_{L^\infty}d\tau), ~~\forall t\in[0,T],
\end{equation}}
and
{\small\begin{equation}\label{3.0-5}
\t_0-\int_0^t\|\t(\tau)\|_{L^\i}d\tau\leq\t(x,t)\leq \t_0+\int_0^t\|\t(\tau)\|_{L^\i}d\tau, ~~\forall t\in[0,T].
\end{equation}}
In addition, if it holds that
{\small\begin{equation}\label{3.0-4}
0<c_0\leq \r(x,t), \t(x,t)\leq \f1{c_0}<\infty,~~\forall t\in[0,T],
\end{equation}}
where $c_0$ is any given small positive constant, then we have the a priori estimate
{\small\begin{align}\label{3.0-1}
&\Upsilon_m(\r,u,\t)\triangleq\mathcal{N}_m(t)
+\int_0^t\|\nabla\partial_t^{m-1}(\r,~\t)(\tau)\|^2d\tau+\int_0^t\|\nabla (\r\t)(\tau)\|^2_{\MH^{2,\infty}}d\tau+\int_0^t\v\|\nabla u(\tau)\|^2_{\mathcal{H}^{m}}d\tau\nonumber\\
&~~~~~~~~~~~~~~~~+\k(\v)\int_0^t\|\nabla \t(\tau)\|^2_{\mathcal{H}^{m}}d\tau+\sum_{k=0}^{m-2}\int_0^t\v\|\nabla^2\partial_t^{k}u(\tau)\|^2_{m-1-k}+\k(\v)\|\Delta\partial_t^{k}\t(\tau)\|^2_{m-1-k}d\tau
\nonumber\\
&~~~~~~~~~~~~~~~~+\v^2\int_0^t\|\nabla^2\partial_t^{m-1}u(\tau)\|^2d\tau+\k(\v)^2\int_0^t\|\Delta\partial_t^{m-1}\t(\tau)\|^2+\|\nabla\MZ^{m-2}\Delta\t(\tau)\|^2d\tau\nonumber\\
&~~~~~~~~~~~~~~\leq \tilde{C}_2C_{m+2}\Big\{P(\mathcal{N}_{m}(0)) +  tP(\mathcal{N}_{m}(t))\Big\}.
\end{align}}
where $\tilde{C}_2$ depends only on $\f1{c_0}$, $P(\cdot)$ is a polynomial and
{\small\begin{align}\label{3.0-2}
\mathcal{N}_m(t)\triangleq\mathcal{N}_m(\r,u,\t)(t)=\sup_{0\leq\tau\leq t}\Big\{1+
\|(\r,u,\t)(\tau)\|^2_{X^\v_m}\Big\}.
\end{align}}
\end{theorem}

Throughout this section, we shall work on the interval of time $[0,T]$ such that $c_0\leq \r, \t\leq \f1{c_0}$. And we point out that the generic constant $C$ may depend on $\f1{c_0}$ in this section.  Since the proof of Theorem \ref{thm3.1} is very complicated, we shall divide the proof into the following several subsections.

\subsection{Conormal Energy Estimates for $\r,u$ and $\t$}

\renewcommand{\theequation}{\arabic{section}.\arabic{subsection}.\arabic{equation}}

Notice that
{\small\begin{equation}\label{2.4}
\Delta u=\nabla\mbox{div} u-\nabla\times\nabla\times u,
\end{equation}}
then $\eqref{1.1}$ is rewritten as
{\small\begin{eqnarray}\label{2.5}
\begin{cases}
\r_t+\mbox{div}(\r u)=0,\\
\r u_t+\r u\cdot\nabla u+\nabla p=-\mu\v \nabla\times\omega+(2\mu+\l)\v \nabla div u,\\
\r \t_t+\r u\cdot\nabla\t+p \mbox{div}u=\k(\v)\Delta\t+2\mu\v|Su|^2+\l\v|\mbox{div}u|^2,
\end{cases}
\end{eqnarray}}
where $\omega=\nabla \times u$ is the vorticity.

\begin{lemma}\label{lem2.2}
For a smooth solution of \eqref{1.1} and \eqref{1.9}, it holds that for $\v\in(0,1]$
{\small\begin{align}\label{2.13}
&\sup_{0\leq\tau\leq t}\Big(\int R f_1(\r)+\r f_2(\t)+ \f12\r|u|^2  dx\Big)+c_1\int_0^t\v\|\nabla u(\tau)\|^2+\k(\v)\|\nabla \t(\tau)\|^2d\tau\nonumber\\
&\leq \int R f(\r_0)+\r_0 h(\t_0)+ \f12\r_0|u_0|^2  dx
+ C \int_0^t\|u(\tau)\|^2d\tau+Ct,
\end{align}}
where $c_1>0$ and $f_1(t)=t\ln t-t+1$ and $f_2(t)=t-\ln t+1$ for $t>0$.
\end{lemma}
\noindent\textbf{Proof}.  Multiplying \eqref{2.5} by $\f1\t$, one gets that
{\small\begin{align}\label{2.14}
&\f{d}{dt}\int \r\ln\t dx+\int R\r\mbox{div}u dx\nonumber\\
&=\k(\v)\int\f{|\nabla\t|^2}{\t^2}dx
+\v\int\f1\t(2\mu|Su|^2+\l|\mbox{div}u|^2)dx+\k(\v)\int_{\partial\Omega}\f{\nabla\t\cdot n}{\t}d\sigma.
\end{align}}
It follows from the boundary condition \eqref{1.9} that
{\small\begin{eqnarray}\label{2.15}
\k(\v)\int_{\partial\Omega}\f{\nabla\t\cdot n}{\t}d\sigma=
\nu\k(\v)\int_{\partial\Omega} d\sigma\leq C\k(\v)\leq C.
\end{eqnarray}}
We rewrite the mass equation $\eqref{2.5}_1$ to be
{\small\begin{equation}\label{2.16}
\r_t+u\cdot\nabla\r+\r \mbox{div}u=0,
\end{equation}}
which yields immediately  that
{\small\begin{align}\label{2.17}
\int\r\mbox{div}u dx=-\int\r_t+u\cdot\nabla\r dx=-\int\r((\ln\r)_t+u\cdot\nabla\ln\r) dx
=-\f{d}{dt}\int\r\ln\r dx.
\end{align}}
Substituting \eqref{2.17} and \eqref{2.15} into \eqref{2.14}, one obtains that
{\small\begin{align}\label{2.18}
\f{d}{dt}\int (R\r\ln\r-\r\ln\t) dx+\int\k(\v)\frac{|\nabla\t|^2}{\t^2}+\f\v\t(2\mu|Su|^2+\l|\mbox{div}u|^2)dx\leq C.
\end{align}}
On the other hand, it follows from $\eqref{1.1}_1$ and $\eqref{1.1}_3$  that
{\small\begin{align}\label{2.17-1}
\f{d}{dt}\int \r dx=0,
\end{align}}
and
{\small\begin{align}\label{2.17-2}
\f{d}{dt}\int \r\t+\f12\r|u|^2 dx=\k(\v)\int_{\partial\Omega} n\cdot\nabla\t d\sigma+2\mu\v\int_{\partial\Omega} ((Su)u) n d\sigma\leq C+\d\v\|\nabla u\|^2+C_\d\|u\|^2,
\end{align}}
where we have used the following facts in the estimates of \eqref{2.17-2}
{\small\begin{align}\nonumber
&\k(\v)\Big|\int_{\partial\Omega} n\cdot\nabla\t d\sigma\Big| +2\mu\v\Big|\int_{\partial\Omega} ((Su)u) n d\sigma\Big|=\nu\k(\v)\Big|\int_{\partial\Omega} \t d\sigma\Big|+2\mu\v\Big|\int_{\partial\Omega} ((Su)n)_\tau u_\tau d\sigma\Big|\nonumber\\
&\leq C+C\v|u|^2_{L^2}\leq C+\d\v\|\nabla u\|^2+C_\d\|u\|^2.\nonumber
\end{align}}
Then, combining \eqref{2.18}, \eqref{2.17-1} and \eqref{2.17-2}, one  obtains that
{\small\begin{align}\label{2.19}
&\f{d}{dt}\Big(\int R(\r\ln\r-\r+1)+\r(\t-\ln\t+1)+ \f12\r|u|^2  dx\Big)+\int\k(\v)\frac{|\nabla\t|^2}{\t^2}+\f\v\t(2\mu|Su|^2+\l|\mbox{div}u|^2)dx\nonumber\\
&\leq  \d\v^2\|\nabla u\|^2
+ C_{\d}\|u\|^2+C.
\end{align}}
It follows from the Korn's inequality and the fact $2\mu+\l>0$  that
{\small\begin{eqnarray}\label{2.20}
&&\int\k(\v)\frac{|\nabla\t|^2}{\t^2}+\f\v\t(2\mu|Su|^2+\l|\mbox{div}u|^2)dx\geq  2c_1\Big[\k(\v)\|\nabla\t\|^2+\v\|\nabla u\|^2\Big]-C \|u\|^2.
\end{eqnarray}}
where $c_1>0$ is a given positive constant depends on $c_0,\mu,\l,\k$. Thus, choosing $\d$ small, it holds that
{\small\begin{equation}\label{2.21}
\f{d}{dt}\Big(\int R(\r\ln\r-\r+1)+\r(\t-\ln\t+1)+ \f12\r|u|^2  dx\Big)+c_1\k(\v)\|\nabla\t\|^2
+c_1\v\|\nabla u\|^2\leq C+ C \|u\|^2.
\end{equation}}
Integrating \eqref{2.21} over $[0,t]$, one obtains \eqref{2.13}. Thus the proof of the Lemma \ref{lem2.2} is completed.
$\hfill\Box$

\

However, the above basic energy estimates is far from enough to get the vanishing dissipation limit. One needs to get some conormal derivative estimates. Set
{\small\begin{equation}\label{2.15-2}
Q(t)\triangleq \sup_{0\leq\tau\leq t}\Big\{\|\nabla( \r, u,\t)(t)\|^2_{\mathcal{H}^{1,\infty}}
+\|(\t,u,\t, \r_t,u_t,\t_t)(t)\|^2_{L_x^\infty}+\v\|\nabla^2u\|^2_{L^\infty}\Big\}.
\end{equation}}
It follows from  Proposition \ref{prop3.3}  that
{\small\begin{align}\label{2.17-3}
Q(t)\leq C P(\mathcal{N}_m(t))~~\mbox{for}~~m\geq3.
\end{align}}

\

\begin{lemma}\label{lem2.3}
For $m\geq3$, it holds that
{\small\begin{align}\label{2.16-1}
&\sup_{0\leq \tau\leq t}\|(\r, \t,  u)(\tau)\|^2_{\mathcal{H}^m}+\v\int_0^t\|\nabla u(\tau)\|^2_{\mathcal{H}^m}d\tau+\k(\v)\int_0^t\|\nabla\t(\tau)\|^2_{\mathcal{H}^m}d\tau\nonumber\\
&\leq
CC_{m+2} \Big\{1+\|(\r_0,u_0,\t_0)\|^2_{\mathcal{H}^m}+\d \v^2\int_0^t\|\nabla^2  u(\tau)\|^2_{\mathcal{H}^{m-1}}d\tau+\d\k(\v)^2\int_0^t\|\Delta \t(\tau)\|^2_{\mathcal{H}^{m-1}}d\tau\nonumber\\
&~~~~~~~~~~~~~~~~~~~~+\d\int_0^t\|\nabla\partial_t^{m-1}(\r,\t)(\tau)\|^2d\tau+C_\d t P(\mathcal{N}_m(t))\Big\},
\end{align}}
where $\d$ is small which will be chosen later and  $C_\d$ is a polynomial function of $\f1\d$ which may vary from line to line.
\end{lemma}
\noindent\textbf{Proof}.
The estimate for $k=0$ is already given in Lemma \ref{lem2.2}. Assuming that it is proven for $k\leq m-1$. We shall prove it for $ k=m\geq 1$. Applying $\mathcal{Z}^\a$ with $|\a|=m$ to \eqref{2.5}, one obtains that
{\small\begin{align}\label{2.26}
\begin{cases}
\r\mathcal{Z}^\a u_t+\r u\cdot\nabla\mathcal{Z}^\a u+\mathcal{Z}^\a\nabla p=-\mu\v\mathcal{Z}^\a\nabla\times\omega+(2\mu+\l)\v\mathcal{Z}^\a\nabla\mbox{div}u
+\mathcal{C}_1^{\a}+\mathcal{C}_2^{\a},\\[2mm]
\r\mathcal{Z}^\a \t_t+\r u\cdot\nabla\mathcal{Z}^\a \t+p\mathcal{Z}^\a\mbox{div}u
-\k(\v)\MZ^\a\Delta\t\\
~~~~~~~~~~~~~~~~~~~~~~~~~~~~~~~~~~~~~~ =2\mu\v\mathcal{Z}^\a(|Su|^2)+\l\v\MZ^\a(|\mbox{div}u|^2)
+\mathcal{C}_3^{\a}+\mathcal{C}_4^{\a}+\mathcal{C}_5^{\a},
\end{cases}
\end{align}}
with
{\small\begin{align}\label{2.30}
\begin{cases}
\mathcal{C}_1^{\a}=-[\mathcal{Z}^\a, \r]u_t=-\sum_{|\b|\geq 1,\b+\g=\a}C_{\a,\b}\mathcal{Z}^\b\r\mathcal{Z}^\g u_t,\\
\mathcal{C}_2^{\a}=-[\mathcal{Z}^\a, \r u\cdot\nabla]u=-\sum_{|\b|\geq 1,\b+\g=\a}C_{\a,\b}\mathcal{Z}^\b(\r u)\mathcal{Z}^\g\nabla u-\r u\cdot [\mathcal{Z}^\a,\nabla]u,\\
\mathcal{C}_3^{\a}=-[\mathcal{Z}^\a, \r]\t_t=-\sum_{|\b|\geq 1,\b+\g=\a}C_{\a,\b}\mathcal{Z}^\b\r\mathcal{Z}^\g \t_t,\\
\mathcal{C}_4^{\a}=-[\mathcal{Z}^\a, \r u\cdot\nabla]\t=-\sum_{|\b|\geq 1,\b+\g=\a}C_{\a,\b}\mathcal{Z}^\b(\r u)\mathcal{Z}^\g\nabla \t-\r u\cdot [\mathcal{Z}^\a,\nabla]\t,\\
\mathcal{C}_5^{\a}=-[\mathcal{Z}^\a, p]\mbox{div}u=-\sum_{|\b|\geq 1,\b+\g=\a}C_{\a,\b}\mathcal{Z}^\b p\mathcal{Z}^\g\mbox{div}u,
\end{cases}
\end{align}}
where $C_{\a,\b}$ are the corresponding  binomial coefficients. Multiplying \eqref{2.26} by $\mathcal{Z}^\a u$ and  integrating by parts, one has  that
{\small\begin{eqnarray}\label{2.32}
&&\f{d}{dt}\int\f12\r|\mathcal{Z}^\a u|^2dx+\int\mathcal{Z}^\a\nabla p\mathcal{Z}^\a udx=-\mu\v\int\mathcal{Z}^\a\nabla\times\omega\cdot\mathcal{Z}^\a udx\nonumber\\
&&~~~~~~+(2\mu+\l)\v\int\mathcal{Z}^\a\nabla\mbox{div}u\cdot\mathcal{Z}^\a udx
+\int(\mathcal{C}_1^{\a}+\mathcal{C}_2^{\a})\mathcal{Z}^\a udx
\end{eqnarray}}
Using the same arguments as Lemma 3.3 of \cite{Wang-Xin-Yong}, one can get that
{\small\begin{align}\label{2.36}
&-\v\int\mathcal{Z}^\a\nabla\times\omega\cdot\mathcal{Z}^\a udx
\leq -\f{3\v}{4}\|\nabla\times\mathcal{Z}^\a u\|^2+\d\v \|\nabla u\|^2_{\mathcal{H}^{m}}+\d\v^2\|\nabla^2u\|^2_{\mathcal{H}^{m-1}}\nonumber\\
&~~~~~~~~~~~~~~~~~~~~~~~~~~~~~~~~~~~~~+C_\d C_{m+2} (\|\nabla u\|^2_{\mathcal{H}^{m-1}}+\|u\|^2_{\mathcal{H}^{m}}).
\end{align}}
and
{\small\begin{align}\label{2.40}
&\v\int\mathcal{Z}^\a\nabla\mbox{div}u\cdot\mathcal{Z}^\a udx
\leq -\f{3\v}{4}\|\mbox{div}\mathcal{Z}^\a u\|^2
+\d\v\|\nabla u\|^2_{\mathcal{H}^{m}}+\d\v^2\|\nabla^2u\|^2_{\mathcal{H}^{m-1}}\nonumber\\
&~~~~~~~~~~~~~~~~~~~~~~~~~~~~~~~~~~+C_\d C_{m+2} (\|\nabla u\|^2_{\mathcal{H}^{m-1}}
+\|u\|^2_{\mathcal{H}^{m}}).
\end{align}}
It follows from  Proposition \ref{prop3.1} that
{\small\begin{align}\label{2.41}
 2c_1\|\nabla\mathcal{Z}^\a u\|_{L^2}
&\leq  \Big(\|\nabla\times\mathcal{Z}^\a u\|_{L^2}
+\|\mbox{div}\mathcal{Z}^\a u\|_{L^2}+\|\mathcal{Z}^\a u\|_{L^2}+|\mathcal{Z}^\a u\cdot n|_{H^{\f12}(\partial\Om)}\Big)\nonumber\\
&\leq \Big(\|\nabla\times\mathcal{Z}^\a u\|_{L^2}
+\|\mbox{div}\mathcal{Z}^\a u\|_{L^2}\Big)+C_{m+2} \Big(\|u\|^2_{\mathcal{H}^{m}}+\|\nabla u\|^2_{\mathcal{H}^{m-1}}\Big).
\end{align}}
Substituting \eqref{2.36}-\eqref{2.40}  into \eqref{2.32} and using \eqref{2.41}, then integrating the resultant inequality over $[0,t]$, one obtains that
{\small\begin{eqnarray}\label{2.42}
&&\f12\int\r|\mathcal{Z}^\a u|^2dx+\int_0^t\int\mathcal{Z}^\a\nabla p\cdot\mathcal{Z}^\a udxd\tau
+2c_1\v\int_0^t\|\nabla\mathcal{Z}^\a u(\tau)\|^2_{L^2}d\tau\nonumber\\
&&\leq \f12\int\r_0|\mathcal{Z}^\a u_0|^2dx+C\d\v^2\int_0^t\|\nabla^2u(\tau)\|^2_{\mathcal{H}^{m-1}}d\tau
+C\d\v\int_0^t\|\nabla u(\tau)\|^2_{\mathcal{H}^{m}}d\tau\nonumber\\
&&~~~~+\int_0^t\int(\mathcal{C}_1^{\a}+\mathcal{C}_2^{\a})\cdot\mathcal{Z}^\a udxd\tau+C_{m+2}C_\d \int_0^t\|\nabla u(\tau)\|^2_{\mathcal{H}^{m-1}}+\|u(\tau)\|^2_{\mathcal{H}^{m}}d\tau.
\end{eqnarray}}

On the other hand, multiplying \eqref{2.26} by $\frac{\mathcal{Z}^\a \t}{\t}$ and  integrating by parts, one has  that
{\small\begin{align}\label{2.43}
&\f{d}{dt}\int\f{\r}{2\t}|\mathcal{Z}^\a \t|^2dx+\int R\r\mathcal{Z}^\a\mbox{div}u \mathcal{Z}^\a \t dx-\k(\v)\int\mathcal{Z}^\a\Delta\t\cdot\frac{\mathcal{Z}^\a \t}{\t}dx\nonumber\\
&\leq \v\int\Big(2\mu\mathcal{Z}^\a(|Su|^2)+\l\mathcal{Z}^\a(|\mbox{div}u|^2)\Big)\frac{\mathcal{Z}^\a \t}{\t} dx
+\int(\mathcal{C}_3^{\a}+\mathcal{C}_4^{\a}+\mathcal{C}_5^{\a})\frac{\mathcal{Z}^\a \t}{\t}dx+CP(\mathcal{N}_m(t)).
\end{align}}
It follows from the boundary condition  $\eqref{1.9}_3$ that
{\small\begin{align}\label{2.44}
&|n\cdot\MZ^\a\nabla\t|_{L^2}\leq C_{m+1}\Big(|\MZ^\a\t|_{L^2}+|\MZ^{m-1}\nabla\t|_{L^2}\Big)\nonumber\\
&\leq C_{m+1}\Big(\|\MZ^m\t\|^{\f12}\|\nabla\MZ^m\t\|^{\f12}+\|\MZ^m\t\|
+\|\nabla^2\t\|^{\f12}_{\mathcal{H}^{m-1}}\|\nabla\t\|^{\f12}_{\mathcal{H}^{m-1}}+\|\nabla\t\|^{\f12}_{\mathcal{H}^{m-1}}\Big),
\end{align}}
which, together with integrating by parts, yields that
{\small\begin{align}\label{2.45}
&-\k(\v)\int\mathcal{Z}^\a\Delta\t\cdot\frac{\mathcal{Z}^\a \t}{\t}dx
=-\k(\v)\int\MZ^\a\mbox{div}\nabla\t\frac{\mathcal{Z}^\a \t}{\t}dx\nonumber\\
&=-\k(\v)\int\mbox{div}\MZ^\a\nabla\t\frac{\mathcal{Z}^\a \t}{\t}dx
-\k(\v)\int[\mbox{div},\MZ^\a]\nabla\t \frac{\mathcal{Z}^\a \t}{\t}dx\nonumber\\
&=\k(\v)\int\MZ^\a\nabla\t\frac{\nabla\mathcal{Z}^\a \t}{\t}dx
-\k(\v)\int[\mbox{div},\MZ^\a]\nabla\t \frac{\mathcal{Z}^\a \t}{\t}dx
-\k(\v)\int_{\partial\Omega}\frac{\mathcal{Z}^\a \t}{\t} \MZ^\a\nabla\t\cdot n d\sigma\\
&\geq \k(\v)\int\f{|\nabla\mathcal{Z}^\a \t|^2}{\t}dx-C\v^a\Big(\|\nabla\t\|_{\mathcal{H}^{m-1}}\|\nabla\MZ^\a\t\|
+\|\t\|_{\mathcal{H}^{m}}\|\nabla^2\t\|_{\mathcal{H}^{m-1}}+\|\MZ^\a\t\|_{L^2}^\f12\|\nabla\MZ^\a\t\|^\f12|\MZ^\a\nabla\t\cdot n|_{L^2}\Big)\nonumber\\
&\geq \f{3\k(\v)}{4}\int\f{|\nabla\mathcal{Z}^\a \t|^2}{\t}dx-\d\k(\v)^2\|\Delta\t\|^2_{\mathcal{H}^{m-1}}-\d\k(\v)\|\nabla\MZ^m\t\|^2-\d\|\nabla\partial_t^{m-1}\t\|^2-C_\d C_{m+1}P(\mathcal{N}_m(t)),\nonumber
\end{align}}
where we have used the Young's inequality in the last inequality of \eqref{2.45}.
It is easy to calculate that
{\small\begin{equation}\label{2.46}
\v\int_0^t\int \frac{\mathcal{Z}^\a \t}{\t}\Big[2\mu\MZ^\a(|Su|^2)+\l\MZ^\a(|\mbox{div}u|^2)\Big]dxd\tau
\leq \d\v^2\int_0^t\|\nabla u\|^2_{\mathcal{H}^{m}}d\tau+C_\d tP(\mathcal{N}_m(t)).
\end{equation}}
Substituting \eqref{2.46} and \eqref{2.45} into \eqref{2.43} and integrating the resultant inequality over $[0,t]$, one gets that
{\small\begin{align}\label{2.47}
&\int\f{\r}{2\t}|\mathcal{Z}^\a \t|^2dx+\int_0^t\int R\r\mathcal{Z}^\a\mbox{div}u \mathcal{Z}^\a \t dxd\tau+\f{3\k(\v)}{4}\int_0^t\int\f{|\nabla\mathcal{Z}^\a \t|^2}{\t}dxd\tau\nonumber\\
&\leq \int\f{\r_0}{2\t_0}|\mathcal{Z}^\a \t_0|^2dx+C\d\int_0^t\k(\v)\|\nabla\MZ^m\t(\tau)\|^2+\k(\v)^2\|\Delta\t(\tau)\|^2_{\mathcal{H}^{m-1}}d\tau+\d\v\int_0^t\|\nabla u(\tau)\|^2_{\mathcal{H}^m}d\tau\nonumber\\
&~~+\d\int_0^t\|\nabla\partial_t^{m-1}\t(\tau)\|^2d\tau+C_\d t P(\mathcal{N}_m(t))+\int_0^t\int(\mathcal{C}_3^{\a}+\mathcal{C}_4^{\a}+\mathcal{C}_5^{\a})\frac{\mathcal{Z}^\a \t}{\t}dxd\tau.
\end{align}}
Combining \eqref{2.42} and \eqref{2.47}, one obtains that
{\small\begin{align}\label{2.48}
&\f12\int\r|\mathcal{Z}^\a u|^2+\f{\r}{2\t}|\mathcal{Z}^\a \t|^2dx
+2c_1\v\int_0^t\|\nabla\mathcal{Z}^\a u(\tau)\|^2_{L^2}d\tau+\f{3\k(\v)}{4}\int_0^t\int\f{|\nabla\mathcal{Z}^\a \t|^2}{\t}dxd\tau\nonumber\\
&+R\int_0^t\int\mathcal{Z}^\a\nabla(\r\t)\cdot\mathcal{Z}^\a u+\r\mathcal{Z}^\a\mbox{div}u\cdot\mathcal{Z}^\a \t dxd\tau \nonumber\\
&\leq \f12\int\r_0|\mathcal{Z}^\a u_0|^2+\f{\r_0}{2\t_0}|\mathcal{Z}^\a \t_0|^2dx+C\d\int_0^t\v\|\nabla u(\tau)\|^2_{\mathcal{H}^{m}}+\v^2\|\nabla^2u(\tau)\|^2_{\mathcal{H}^{m-1}}\nonumber\\
&~+C\d\int_0^t
\k(\v)\|\nabla\t(\tau)\|^2_{\mathcal{H}^{m}}+\k(\v)^2\|\Delta\t(\tau)\|^2_{\mathcal{H}^{m-1}}d\tau
+C\d\int_0^t\|\nabla\partial_t^{m-1}\t(\tau)\|^2d\tau\nonumber\\
&~~~~+C_\d C_{m+2}t P(\mathcal{N}_m(t))+C\Big(\int_0^t\|(\mathcal{C}_1^{\a},\mathcal{C}_2^{\a},\mathcal{C}_3^{\a},\mathcal{C}_4^{\a},\mathcal{C}_5^{\a})\|^2d\tau\Big)^{\f12}\Big(\int_0^t\|\mathcal{Z}^\a(u,\t)\|^2d\tau\Big)^{\f12}.
\end{align}}
Now we estimate the fourth term on the LHS of \eqref{2.48}. Note that
{\small\begin{align}\label{2.49}
&\MZ^\a\nabla(\r\t)=\t\cdot\nabla\MZ^\a\r+\r\cdot\nabla\MZ^\a\t+[\MZ^\a,\nabla\r]\t+[\MZ^\a,\nabla\t]\r+\t\cdot[\MZ^\a,\nabla]\r+\r\cdot[\MZ^\a,\nabla]\t,
\end{align}}
which, together with Proposition \ref{prop3.2} and Holder inequality, yields that
{\small\begin{align}\label{2.50}
& I=\int_0^t\int\mathcal{Z}^\a\nabla(\r\t)\cdot\mathcal{Z}^\a u+\r\mathcal{Z}^\a\mbox{div}u\cdot\mathcal{Z}^\a \t dxd\tau \nonumber\\
&\geq \int_0^t\int(\t\cdot\nabla\MZ^\a\r+\r\cdot\nabla\MZ^\a\t)\mathcal{Z}^\a u dxd\tau+\int_0^t\int\r\mathcal{Z}^\a\mbox{div}u\cdot\mathcal{Z}^\a \t dxd\tau \nonumber\\
&~~~~-C\d \int_0^t\|\nabla\partial_t^{m-1}(\r,\t)\|^2d\tau-C_\d (1+Q(t))\int_0^t\mathcal{N}_m(\tau)d\tau\nonumber\\
&\geq -\int_0^t\int(\t\MZ^\a\r+\r\MZ^\a\t)\mbox{div}\mathcal{Z}^\a u dxd\tau+\int_0^t\int\r\mathcal{Z}^\a\mbox{div}u\cdot\mathcal{Z}^\a \t dxd\tau \nonumber\\
&~~~~-C\d \int_0^t\|\nabla\partial_t^{m-1}(\r,\t)\|^2d\tau+\int_0^t\int_{\partial\Omega}(\t\MZ^\a\r+\r\MZ^\a\t)\mathcal{Z}^\a u\cdot n d\sigma d\tau-C_\d P(\mathcal{N}_m(t)) \nonumber\\
&\geq -\int_0^t\int \t\MZ^\a\r\cdot\mbox{div}\mathcal{Z}^\a u dxd\tau
+\int_0^t\int_{\partial\Omega}(\t\MZ^\a\r+\r\MZ^\a\t)\mathcal{Z}^\a u\cdot n d\sigma d\tau\nonumber\\
&~~~~~~~~~~~~~~~~~~~~~-C\d \int_0^t\|\nabla\partial_t^{m-1}(\r,\t)\|^2d\tau-C_\d t P(\mathcal{N}_m(t)).
\end{align}}
We shall calculate the  boundary term in \eqref{2.50} when $\a_{13}=0$(for $\a_{13}\neq0$, we have that $\mathcal{Z}^\a u=0$ on the boundary) in the right hand side of \eqref{2.50}.
It follows from that \eqref{1.9} and \eqref{3.4}, for $k\leq m$, that
{\small\begin{eqnarray}\label{2.51}
|Z_y^{m-k}\partial_t^k u\cdot n|_{H^{\f12}}\leq
\begin{cases}
0, ~\mbox{if}~~k=m,\\
 C_{m+2}\{\|\nabla u\|_{\mathcal{H}^{m-1}}+\|u\|_{\mathcal{H}^{m}}\}, ~\mbox{if}~~k\leq m-1.
\end{cases}
\end{eqnarray}}
If $|\a_0|=|\a|$,   it follows from \eqref{2.51} that
{\small\begin{equation}\label{2.52}
\int_0^t\int_{\partial\Omega}(\t\MZ^\a\r+\r\MZ^\a\t)\mathcal{Z}^\a u\cdot n d\sigma d\tau=0.
\end{equation}}
If $|\a_1|\geq 1$, integrating by parts along the boundary and using \eqref{2.52}, one has that
{\small\begin{align}\label{2.53}
&|\int_0^t\int_{\partial\Omega}(\t\MZ^\a\r+\r\MZ^\a\t)\mathcal{Z}^\a u\cdot n d\sigma d\tau|=|\int_0^t\int_{\partial\Om}(\t Z_y^{\a_1}\partial_t^{\a_0}\r +
\r Z_y^{\a_1}\partial_t^{\a_0}\t) Z^{\a}u\cdot n d\sigma d\tau|\nonumber\\
&\leq C Q(t) \int_0^t\Big(|Z_y^{\a_1-1}\partial_t^{\a_0}\r|_{H^{\f12}} +
|Z_y^{\a_1-1}\partial_t^{\a_0}\t|_{H^{\f12}}\Big) |Z^{\a}u\cdot n|_{H^{\f12}} d\sigma d\tau\nonumber\\
&\leq \d\int_0^t \|\nabla(\r,\t)\|^2_{\mathcal{H}^{m-1}}d\tau+C_\d C_{m+2}t P(\mathcal{N}_m(t)).
\end{align}}
Therefore from \eqref{2.53} and \eqref{2.52}, we obtain that
{\small\begin{equation}\label{2.54}
|\int_0^t\int_{\partial\Omega}(\t\MZ^\a\r+\r\MZ^\a\t)\mathcal{Z}^\a u\cdot n d\sigma d\tau|
\leq \d \int_0^t\|\nabla\partial_t^{m-1}(\r,\t)\|^2d\tau+C_\d C_{m+2}t P(\mathcal{N}_m(t)).
\end{equation}}
In order to estimate the first term on the RHS of \eqref{2.50}, we use the following equation which is derived from $\eqref{2.5}_1$, i.e,
{\small\begin{equation}\label{2.55}
\mbox{div}u=-\f{\r_t}{\r}-\f{u}{\r}\cdot\nabla \r.
\end{equation}}
Applying $\MZ^\a$ to \eqref{2.55}, one immediately obtains that
{\small\begin{equation}\label{2.56}
\mathcal{Z}^\a\mbox{div}u=-\f{1}{\r}\mathcal{Z}^{\a}\r_t-\f{u}{\r}\cdot\mathcal{Z}^{\a}\nabla\r
-\sum_{|\beta|\geq1, \b+\g=\a}C_{\a,\b}\MZ^\b(\f1{\r})\cdot\MZ^{\g}\r_t
-\sum_{|\beta|\geq1, \b+\g=\a}C_{\a,\b}\MZ^\b(\f{u}{\r})\cdot\MZ^{\g}\nabla\r.
\end{equation}}
It is easy to get that
{\small\begin{equation}\label{2.57}
\int_0^t\int\f{\t}{\r}\MZ^{\a}\r\cdot \mathcal{Z}^\a\r_t dxd\tau
\geq \int\f{\t}{2\r}|\MZ^{\a}\r|^2dx-\int\f{\t_0}{2\r_0}|\MZ^{\a}\r_0|^2dx
-C t P(\mathcal{N}_m(t)).
\end{equation}}
Integrating by parts and using boundary condition \eqref{1.9}, one has
{\small\begin{align}\label{2.58}
&\int_0^t\int\f{\t}{\r}u\mathcal{Z}^\a\r \cdot\mathcal{Z}^{\a}\nabla\r dxd\tau
=\int_0^t\int\f{\t}{\r}u\mathcal{Z}^\a\r \Big(\nabla\mathcal{Z}^{\a}\r+[\mathcal{Z}^{\a},\nabla]\r\Big) dxd\tau\nonumber\\
&\geq -\d \int_0^t\|\nabla\partial_t^{m-1}\r(\tau)\|^2d\tau-  C_\d   t P(\mathcal{N}_m(t)).
\end{align}}
It follows from Proposition \ref{prop3.2}  that
{\small\begin{align}\label{2.59}
&\sum_{|\beta|\geq1, \b+\g=\a}\Big|\int_0^t\int C_{\a,\b}\t\mathcal{Z}^\a\r\cdot\Big(\MZ^\b(\f1{ \r})\cdot\MZ^{\g}\r_t+\MZ^\b(\f{u}{\r})\MZ^\g\nabla\r\Big)dxd\tau\Big|\nonumber\\
&\leq \d \int_0^t\|\nabla\partial_t^{m-1}\r(\tau)\|^2d\tau+ C_\d  t P(\mathcal{N}_m(t)).
\end{align}}
Combining  \eqref{2.56}-\eqref{2.59}, one obtains  that
{\small\begin{align}\label{2.60}
&-\int_0^t\int\t\mathcal{Z}^\a\r\cdot\mathcal{Z}^\a \mbox{div}udxd\tau\nonumber\\
&\geq \int\f{1}{2\r}|\MZ^{\a}\r|^2dx-\int\f{1}{2\r_0}|\MZ^{\a}\r_0|^2dx-C\d \int_0^t\|\nabla\partial_t^{m-1}\r(\tau)\|^2d\tau
-C_\d  t P(\mathcal{N}_m(t)).
\end{align}}
Substituting \eqref{2.60} and \eqref{2.54} into \eqref{2.50}, one gets  that
{\small\begin{align}\label{2.61}
& I=\int_0^t\int\mathcal{Z}^\a\nabla(\r\t)\cdot\mathcal{Z}^\a u+\r\mathcal{Z}^\a\mbox{div}u\cdot\mathcal{Z}^\a \t dxd\tau \nonumber\\
&\geq \int\f{1}{2\r}|\MZ^{\a}\r|^2dx-\int\f{1}{2\r_0}|\MZ^{\a}\r_0|^2dx
-C\d \int_0^t\|\nabla\partial_t^{m-1}(\r,\t)\|^2d\tau-C_{m+2}C_\d  t P(\mathcal{N}_m(t))
\end{align}}
In order to complete the estimates in \eqref{2.48}, it remains to estimate the terms involving $\mathcal{C}_i^\a,~i=1\cdots 5$.  It follows from  Proposition \ref{prop3.2} and \eqref{2.41}  that
{\small\begin{eqnarray}\label{2.62}
&&\sum_{i=1}^5\int_0^t\|\mathcal{C}_i^\a\|^2 dxd\tau
\leq  C[1+P(Q(t))]\int_0^t\|\nabla\partial_t^{m-1}\t\|^2d\tau+Ct P(\mathcal{N}_m(t)),
\end{eqnarray}}
Therefore, substituting \eqref{2.62} and \eqref{2.61} into \eqref{2.48} and using Cauchy inequality, one proves \eqref{2.16}. Thus,  the proof of the Lemma \ref{lem2.3} is completed.
$\hfill\Box$

\subsection{ Conormal Estimates for $\mbox{div}u$, $\nabla{\r}$ and $\nabla\t$}
\setcounter{equation}{0}

In order to use the compactness argument in the proof of the vanishing dissipation limit, one needs some uniform spatial derivative estimates. In this subsection, we shall get some uniform estimates on $\mbox{div}u$, $\nabla{\r}$ and $\nabla\t$.   In fact, in order to get the uniform estimate of $\|\nabla u\|_{\mathcal{H}^{m-1}}$, one needs  the uniform estimate of $\|\mbox{div} u\|_{\mathcal{H}^{m-1}}$ since we consider the compressible flow in this paper.
\begin{lemma}\label{lem4.1}
For  $m\geq3$, it holds that
{\small\begin{align}\label{4.1}
&\sup_{0\leq \tau\leq t}\|(\mbox{div}u,\nabla\r,\nabla\t)(\tau)\|^2+\v\int_0^t\|\nabla \mbox{div}u(\tau)\|^2d\tau+
\k(\v)\int_0^t\|\Delta\t(\tau)\|^2d\tau\nonumber\\
&\leq C\Big\{\|(\mbox{div}u_0,\nabla\r_0,\nabla\t_0)\|^2 +\d\v^2\int_0^t\|\nabla^2u(\tau)\|^2d\tau
+\d\k(\v)^2\int_0^t\|\nabla\Delta\t(\tau)\|^2d\tau
+C_3tP(\mathcal{N}_m(t))\Big\}.
\end{align}}
\end{lemma}
\noindent\textbf{Proof}. Multiplying $\eqref{2.5}_2$ by $\nabla\mbox{div}u$, one has  that
{\small\begin{align}\label{4.2}
&\int_0^t\int(\r u_t+\r u\cdot\nabla u)\cdot\nabla\mbox{div}u dxd\tau
+\int_0^t\int\nabla p\cdot\nabla\mbox{div}u dxd\tau\nonumber\\
&=-\mu\v\int_0^t\int\nabla\times\omega\cdot\nabla\mbox{div}u dxd\tau
+(2\mu+\l)\v\int_0^t\|\nabla\mbox{div}u\|^2 d\tau.
\end{align}}
It follows from integrating by parts and  the boundary conditions  \eqref{1.9}  that
{\small\begin{align}\label{4.3}
&\int_0^t\int(\r u_t+\r u\cdot\nabla u)\cdot\nabla\mbox{div}u dxd\tau=-\int_0^t\int(\r \mbox{div}u_t+\r u\cdot\nabla \mbox{div}u)\mbox{div}u dxd\tau\nonumber\\
&~~~-\int_0^t\int(\nabla\r\cdot u_t+\nabla(\r u)^t\cdot\nabla u) \mbox{div}u dxd\tau
+\int_0^t\int_{\partial\Om}\r(u\cdot\nabla)u\cdot n\mbox{div}ud\sigma d\tau\nonumber\\
&\leq -\f12\int\r|\mbox{div}u|^2dx+\f12\int\r_0|\mbox{div}u_0|^2dx+C[1+P(Q(t))]\int_0^t\|(u_t,\nabla u)\|^2d\tau\nonumber\\
&~~~~~~~~~~~~~~~~~~~+\Big|\int_0^t\int_{\partial\Om}\r(u\cdot\nabla)n\cdot u\mbox{div}ud\sigma d\tau\Big|\nonumber\\
&\leq -\f12\int\r|\mbox{div}u|^2dx+\f12\int\r_0|\mbox{div}u_0|^2dx+C[1+P(Q(t))]\int_0^t\|(u_t,\nabla u)\|^2+|u|^2_{L^2}d\tau\nonumber\\
&\leq -\f12\int\r|\mbox{div}u|^2dx+\f12\int\r_0|\mbox{div}u_0|^2dx+CtP(\mathcal{N}_m(t)).
\end{align}}
By using the boundary condition \eqref{1.9} and integrating by parts along the boundary, one has  that
{\small\begin{align}\label{4.5}
&\v\int_0^t\int\nabla\times\omega\cdot\nabla\mbox{div}u dxd\tau
=\v\int_0^t\int_{\partial\Om} n\times\omega\cdot\nabla\mbox{div}u d\sigma d\tau\nonumber\\
&=\v\int_0^t\int_{\partial\Om} (Bu)\cdot \Pi(\nabla\mbox{div}u) d\sigma d\tau
=\v\int_0^t\int_{\partial\Om} (Bu)\cdot Z_y\mbox{div}u d\sigma d\tau\leq  C_3\v\int_0^t |u|_{H^{\f12}}|\mbox{div}u|_{H^{\f12}} d\tau\nonumber\\
&\leq C_3\v\int_0^t\|u\|_{H^1}\|\mbox{div}u\|_{H^1}d\tau\leq \f\v4\int_0^t\|\nabla\mbox{div}u\|^2d\tau+C_3\v\int_0^t\|(\nabla u,u)\|^2 d\tau.
\end{align}}
Substituting \eqref{4.3} and \eqref{4.5} into \eqref{4.2}, one obtains that
{\small\begin{align}\label{4.5-1}
&\f12\int\r|\mbox{div}u|^2dx+\f78(2\mu+\l)\v\int_0^t\|\nabla\mbox{div}u\|^2d\tau
-R\int_0^t\int\nabla(\r\t)\cdot\nabla\mbox{div}u dxd\tau\nonumber\\
&\leq \f12\int\r_0|\mbox{div}u_0|^2dx+CtP(\mathcal{N}_m(t)).
\end{align}}

On the other hand, applying $\nabla$ to $\eqref{2.5}_3$ yields that
{\small\begin{equation}\label{4.6-1}
\r\nabla\t_t+\r (u\cdot\nabla)\nabla\t+p\nabla\mbox{div}u=\k(\v)\Delta\nabla\t+\v\nabla(2\mu|Su|^2+\l|\mbox{div}u|^2)
-\Big[\nabla\r\cdot\t_t+\nabla p\cdot\mbox{div}u+\nabla(\r u)^t\nabla\t\Big].
\end{equation}}
Multiplying \eqref{4.6-1} by $\f{\nabla\t}{\t}$, one obtains that
{\small\begin{align}\label{4.6-2}
&\int_0^t\int[\r\nabla\t_t+\r (u\cdot\nabla)\nabla\t]\cdot\f{\nabla\t}{\t}dxd\tau+R\int_0^t\int\r\nabla\t\cdot\nabla\mbox{div}udxd\tau\nonumber\\
&=\k(\v)\int_0^t\int\Delta\nabla\t\cdot\f{\nabla\t}{\t}dxd\tau
+\v\int_0^t\int\nabla(2\mu|Su|^2+\l|\mbox{div}u|^2)\cdot\f{\nabla\t}{\t}dxd\tau\nonumber\\
&~~~-\int_0^t\int(\nabla\r\cdot\t_t+\nabla p\cdot\mbox{div}u+\nabla(\r u)^t\nabla\t)\cdot\f{\nabla\t}{\t}dxd\tau.
\end{align}}
For the first terms on the LHS of \eqref{4.6-2}, it follows from  integrating by parts that
{\small\begin{equation}\label{4.6-3}
\int_0^t\int[\r\nabla\r_t+\r (u\cdot\nabla)\nabla\t]\cdot\f{\nabla\t}{\t}dxd\tau\geq \int\f{\r}{2\t}|\nabla\t|^2dx-\int\f{\r_0}{2\t_0}|\nabla\t_0|^2dx
-CtP(\mathcal{N}_m(t)).
\end{equation}}
For the last two terms on the right hand side of \eqref{4.6-2}, it follows from the Cauchy inequality that
{\small\begin{align}\label{4.6-5}
&\v\Big|\int_0^t\int\nabla(2\mu|Su|^2+\l|\mbox{div}u|^2)\cdot\f{\nabla\t}{\t}dxd\tau\Big|+\Big|\int_0^t\int(\nabla\r\cdot\t_t+\nabla p\cdot\mbox{div}u+\nabla(\r u)^t\nabla\t)\cdot\f{\nabla\t}{\t}dxd\tau\Big|\nonumber\\
&\leq \d\v^2\int_0^t\|\nabla^2u\|^2d\tau+C_\d tP(\mathcal{N}_m(t)).
\end{align}}
Using the boundary condition $\eqref{1.9}_3$, one gets that
{\small\begin{align}\label{4.6-4}
&\k(\v)\int_0^t\int\nabla\Delta\t\cdot \f{\nabla\t}{\t} dxd\tau=-\k(\v)\int_0^t\int\f{|\Delta\t|^2}{\t}dxd\tau+\k(\v)\int_0^t\int\Delta\t\f{|\nabla\t|^2}{\t^2}dxd\tau
\nonumber\\
&+\k(\v)\int_0^t\int_{\partial\Omega}\Delta\t\f{\nabla\t\cdot n}{\t}d\sigma d\tau\leq -\f78\k(\v)\int_0^t\int\f{|\Delta\t|^2}{\t}dxd\tau
+C\k(\v)\int_0^t|\Delta\t|_{L^2}d\tau+CtP(\mathcal{N}_m(t))\nonumber\\
&\leq -\f78\k(\v)\int_0^t\int\f{|\Delta\t|^2}{\t}dxd\tau+C\k(\v)\int_0^t\|\Delta\t\|_{H^1}d\tau+CtP(\mathcal{N}_m(t))\\
&\leq -\f34\k(\v)\int_0^t\int\f{|\Delta\t|^2}{\t}dxd\tau+\d\k(\v)^2\int_0^t\|\nabla\Delta\t\|^2d\tau+C_\d tP(\mathcal{N}_m(t)).\nonumber
\end{align}}
Substituting \eqref{4.6-5}-\eqref{4.6-4} into \eqref{4.6-2}, one obtains that
{\small\begin{align}\label{4.6-6}
&\int\f{\r}{2\t}|\nabla\t|^2dx+R\int_0^t\int\r\nabla\t\cdot\nabla\mbox{div}udxd\tau+\f34\k(\v)\int_0^t\int\f{|\Delta\t|^2}{\t}dxd\tau \nonumber\\
&\leq \int\f{\r_0}{2\t_0}|\nabla\t_0|^2dx+ C\d\int_0^t\v^2\|\nabla^2u\|^2+\k(\v)^2\|\nabla\Delta\t\|^2d\tau+C_\d tP(\mathcal{N}_m(t)).
\end{align}}
Combining  \eqref{4.5-1} and \eqref{4.6-6}, it holds that
{\small\begin{align}\label{4.6-7}
&\f12\int\f{\r}{\t}|\nabla\t(\tau)|^2+\r|\mbox{div}u(\tau)|^2dx-R\int_0^t\int \t\nabla\r\cdot\nabla\mbox{div}udxd\tau\nonumber\\
&~~~~~~~~~~~~~~~+\f{3\k(\v)}4\int_0^t\int\f{|\Delta\t|^2}{\t}dxd\tau+\f34(2\mu+\l)\v\int_0^t\|\nabla\mbox{div}u\|^2d\tau \nonumber\\
&\leq C\Big\{\|(\nabla\t_0,\mbox{div}u_0)\|^2+ \d\int_0^t\v^2\|\nabla^2u\|^2+\k(\v)^2\|\nabla\Delta\t\|^2d\tau+C_\d t P(\mathcal{N}_m(t))\Big\}.
\end{align}}
Finally, it follows from \eqref{2.55}  that
{\small\begin{equation}\label{4.6-8}
I=-\int_0^t\int \t\nabla\r\cdot\nabla\mbox{div}udxd\tau\geq  \int\f{\t}{2\r}|\nabla\r|^2dx-\int\f{\t_0}{2\r_0}|\nabla\r_0|^2dx
-C tP(\mathcal{N}_m(t)).
\end{equation}}
Substituting \eqref{4.6-8} into \eqref{4.6-7}, one proves  \eqref{4.1}. Thus, the proof of Lemma \ref{lem4.1} is completed.
$\hfill\Box$

 \

Next we consider the higher order estimates. Firstly, we  estimate  $\MZ^\a\mbox{div}u$ for $|\a_0|\leq m-2$ with $|\a|=m-1$.
\begin{lemma}\label{lem4.2}
For every $m\geq3$ and $|\a|\leq m-1$ with $|\a_0|\leq m-2$, it holds that
{\small\begin{align}\label{4.6}
&\sup_{0\leq \tau\leq t}\|(\MZ^\a\mbox{div}u,\MZ^\a\nabla \r,\nabla \MZ^\a\t)(\tau)\|^2+\v\int_0^t\|\nabla\MZ^\a\mbox{div}u(\tau)\|^2d\tau
+\k(\v)\int_0^t\|\MZ^\a\Delta\t(\tau)\|^2d\tau\nonumber\\
&\leq C C_{m+2}\Big\{\mathcal{N}_m(0)+(\d+\v)\int_0^t\|\nabla\MZ^{m-2}\mbox{div}u(\tau)\|^2d\tau
+\d\int_0^t\|\nabla\partial_t^{m-1}(\r,\t)(\tau)\|^2d\tau\\
&+\d\int_0^t\v^2\|\nabla^2u(\tau)\|_{\mathcal{H}^{m-1}}^2
+\k(\v)^2\|\nabla\MZ^{m-2}\Delta\t(\tau)\|^2d\tau+\v\int_0^t\|\nabla^2\MZ^{m-2}u(\tau)\|^2d\tau+C_\d tP(\mathcal{N}_m(t))\Big\}.\nonumber
\end{align}}
\end{lemma}
\noindent\textbf{Proof}. The estimate for $|\a|=0$ is already given in Lemma \ref{lem4.1}. Assuming that it is proven for $|\a|\leq m-2$. We shall prove it for $|\a|=m-1\geq 1$ with $|\a_0|\leq m-2$.
Multiplying \eqref{2.26} by $\nabla\MZ^\a\mbox{div}u$ yields  that
{\small\begin{align}\label{4.7}
&\int_0^t\int(\r\mathcal{Z}^\a u_t+\r u\cdot\nabla\mathcal{Z}^\a u)\cdot\nabla\MZ^\a\mbox{div}udxd\tau+\int_0^t\int\mathcal{Z}^\a\nabla p\cdot\nabla\MZ^\a\mbox{div}udxd\tau\nonumber\\
&=-\mu\v\int_0^t\int\mathcal{Z}^\a\nabla\times\omega\cdot\nabla\MZ^\a\mbox{div}udxd\tau
+(2\mu+\l)\v\int_0^t\int\mathcal{Z}^\a\nabla\mbox{div}u\cdot\nabla\MZ^\a\mbox{div}udxd\tau\nonumber\\
&~~~~~~~~~~~~~~~~~+\int_0^t\int(\mathcal{C}_1^{\a}+\mathcal{C}_2^{\a})\cdot\nabla\MZ^\a\mbox{div}udxd\tau.
\end{align}}
Since
{\small\begin{align}\label{4.8}
&\int_0^t\int(\r\mathcal{Z}^\a u_t+\r u\cdot\nabla\mathcal{Z}^\a u)\cdot\nabla\MZ^\a\mbox{div}udxd\tau\nonumber\\
&=-\int_0^t\int(\r\mbox{div}\mathcal{Z}^\a u_t+\r u\cdot\nabla\mbox{div}\mathcal{Z}^\a u)\MZ^\a\mbox{div}udxd\tau-\int_0^t\int(\nabla\r\cdot\mathcal{Z}^\a u_t+\nabla(\r u)^t\cdot\nabla\mathcal{Z}^\a u)\MZ^\a\mbox{div}udxd\tau\nonumber\\
&~~~~+\int_0^t\int_{\partial\Om}(\r\mathcal{Z}^\a u_t\cdot n+\r(u\cdot\nabla)\mathcal{Z}^\a u\cdot n)\MZ^\a\mbox{div}ud\sigma d\tau\triangleq I_1+I_2+I_3.
\end{align}}
For $I_1$ and $I_2$, one can  obtains easily that
{\small\begin{align}\label{4.9}
I_1&=-\int_0^t\int(\r\mathcal{Z}^\a \mbox{div}u_t+\r u\cdot\nabla\mathcal{Z}^\a \mbox{div}u)\MZ^\a\mbox{div}udxd\tau\nonumber\\
&~~~~~~~~~-\int_0^t\int\Big(\r[\mbox{div},\mathcal{Z}^\a]u_t+\r(u_1Z_{y^1}+u_2Z_{y^2}+\frac{u\cdot n}{\varphi(z)}Z_3)[\mbox{div},\mathcal{Z}^\a]u\Big)\MZ^\a\mbox{div}udxd\tau\nonumber\\
&\leq -\int\frac{\r}{2}|\MZ^\a\mbox{div}u(t)|^2dx+\int\frac{\r_0}{2}|\MZ^\a\mbox{div}u_0|^2dx
+C_2[1+P(Q(t))]\int_0^t\|\nabla u\|^2_{\mathcal{H}^{m-1}}d\tau,
\end{align}}
and
{\small\begin{eqnarray}\label{4.10}
I_2\leq C[1+P(Q(t))]\int_0^t\|\nabla u\|^2_{\mathcal{H}^{m-1}}d\tau,
\end{eqnarray}}
where $\varphi(z)=\f{z}{1+z}$.
Noting that $\MZ^\a$ contains at least one tangential derivative $Z_y$, integrating by parts along the boundary and using \eqref{3.4}, \eqref{2.51}, one has  that
{\small\begin{align}\label{4.12}
 I_3&=\int_0^t\int_{\partial\Om}[\r\mathcal{Z}^\a u_t\cdot n-\r(u\cdot\nabla)n\cdot\mathcal{Z}^\a u+\r(u\cdot\nabla)(\mathcal{Z}^\a u\cdot n)]\MZ^\a\mbox{div}ud\sigma d\tau\nonumber\\
&~~~=\int_0^t\int_{\partial\Om}[\r\mathcal{Z}^\a u_t\cdot n-\r(u\cdot\nabla)n\cdot\mathcal{Z}^\a u+\r(u_1\partial_{y^1}+u_2\partial_{y^2})(\mathcal{Z}^\a u\cdot n)]\MZ^\a\mbox{div}ud\sigma d\tau\nonumber\\
&~~~\leq  C[1+P(Q(t))]\int_0^t\Big(|\mathcal{Z}^\a u_t\cdot n|_{H^{\f12}}+|\mathcal{Z}^\a u\cdot n|_{H^{\f32}}+|\mathcal{Z}^\a u|_{H^{\f12}}\Big)\cdot|\MZ^{m-2}\mbox{div}u|_{H^{\f12}}d\tau\nonumber\\
&~~\leq  \d\int_0^t\|\nabla\MZ^{m-2}\mbox{div}u\|^2d\tau+ C_{m+2}C_\d t P(\mathcal{N}_m(t)).
\end{align}}
Substituting \eqref{4.9}, \eqref{4.10} and \eqref{4.12} into \eqref{4.8}, one gets that
{\small\begin{align}\label{4.13}
&\int_0^t\int(\r\mathcal{Z}^\a u_t+\r u\cdot\nabla\mathcal{Z}^\a u)\nabla\MZ^\a\mbox{div}udxd\tau\\
&\leq -\int\frac{\r}{2}|\MZ^\a\mbox{div}u(t)|^2dx+\int\frac{\r_0}{2}|\MZ^\a\mbox{div}u_0|^2dx
+\d\int_0^t\|\nabla\MZ^{m-2}\mbox{div}u\|^2d\tau+C_\d C_{m+2}t P(\mathcal{N}_m(t)).\nonumber
\end{align}}
By the same argument as Lemma 3.6 of \cite{Wang-Xin-Yong}, one can obtain that
{\small\begin{align}\label{4.14}
\v\int_0^t\int\mathcal{Z}^\a\nabla\mbox{div}u\cdot\nabla\MZ^\a\mbox{div}udxd\tau
\geq \f34\v\int_0^t\|\nabla\MZ^\a\mbox{div}u\|^2d\tau
-C\v\int_0^t\|\nabla\MZ^{m-2}\mbox{div}u\|^2d\tau,
\end{align}}
and
{\small\begin{align}\label{4.16}
&-\v\int_0^t\int\mathcal{Z}^\a\nabla\times\omega\cdot\nabla\MZ^\a\mbox{div}udxd\tau\nonumber\\
&\geq -\f\v4\int_0^t\|\nabla\MZ^\a\mbox{div}u\|^2d\tau -C\v\int_0^t\|\nabla\MZ^{m-2}\omega\|^2d\tau - C_{m+2}t P(\mathcal{N}_m(t)).
\end{align}}
It follows from  Proposition \ref{prop3.2} that
{\small\begin{equation}\label{4.24}
\int_0^t\|\mathcal{C}_1^{\a}\|^2_1+\|\mathcal{C}_2^{\a}\|^2_1d\tau
\leq C(1+P(Q(t)))\int_0^tP(\mathcal{N}_m(\tau))d\tau\leq CtP(\mathcal{N}_m(t)).
\end{equation}}
which, together with integrating by parts, yields  that
{\small\begin{align}\label{4.23}
&|\int_0^t\int(\mathcal{C}_1^{\a}+\mathcal{C}_2^{\a})\nabla\MZ^\a\mbox{div}udxd\tau|\nonumber\\
&\leq |\int_0^t\int(\mathcal{C}_1^{\a}+\mathcal{C}_2^{\a})Z\nabla\MZ^{\a-1}\mbox{div}udxd\tau|
+C\int_0^t\int(|\mathcal{C}_1^{\a}|+|\mathcal{C}_2^{\a}|)|\nabla\MZ^{\a-1}\mbox{div}u|dxd\tau\\
&\leq |\int_0^t\int(|Z\mathcal{C}_1^{\a}|+|Z\mathcal{C}_2^{\a}|+|\mathcal{C}_1^{\a}|+|\mathcal{C}_2^{\a}|)
\cdot|\nabla\MZ^{\a-1}\mbox{div}u|dxd\tau\nonumber\\
&\leq \d\int_0^t\|\nabla\MZ^{m-2}\mbox{div}u\|^2d\tau+C\int_0^t\|\mathcal{C}_1^{\a}\|^2_{H^1_{co}}+\|\mathcal{C}_2^{\a}\|^2_{H^1_{co}}d\tau\leq \d\int_0^t\|\nabla\MZ^{m-2}\mbox{div}u\|^2d\tau+CtP(\mathcal{N}_m(t)),\nonumber
\end{align}}
Substituting \eqref{4.13}-\eqref{4.23} into \eqref{4.7}, one obtains that
{\small\begin{align}\label{4.13-1-1}
&\int\frac{\r}{2}|\MZ^\a\mbox{div}u(t)|^2dx-R\int_0^t\int\mathcal{Z}^\a\nabla (\r\t)\cdot\nabla\MZ^\a\mbox{div}udxd\tau+\f34(2\mu+\l)\v\int_0^t\|\nabla\MZ^\a\mbox{div}u\|^2d\tau\\
&\leq \int\frac{\r_0}{2}|\MZ^\a\mbox{div}u_0|^2dx
+C\Big\{(\d+\v)\int_0^t\|\nabla\MZ^{m-2}\mbox{div}u\|^2d\tau
+\v\int_0^t\|\nabla\MZ^{m-2}\omega\|^2d\tau+C_{m+2}C_\d tP(\mathcal{N}_m(t))\Big\}.\nonumber
\end{align}}

Next we shall estimate the temperature part. Applying $\nabla$ to $\eqref{2.26}_2$, one obtains that
{\small\begin{align}\label{4.13-1}
&\r\nabla\mathcal{Z}^\a \t_t+\r (u\cdot\nabla)\nabla\mathcal{Z}^\a \t+p\nabla\mathcal{Z}^\a\mbox{div}u
-\k(\v)\nabla\MZ^\a\Delta\t\nonumber\\
&=\nabla\r \MZ^\a\t_t+\nabla(\r u)^t\nabla\MZ^\a\t+\nabla p \MZ^\a\mbox{div}u +\v\nabla\mathcal{Z}^\a\Big(2\mu |Su|^2+\l|\mbox{div}u|^2\Big)+\nabla(\mathcal{C}_3^{\a}+\mathcal{C}_4^{\a}+\mathcal{C}_4^{\a}),
\end{align}}
with $|\a|=m-1$ and $|\a_0|\leq m-2$. Multiplying \eqref{4.13-1} by $\f{\nabla\MZ^\a\t}{\t}$ yields that
{\small\begin{align}\label{4.13-2}
&\int\f{\r}{2\t}|\nabla\MZ^\a\t|^2dx+R\int_0^t\int\r\nabla\MZ^\a\t\cdot\nabla\MZ^\a\mbox{div}udxd\tau
-\k(\v)\int_0^t\int\nabla\MZ^\a\Delta\t\f{\nabla\MZ^\a\t}{\t}dxd\tau\\
&\leq \v\int_0^t\int\nabla\mathcal{Z}^\a\Big(2\mu |Su|^2+\l|\mbox{div}u|^2\Big)\f{\nabla\MZ^\a\t}{\t}dxd\tau
+\int_0^t\int\nabla(\mathcal{C}_3^{\a}+\mathcal{C}_4^{\a}+\mathcal{C}_5^{\a})\f{\nabla\MZ^\a\t}{\t}dxd\tau+CtP(\mathcal{N}_m).\nonumber
\end{align}}
It follows from integrating by parts that
{\small\begin{align}\label{4.13-3}
&\k(\v)\int_0^t\int\nabla\MZ^\a\Delta\t\f{\nabla\MZ^\a\t}{\t}dxd\tau=-\k(\v)\int_0^t\int\MZ^\a\Delta\t\f{\Delta\MZ^\a\t}{\t}dxd\tau+\k(\v)\int_0^t\int \MZ^\a\Delta\t\f{\nabla\MZ^\a\t}{\t^2}\nabla\t dxd\tau\nonumber\\
&+\k(\v)\int_0^t\int_{\partial\Omega}\MZ^\a\Delta\t\f{ n\cdot\nabla\MZ^\a\t}\t d\sigma d\tau\leq -\f{3\k(\v)}4\int_0^t\int\f{|\MZ^\a\Delta\t|^2}{\t}dxd\tau +\k(\v)\int_0^t\int_{\partial\Omega}\MZ^\a\Delta\t\f{ n\cdot\nabla\MZ^\a\t}\t d\sigma d\tau\nonumber\\
&~~~~~~~~~~~~~~~~~~~~~~~~~~~~~~~~~~~~~~~~~~~~~~~~~~~+C\k(\v)\int_0^t\|\MZ^{m-2}\Delta\t\|^2d\tau.
+CtP(\mathcal{N}_{m}(t)).
\end{align}}
Noting that $\MZ^\a$ contains at least one tangential derivative $Z_y$, integrating by parts along the boundary yields   that
{\small\begin{align}\label{4.13-6}
&\k(\v)\Big|\int_0^t\int_{\partial\Omega}\MZ^\a\Delta\t\f{ n\cdot\nabla\MZ^\a\t}\t d\sigma d\tau\Big|\nonumber\\
&=\k(\v)\Big|\int_0^t\int_{\partial\Omega}\MZ^{\a-1}\Delta\t\cdot
\Big(\f1\t Z_y(n\cdot\nabla\MZ^\a\t)+n\cdot\nabla\MZ^\a\t Z_y(\f1\t)\Big) d\sigma d\tau\Big|\nonumber\\
&\leq CP(Q(t))\k(\v)\int_0^t\Big[\|\MZ^{m-2}\Delta\t\|^{\f12}_{L^2}\|\nabla\MZ^{m-2}\Delta\t\|^{\f12}_{L^2}+\|\MZ^{m-2}\Delta\t\|_{L^2}\Big]\Big( |n\cdot\nabla\MZ^\a\t|_{H^1}+|n\cdot\nabla\MZ^\a\t|_{L^2} \Big)d\tau\nonumber\\
&\leq \d\k(\v)^2\int_0^t\|\nabla\MZ^{m-2}\Delta\t\|^2d\tau+\sum_{|\b|\leq m-1}^{|\b_0|\leq m-2}\d\k(\v)^2\int_0^t\|\MZ^{\b}\Delta\t\|^2d\tau
+\d\k(\v)\int_0^t\|\nabla\MZ^m\t\|^2d\tau\nonumber\\
&~~~~~~~+C\k(\v)\int_0^t\|\MZ^{m-2}\Delta\t\|^2d\tau+C_\d C_{m+1}tP(\mathcal{N}_m(t)),
\end{align}}
where we have used the following boundary estimates, for $|\a|\leq m-1$ with $|\a_0|\leq m-2$,
{\small\begin{align}\label{4.13-5}
\begin{cases}
|n\cdot\nabla\MZ^\a\t|_{L^2(\partial\Omega)}
\leq CC_{m}\Big(\mathcal{N}_m^{\f12}+\|\MZ^{m-2}\Delta\t\|^{\f12}\mathcal{N}_m^{\f14}\Big),\\
|n\cdot\nabla\MZ^\a\t|_{H^1(\partial\Omega)}
\leq CC_{m+1}\Big(\|\nabla\MZ^m\t\|^{\f12}+\sum_{|\b|\leq m-1}^{|\b_0|\leq m-2}\|\MZ^{\b}\Delta\t\|^{\f12}\Big)\mathcal{N}_m^{\f14}+\mathcal{N}_m^{\f12}, 
\end{cases}
\end{align}}
which follows from  the boundary condition \eqref{1.9} and \eqref{3.4}.
Substituting \eqref{4.13-6} into \eqref{4.13-3}, one obtains that
{\small\begin{align}\label{4.13-7}
&-\k(\v)\int_0^t\int\nabla\MZ^\a\Delta\t\f{\nabla\MZ^\a\t}{\t}dxd\tau\nonumber\\
&\geq \f{\k(\v)}2\int_0^t\int\f{|\MZ^\a\Delta\t|^2}{\t}dxd\tau-
\d\k(\v)^2\int_0^t\|\nabla\MZ^{m-2}\Delta\t\|^2d\tau-\sum_{|\b|\leq m-1}^{|\b_0|\leq m-2}\d\k(\v)^2\int_0^t\|\MZ^{\b}\Delta\t\|^2d\tau
\nonumber\\
&~~~-\d\k(\v)\int_0^t\|\nabla\MZ^m\t\|^2d\tau-C\k(\v)\int_0^t\|\MZ^{m-2}\Delta\t\|^2d\tau-C_\d C_{m+1}tP(\mathcal{N}_m(t)).
\end{align}}

For the  terms on the RHS of \eqref{4.13-2}, it follows from  Proposition \ref{prop3.2} and Holder inequality that
{\small\begin{align}\label{4.13-9}
&\v\Big|\int_0^t\int\nabla\mathcal{Z}^\a\Big(2\mu |Su|^2+\l|\mbox{div}u|^2\Big)\cdot\f{\nabla\MZ^\a\t}{\t}dxd\tau\Big|
+\int_0^t\|(\nabla C_3^{\a},\nabla C_5^{\a})\|^2d\tau\leq \d\v^2\int_0^t\|\nabla^2u\|^2_{\mathcal{H}^{m-1}}d\tau\nonumber\\
&+C[P(\mathcal{N}_m(t))+\|\nabla\mbox{div}u\|^2_{L^\infty}]
\int_0^tP(\mathcal{N}_m(\tau))+\|\nabla\MZ^{m-2}\mbox{div}u,\nabla\partial_t^{m-1}(\r,\t)\|^2 d\tau+C_\d C_{m+1}tP(\mathcal{N}_m(t))\nonumber\\
&\leq \d\v^2\int_0^t\|\nabla^2u\|^2_{\mathcal{H}^{m-1}}+CP(\mathcal{N}_m(t))
\int_0^t\|\nabla\MZ^{m-2}\mbox{div}u,\nabla\partial_t^{m-1}(\r,\t)\|^2 d\tau+C_\d C_{m+1}tP(\mathcal{N}_m(t)),
\end{align}}
where we have used the fact
$\|\nabla\mbox{div}u\|^2_{L^\infty}\leq C P(\mathcal{N}_m(t))$
which will be proved in Lemma \ref{lem6.1} below. For the term $\nabla C_4^{\a}$, one needs to be more careful. First, one notices that
{\small\begin{equation}\label{4.13-10}
C_4^{\a}=\sum_{|\b|\geq1,\b+\g=\a}\sum_{i=1}^2C_{\a,\b}\MZ^{\b}(\r u_i)\MZ^\g\partial_{y^i}\t
+\sum_{|\b|\geq1,\b+\g=\a}C_{\a,\b}\MZ^{\b}(\r u\cdot N)\MZ^\g\partial_{z}\t+\r u\cdot N [\MZ^\a,\partial_z]\t,
\end{equation}}
then, from Proposition \ref{prop3.2}, it holds that
{\small\begin{align}\label{4.13-11}
&\int_0^t\|\nabla(\r u\cdot N[\MZ^\a,\partial_z]\t)\|^2+\sum_{|\b|\geq1,\b+\g=\a}\Big[\sum_{i=1}^2\|\nabla(\MZ^{\b}(\r u_i)\MZ^\g\partial_{y^i}\t)\|^2
+\|\MZ^\g\partial_z\t\cdot\nabla\MZ^\b(\r u\cdot N)\|^2\Big]d\tau
\nonumber\\
&\leq C(1+P(Q(t)))\int_0^t P(\mathcal{N}_m(\tau))+\|\nabla\partial_t^{m-1}\r\|^2d\tau.
\end{align}}
For $|\b|\geq1,\b+\g=\a$, and $|\a|=m-1$, one notices that
{\small\begin{align}\label{4.13-12}
\MZ^\b(\r u\cdot N) \nabla\MZ^\g\partial_z\t
=\sum_{\tilde{\b}\leq \b}C_{\tilde{\b}}(z)
\MZ^{\tilde{\b}}(\f{\r u\cdot N}{\varphi(z)})\cdot\varphi(z)\nabla\MZ^{\g}\partial_z\t,
\end{align}}
where $C_{\tilde{\b}}(z)$ is bounded smooth function of $z$.
If $\tilde\b=0$ and $|\g|\leq m-2$, one gets  that
{\small\begin{align}\label{4.13-13}
\int_0^t\|\MZ^{\tilde{\b}}(\f{\r u\cdot N}{\varphi(z)})\cdot\varphi(z)\nabla\MZ^{\g}\partial_z\t\|^2d\tau
&\leq \|\f{\r u\cdot N}{\varphi(z)}\|^2_{L^\infty}\int_0^t\|\varphi(z)\nabla\MZ^{\g}\partial_z\t\|^2d\tau \leq  CtP(\mathcal{N}_m(t)),
\end{align}}
If $|\tilde\b|\neq0$,  from Proposition \ref{prop3.2}, one obtains that
{\small\begin{equation}\label{4.13-14}
\int_0^t\|\MZ^\b(\f{\r u\cdot N}{\varphi(z)})\nabla\MZ^\g\partial_z\t\|^2d\tau
\leq C(1+P(Q(t)))\int_0^t P(\mathcal{N}_m)+\|\f{\r u\cdot N}{\varphi(z)}\|_{\mathcal{H}^{m-1}}^2d\tau
\leq CtP(\mathcal{N}_m(t)),
\end{equation}}
where in the last inequality, we have used the  Hardy inequality
{\small\begin{equation}\label{4.13-15}
\|\f{u\cdot N}{\varphi(z)}\|_{\mathcal{H}^{m-1}}^2\leq C_{m+1}\|\nabla u\|_{\mathcal{H}^{m-1}}^2,
\end{equation}}
which has already been proved in page 543 of \cite{Masmoudi-R}. Then, combining \eqref{4.13-12}-\eqref{4.13-14}, one obtains that
{\small\begin{equation}\label{4.13-16}
\sum_{|\b|\geq1,\b+\g=\a}\int_0^t\|\MZ^\b(\r u\cdot N) \nabla\MZ^\g\partial_z\t\|^2d\tau\leq Ct P(\mathcal{N}_m(t)).
\end{equation}}
Thus, it follows from \eqref{4.13-10}, \eqref{4.13-11} and \eqref{4.13-16}  that
{\small\begin{eqnarray}\label{4.13-17}
&&\int_0^t\|\nabla C_4^{\a}\|^2d\tau\leq C P(\mathcal{N}_m(\tau))
\int_0^t\|\nabla\partial_t^{m-1}(\r,\t)\|^2 d\tau+CtP(\mathcal{N}_m(t))
\end{eqnarray}}
Substituting \eqref{4.13-7}, \eqref{4.13-9} and \eqref{4.13-17} into  \eqref{4.13-2}  and using Holder inequality, one obtains  that
{\small\begin{align}\label{4.13-112}
&\int\f{\r}{2\t}|\nabla\MZ^\a\t|^2dx+R\int_0^t\int\r\nabla\MZ^\a\t\cdot\nabla\MZ^\a\mbox{div}udxd\tau
+\f{\k(\v)}2\int_0^t\int\f{|\MZ^\a\Delta\t|^2}{\t}dxd\tau\nonumber\\
&\leq CC_{m+1}\Big\{\|\nabla\MZ^\a\t_0\|^2+\d\v^2\int_0^t\|\nabla^2u\|^2_{\mathcal{H}^{m-1}}d\tau
+\d\k(\v)^2\int_0^t\|\nabla\MZ^{m-2}\Delta\t\|^2+\sum_{|\b|\leq m-1}^{|\b_0|\leq m-2}\|\MZ^{\b}\Delta\t\|^2d\tau
\nonumber\\
&~+\d\k(\v)\int_0^t\|\nabla\MZ^m\t\|^2d\tau+\k(\v)\int_0^t\|\MZ^{m-2}\Delta\t\|^2d\tau
+\d\int_0^t \|\nabla\MZ^{m-2}\mbox{div}u\|^2+\|\nabla\partial_t^{m-1}(\r,\t)\|^2 d\tau \nonumber\\
&~~~~~~~~~~~~~~~~~~~~~~~~+C_\d t P(\mathcal{N}_m(t))\Big\}.
\end{align}}
It follows from  \eqref{4.13-1-1} and \eqref{4.13-112} that
{\small\begin{align}\label{4.13-18}
&\int\frac{\r}{2}|\MZ^\a\mbox{div}u(t)|^2+\f{\r}{2\t}|\nabla\MZ^\a\t|^2dx
+R\int_0^t\int(\r\nabla\MZ^\a\t-\mathcal{Z}^\a\nabla (\r\t))\cdot\nabla\MZ^\a\mbox{div}udxd\tau\nonumber\\
&+\f34(2\mu+\l)\v\int_0^t\|\nabla\MZ^\a\mbox{div}u\|^2d\tau+\f{\k(\v)}2\int_0^t\int\f{|\MZ^\a\Delta\t|^2}{\t}dxd\tau\nonumber\\
&\leq CC_{m+1}\Big\{\mathcal{N}_m(0)+\d\v^2\int_0^t\|\nabla^2u\|^2_{\mathcal{H}^{m-1}}d\tau
+\d\k(\v)^2\int_0^t\|\nabla\MZ^{m-2}\Delta\t\|^2+\sum_{|\b|\leq m-1}^{|\b_0|\leq m-2}\|\MZ^{\b}\Delta\t\|^2d\tau
\nonumber\\
&~~~~~~~~~+\d\k(\v)\int_0^t\|\nabla\MZ^m\t\|^2d\tau+C\int_0^t\k(\v)\|\MZ^{m-2}\Delta\t\|^2+\v\|\nabla^2\MZ^{m-2} u\|^2d\tau
 \nonumber\\
&~~~~~~~~~~~~~+(\d+\v)\int_0^t \|\nabla\MZ^{m-2}\mbox{div}u\|^2d\tau+\d\int_0^t\|\nabla\partial_t^{m-1}(\r,\t)\|^2 d\tau+C_\d t P(\mathcal{N}_m(t))\Big\}.
\end{align}}
In order to close the estimate of \eqref{4.13-18}, one notes that
{\small\begin{align}\label{4.18}
	&\r\nabla\MZ^\a\t-\mathcal{Z}^\a\nabla (\r\t)=-\t\MZ^\a\nabla\r-\r[\MZ^\a,\nabla]\t-\sum_{|\b|\geq1,\beta+\g=\a}C_{\a,\b}\Big(\MZ^\b\t\cdot\MZ^\g\nabla\r+\MZ^\b\r\cdot\MZ^\g\nabla\t\Big).
\end{align}}
Since $\MZ^\a\neq\partial_t^{m-1}$, it follows from \eqref{3.2}, \eqref{4.18} and  integrating by parts that
{\small\begin{align}\label{4.18-1}
& \int_0^t\int(\r\nabla\MZ^\a\t-\mathcal{Z}^\a\nabla (\r\t))\cdot\nabla\MZ^\a\mbox{div}udxd\tau=\int_0^t\int(\r\nabla\MZ^\a\t-\mathcal{Z}^\a\nabla (\r\t))\cdot\Big(\MZ^\a\nabla\mbox{div}u+[\nabla,\MZ^\a]\mbox{div}u\Big)dxd\tau \nonumber\\
&\geq \int_0^t\int(\r\nabla\MZ^\a\t-\mathcal{Z}^\a\nabla (\r\t))\cdot \MZ^\a\nabla\mbox{div}u dxd\tau -\d\int_0^t\|\nabla\MZ^{m-2}\mbox{div}u\|^2d\tau-C_\d tP(\mathcal{N}_m(t))
\nonumber\\
&\geq -\int_0^t\int\t\MZ^\a\nabla\r\cdot \MZ^\a\nabla\mbox{div}u dxd\tau-\d\int_0^t\|\nabla\MZ^{m-2}\mbox{div}u\|^2d\tau-C_\d tP(\mathcal{N}_m(t))
\nonumber\\
&~~~~~~-\sum_{|\b|\geq1,\beta+\g=\a}C_{\a,\b}\int_0^t\int\left| Z\Big(\MZ^\b\t\cdot\MZ^\g\nabla\r+\MZ^\b\r\cdot\MZ^\g\nabla\t\Big)\cdot \MZ^{\a-1}\nabla\mbox{div}u\right| dxd\tau\\
&\geq \int_0^t\int\t\MZ^\a\nabla\r\cdot \MZ^\a\nabla(\f{\r_t}{\r}) +\int_0^t\int\t\MZ^\a\nabla\r\cdot \MZ^\a\nabla(\f{u\cdot\nabla\r}{\r})-C\d\int_0^t\|\nabla\MZ^{m-2}\mbox{div}u\|^2d\tau-C_\d tP(\mathcal{N}_m(t)).\nonumber
\end{align}}
where we have used \eqref{2.55} in the last inequality. For the first term on the right hand side of \eqref{4.17}, one notices that
{\small\begin{equation}\label{4.18-3}
\MZ^\a\nabla(\f{\r_t}{ \r})=\f{1}{ \r}\mathcal{Z}^{\a}\nabla \r_t
+\sum_{|\beta|\geq1, \b+\g=\a}C_{\a,\b}\MZ^\b(\f1{\r})\cdot\MZ^{\g}\nabla\r_t+\sum_{\b+\g=\a}C_{\a,\b}\MZ^\b\r_t\cdot\MZ^{\g}\nabla(\f1{\r}).
\end{equation}}
Therefore, using \eqref{4.18-3} and Proposition \ref{prop3.2}, one has that
{\small\begin{align}\label{4.19}
\int_0^t\int\t\mathcal{Z}^\a\nabla \r\cdot\MZ^\a\nabla(\f{\r_t}{\r})dxd\tau
&\leq \int\f{\t}{2\r}|\mathcal{Z}^\a\nabla \r|^2 dx-\int\f{\t_0}{2\r_0}|\mathcal{Z}^\a\nabla\r_0|^2dx+\d\int_0^t\|\nabla\partial_t^{m-1}\r\|^2d\tau   \nonumber\\
&~~~~~~~+C_\d tP(\mathcal{N}_m(t)).
\end{align}}
To estimate the second term on the RHS of \eqref{4.18-1}, note that
{\small\begin{align}\label{4.20}
&\MZ^\a\nabla(\f{u}{\r}\cdot\nabla\r)=\sum_{i=1,2}\f{u_i}{\r}\mathcal{Z}^{\a}\nabla\partial_{y^i}\r
+\f{u\cdot N}{\r}\mathcal{Z}^{\a}\partial_z\nabla\r\nonumber\\
&~~~~~+\sum_{i=1,2}\sum_{|\beta|\geq1, \b+\g=\a}C_{\a,\b}\MZ^\b(\f{u_i}{\r})\cdot\MZ^{\g}\nabla\partial_{y^i} \r+\sum_{|\beta|\geq1, \b+\g=\a}C_{\a,\b}\MZ^\b(\f{u\cdot N}{\r})\cdot\MZ^{\g}\partial_z\nabla\r\nonumber\\
&~~~~~+\sum_{i=1,2}\sum_{\b+\g=\a}C_{\a,\b}\MZ^\b\nabla(\f{u_i}{\r})\cdot\MZ^{\g}\partial_{y^i}\r
+\sum_{\b+\g=\a}C_{\a,\b}\MZ^\b\nabla(\f{u\cdot N}{\r})\cdot\MZ^{\g}\partial_z\r.
\end{align}}
Integrating by parts yields immediately that
{\small\begin{align}\label{4.21-1}
&\int_0^t\int\t\mathcal{Z}^\a\nabla\r\cdot\Big(\sum_{i=1,2}\f{u_i}{\r}\mathcal{Z}^{\a}\nabla\partial_{y^i}\r+\f{u\cdot N}{\r}\mathcal{Z}^{\a}\partial_z\nabla\r\Big) dxd\tau\nonumber\\
&=-\int_0^t\int\t\mathcal{Z}^\a\nabla \r\cdot\Big(\sum_{i=1,2}\f{u_i}{\r}\partial_{y^i} \mathcal{Z}^{\a}\nabla\r+\f{u\cdot N}{\r}\partial_z\mathcal{Z}^{\a}\nabla\r\Big) dxd\tau\nonumber\\
&~~-\int_0^t\int\t\mathcal{Z}^\a\nabla\r\cdot\Big(\sum_{i=1,2}\f{u_i}{\r}[\mathcal{Z}^{\a}\nabla,\partial_{y^i}] \r+\f{u\cdot N}{\r\varphi(z)}\varphi(z)[\mathcal{Z}^{\a}, \partial_z]\nabla\r\Big) dxd\tau\leq C_2C_\d tP(\mathcal{N}_m(t)).
\end{align}}
It follows from  Proposition \ref{prop3.2}  that
{\small\begin{align}\label{4.21-2}
&\sum_{i=1,2}\sum_{\b+\g=\a}C_{\a,\b}\int_0^t\int\t\mathcal{Z}^\a\nabla\r\cdot\Big(\MZ^\b\nabla(\f{u_i}{\r})\cdot\MZ^{\g}\partial_{y^i} \r+\MZ^\b(\f{u_i}{\r})\cdot\MZ^{\g}\nabla\partial_{y^i}\r\Big) dxd\tau\nonumber\\
&-\sum_{\b+\g=\a}C_{\a,\b}\int_0^t\int\t\mathcal{Z}^\a\nabla\r\cdot\MZ^\b\nabla(\f{u\cdot N}{\r})\cdot\MZ^{\g}\partial_z\r dxd\tau\leq \d\int_0^t\|\nabla\partial_t^{m-1}\r\|^2d\tau +C_\d tP(\mathcal{N}_m(t)),
\end{align}}
On the other hand, note that for $ |\beta|\geq1, \b+\g=\a,$ and $|\a|=m-1$
{\small\begin{equation}\label{4.21-3}
\MZ^\b(\f{u\cdot N}{\r})\cdot\MZ^{\g}\partial_z\nabla\r=\sum_{\tilde{\beta}\leq \beta,\tilde{\g}\leq \g}C_{\a,\tilde\beta,\tilde\g}(z)\MZ^{\tilde\b}(\f{u\cdot N}{ \r\varphi(z)})\cdot\MZ^{\tilde\g}(Z_3\nabla\r),
\end{equation}}
where $|\tilde\beta|+|\tilde\g|\leq m-1$, $|\tilde\g|\leq m-2$ and $C_{\a,\tilde\beta,\tilde\g}(z) $ is some smooth bounded function of $z$. Using \eqref{4.21-3} and   similar arguments as in the proof of \eqref{4.13-16}, one has that
{\small\begin{align}\label{4.21-6}
&\sum_{\tilde{\beta}\leq \beta,\tilde{\g}\leq \g}\int_0^t\int C_{\a,\tilde\beta,\tilde\g}(z)\t\mathcal{Z}^\a\nabla \r\cdot \MZ^\b(\f{u\cdot N}{\r})\cdot\MZ^{\g}\partial_z\nabla\r dxd\tau\\
&\leq C\Big(\int_0^t\|\mathcal{Z}^\a\nabla \r\|^2d\tau\Big)^{\f12}
\cdot
\Big\{\sum_{|\b|\geq1,\beta+\g=\a}\int_0^t\| \MZ^\b(\f{u\cdot N}{\r})\cdot\MZ^{\g}\partial_z\nabla\r\|^2 \Big\}^{\f12} \leq CC_\d tP(\mathcal{N}_m(t)).\nonumber
\end{align}}
Combining  \eqref{4.20}-\eqref{4.21-6}, one obtains that
{\small\begin{equation}\label{4.21}
\int_0^t\int\t\mathcal{Z}^\a\nabla \r\cdot\MZ^\a\nabla(\f{u}{\r}\cdot\nabla\r)dxd\tau\leq C\d\int_0^t\|\nabla\partial_t^{m-1}\r\|^2d\tau+C_\d tP(\mathcal{N}_m(t)).
\end{equation}}
Then, substituting  \eqref{4.19} and \eqref{4.21} into \eqref{4.18-1}, one obtains that
{\small\begin{align}\label{4.17}
& \int_0^t\int(\r\nabla\MZ^\a\t-\mathcal{Z}^\a\nabla (\r\t))\cdot\nabla\MZ^\a\mbox{div}udxd\tau\leq \int\f{\t}{2\r}|\MZ^\a\nabla\r|^2dxd\tau-\int\f{\t_0}{2\r_0}|\MZ^\a\nabla\r_0|^2dxd\tau
\nonumber\\
&~~~~~~~~~~~~~~~~~~~~~-C\d\int_0^t\|\nabla\MZ^{m-2}\mbox{div}u\|^2d\tau-C\d\int_0^t\|\nabla\partial_t^{m-1}\r\|^2d\tau
-C_\d tP(\mathcal{N}_m(t)).
\end{align}}
Substituting \eqref{4.17} into  \eqref{4.13-18}, we proved \eqref{4.6}. Thus,   the proof of Lemma \ref{lem4.2} is completed.  $\hfill\Box$

In the proof of Lemma \ref{lem4.2},  we have used the fact $|\a_0|\leq m-2$ in \eqref{4.12} and \eqref{4.13-6}. However,  for the case  $\MZ^\a=\partial_t^{m-1}$,   arguments as \eqref{4.12} and \eqref{4.13-6} are not available  anymore.  And we can obtain only  the following weak estimate $\v\|\partial_t^{m-1}(\mbox{div}u,\nabla \r)\|^2$, but the control of $\v^2\int_0^t\|\nabla\partial_t^{m-1}\mbox{div}u(\tau)\|^2d\tau$ is crucial for us to close the a priori estimation.
\begin{lemma}\label{lem4.3}
For every $m\geq1$, it holds that
{\small\begin{align}\label{4.25}
&\sup_{0\leq \tau\leq t}\Big(\v\|(\partial_t^{m-1}\mbox{div}u,\nabla\partial_t^{m-1}\r)(\tau)\|^2\Big)+\f{1}{2}(2\mu+\l)\v^2\int_0^t\|\nabla\partial_t^{m-1}\mbox{div}u(\tau)\|^2d\tau\nonumber\\
&\leq C\Big\{ \v\|(\partial_t^{m-1}\mbox{div}u_0,\nabla\partial_t^{m-1}\r_0)\|^2+\int_0^t\|\partial_t^{m-1}\nabla\t\|^2d\tau+C_\d C_{m+1} tP(\mathcal{N}_m(t))\Big\}.
\end{align}}
\end{lemma}
\noindent\textbf{Proof}. First, it follows from the boundary condition \eqref{1.9} that
{\small\begin{equation}\label{4.29}
n\cdot\partial_t^{m-1}u=0,~~n\times\partial_t^{m-1}\omega=[B\partial_t^{m-1}u]_\tau,~~n\cdot\nabla\partial_t^{m-1}\t=\nu\partial_t^{m-1}\t.
\end{equation}}
Multiplying $\eqref{2.26}_1$(with $\MZ^\a=\partial_t^m$) by $\v\nabla\mbox{div}\partial_t^{m-1}u$, one obtains that
{\small\begin{align}\label{4.30}
&\v\int_0^t\int(\r\partial_t^{m-1}u_t+\r u\cdot\nabla\partial_t^{m-1}u)\nabla\mbox{div}\partial_t^{m-1}udxd\tau
+\v\int_0^t\int\partial_t^{m-1}\nabla p\cdot\nabla\mbox{div}\partial_t^{m-1}udxd\tau\nonumber\\
&=-\mu\v^2\int_0^t\int\nabla\times\partial_t^{m-1}\omega\cdot\nabla\mbox{div}\partial_t^{m-1}udxd\tau
+(2\mu+\l)\v^2\int_0^t\|\nabla\partial_t^{m-1}\mbox{div}u\|^2d\tau\\
&~~~~~+\v\int_0^t\int(\mathcal{C}_1^{\a}+\mathcal{C}_2^{\a})\nabla\mbox{div}\partial_t^{m-1}u.\nonumber
\end{align}}
It follows from  \eqref{4.29} and  integrating by parts  that
{\small\begin{align}\label{4.32}
&\v^2|\int_0^t\int\nabla\times\partial_t^{m-1}\omega\cdot\nabla\mbox{div}\partial_t^{m-1}udxd\tau|
=\v^2|\int_0^t\int_{\partial\Omega} n\times\partial_t^{m-1}\omega\cdot\Pi(\nabla\mbox{div}\partial_t^{m-1}u)d\sigma d\tau|\nonumber\\
&\leq C\v^2\int_0^t|n\times\partial_t^{m-1}\omega|_{H^{\f12}}\cdot|\partial_t^{m-1}\mbox{div}u|_{H^{\f12}}d\tau
\leq C  C_3 \v^2\int_0^t|\partial_t^{m-1}u|_{H^{\f12}}\cdot|\partial_t^{m-1}\mbox{div}u|_{H^{\f12}}d\tau\nonumber\\
&\leq \f1{16}\v^2\int_0^t\|\nabla\partial_t^{m-1}\mbox{div}u\|^2d\tau+CC_3tP(\mathcal{N}_m(t)).
\end{align}}
and
{\small\begin{align}\label{4.31}
&\v\int_0^t\int(\r\partial_t^{m-1}u_t+\r u\cdot\nabla\partial_t^{m-1}u)\nabla\mbox{div}\partial_t^{m-1}u=-\v\int_0^t\int(\r\partial_t^{m-1}\mbox{div}u_t+\r u\cdot\nabla\partial_t^{m-1}\mbox{div}u)\mbox{div}\partial_t^{m-1}u\nonumber\\
&-\v\int_0^t\int(\nabla\r\cdot\partial_t^{m-1}u_t+\nabla(\r u)^t\cdot\nabla\partial_t^{m-1}u)\mbox{div}\partial_t^{m-1}udxd\tau-\v\int_0^t\int_{\partial\Om}\r(u\cdot\nabla)n\cdot\partial_t^{m-1}u \partial_t^{m-1}\mbox{div}ud\sigma d\tau\nonumber\\
&\leq -\v\int\f{\r}{2}|\partial_t^{m-1}\mbox{div}u(t)|^2dx+\v\int\f{\r_0}{2}|\partial_t^{m-1}\mbox{div}u_0|^2dx
+CtP(\mathcal{N}_m(t))\nonumber\\
&~~+C[1+P(Q(t))]\v\int_0^t\|\partial_t^{m-1}u\|^{\f12}\|\partial_t^{m-1}u\|^{\f12}_{H^1}
\cdot\|\partial_t^{m-1}\mbox{div}u\|^{\f12}_{H^1}\|\partial_t^{m-1}\mbox{div}u\|^{\f12}d\tau\\
&\leq -\v\int\r|\partial_t^{m-1}\mbox{div}u(t)|^2dx+\v\int\r_0|\partial_t^{m-1}\mbox{div}u_0|^2dx
+\f1{16} \v^2 \int_0^t\|\nabla\partial_t^{m-1}\mbox{div}u\|^2d\tau+ CtP(\mathcal{N}_m(t)).\nonumber
\end{align}}
Using \eqref{4.24}, one obtains that
{\small\begin{align}\label{4.34}
\v|\int_0^t\int(\mathcal{C}_1^\a+\mathcal{C}_2^\a)\nabla\mbox{div}\partial_t^{m-1}udxd\tau|
&\leq \f{1}{16}\v^2\int_0^t\|\nabla\partial_t^{m-1}\mbox{div}u\|^2d\tau
+C\int_0^t\|(\mathcal{C}_1^{\a},\mathcal{C}_2^{\a})\|^2d\tau\nonumber\\
&\leq \f{1}{16}\v^2\int_0^t\|\nabla\partial_t^{m-1}\mbox{div}u\|^2d\tau
+CtP(\mathcal{N}_m(t)).
\end{align}}
Substituting \eqref{4.32}-\eqref{4.34} into \eqref{4.30}, one gets that
{\small\begin{align}\label{4.30-1}
&\v\int\f\r2|\partial_t^{m-1}\mbox{div}u(t)|^2dx
-\v R\int_0^t\int\partial_t^{m-1}\nabla(\r\t)\cdot\nabla\mbox{div}\partial_t^{m-1}udxd\tau\nonumber\\
&~~~~~~~~~~~~~+\f34(2\mu+\l)\v^2\int_0^t\|\nabla\partial_t^{m-1}\mbox{div}u\|^2d\tau\leq\v\int\f{\r_0}2|\partial_t^{m-1}\mbox{div}u_0|^2dx+C_3tP(\mathcal{N}_m(t)).
\end{align}}
Finally, it follows from \eqref{2.55} and the Cauchy  inequality that
{\small\begin{align}\label{4.30-11}
&-\v \int_0^t\int\partial_t^{m-1}\nabla(\r\t)\cdot\nabla\mbox{div}\partial_t^{m-1}udxd\tau\nonumber\\
&=-\v \int_0^t\int\Big\{\t\partial_t^{m-1}\nabla\r+[\partial_t^{m-1}\nabla(\r\t)-\t\partial_t^{m-1}\nabla\r]\Big\}\nabla\mbox{div}\partial_t^{m-1}udxd\tau\nonumber\\
&\geq\v\int\f{\t}{2\r}|\partial_t^{m-1}\nabla\r|^2dx-\v\int\f{\t_0}{2\r_0}|\partial_t^{m-1}\nabla\r_0|^2dx-\f18(2\mu+\l)\v^2\int_0^t\|\nabla\partial_t^{m-1}\mbox{div}u\|^2d\tau\nonumber\\
&~~~~-C\int_0^t\|\partial_t^{m-1}\nabla\t\|^2d\tau-C_{m+1}tP(\mathcal{N}_m(t)).
\end{align}}
Substituting \eqref{4.30-11} into \eqref{4.30-1}, one proves \eqref{4.25}. Thus, the proof of Lemma \ref{lem4.3} is completed. $\hfill\Box$

\

Usually, it is hard to obtain the uniform estimate for the term $\int_0^t\|\nabla\partial_t^{m-2}\mbox{div}u\|^2d\tau$ since it involves two times standard space derivatives.  But we observe that $\mbox{div}u$ can be expressed by some good terms   by using the mass conservation law.

\begin{lemma}\label{lem4.5}
For every $m\geq3$, it holds that
{\small\begin{align}
\int_0^t\|\nabla\MZ^{m-2}\mbox{div}u(\tau)\|^2d\tau
&\leq C\int_0^t\|\nabla\partial_t^{m-1}\r(\tau)\|^2d\tau
+C_m tP(\mathcal{N}_m(t)),\label{4.43}\\
\v^2\int_0^t\|\nabla^2\MZ^{m-2}u(\tau)\|^2d\tau
&\leq C_m tP(\mathcal{N}_m(t)),\label{4.43-3}\\
\k(\v)^2\int_0^t\|\nabla\MZ^{m-2}\Delta\t(\tau)\|^2d\tau
&\leq CP(\mathcal{N}_m(t))\int_0^t1+\|\nabla\partial_t^{m-1}(\r,\t)\|^2+\v^2\|\nabla^2u\|^2_{\mathcal{H}^{m-2}}d\tau,
\label{4.43-1}\\
\k(\v)^2\int_0^t\|\partial_t^{m-1}\Delta\t(\tau)\|^2d\tau&\leq CtP(\mathcal{N}_m(t)).\label{4.43-2}
\end{align}}
\end{lemma}
\noindent\textbf{Proof}.  Applying $\nabla\MZ^{\a}$ to \eqref{2.55} with $|\a|\leq m-2$, one has
{\small\begin{equation}\label{4.44}
\nabla\MZ^{\a}\mbox{div}u=-\nabla\MZ^{\a}(\ln\r)_t-\nabla\MZ^{\a}(u_i\partial_{y^i}\ln\r)
-\nabla\MZ^{\a}(u\cdot N\partial_z\ln\r).
\end{equation}}
By using  Proposition \ref{prop3.2}, it is easy to obtain
{\small\begin{equation}\label{4.45}
\int_0^t\|\nabla\MZ^{\a}(\ln\r)_t\|^2+\|\nabla\MZ^{\a}(u_i\partial_{y^i}\ln\r)\|^2d\tau\leq C\int_0^t\|\nabla\partial_t^{m-1}p\|^2d\tau
+CtP(\mathcal{N}_m(t)),
\end{equation}}
and
{\small\begin{align}\label{4.47}
&\int_0^t\|\nabla\MZ^{\a}(u\cdot N\partial_z\ln\r)\|^2d\tau\leq \int_0^t\|\nabla(u\cdot N)\cdot\partial_z\ln\r)\|^2_{\mathcal{H}^{m-2}}d\tau
+\int_0^t\|(\frac{u\cdot N}{\varphi(z)})\cdot Z_3\nabla\ln\r)\|^2_{\mathcal{H}^{m-2}}d\tau\nonumber\\
&\leq  C[P(\mathcal{N}_m(t))+\sup_{0\leq\tau\leq t}\|\frac{u\cdot N}{\varphi(z)}\|^2_{L^\infty}]
\int_0^t\Big(P(\mathcal{N}_m(\tau))+\|\frac{u\cdot N}{\varphi(z)}\|^2_{\mathcal{H}^{m-2}}\Big)d\tau\leq C_mtP(\mathcal{N}_m(t)).
\end{align}}
where the hardy inequality was used  in the last inequality of \eqref{4.47}. Combining \eqref{4.44}-\eqref{4.47}, we obtain \eqref{4.43}.

Note that
{\small\begin{equation}\label{6.10-1}
\Delta =(1+|\nabla\psi|^2)\partial_{zz}+\sum_{i=1,2}\Big(\partial^2_{y^i}-\partial_{y^i}(\partial_i\psi\partial_z)
-\partial_i\psi\partial_z\partial_{y^i}\Big),
\end{equation}}
which, together with $\eqref{2.5}_2$, yields that
{\small\begin{align}
&\v^2\int_0^t\|\nabla^2\MZ^{m-2}u(\tau)\|^2d\tau\leq C\v^2\int_0^t\|\MZ^{m-2}\partial_z^2u(\tau)\|^2d\tau
+ C_m tP(\mathcal{N}_m(t))\nonumber\\
&\leq C\v^2\int_0^t\|\MZ^{m-2}\Delta u(\tau)\|^2d\tau+ C_m tP(\mathcal{N}_m(t)) \nonumber\\
&\leq C\v^2\int_0^t\|\MZ^{m-2}\nabla\mbox{div} u(\tau)\|^2d\tau+ C_m tP(\mathcal{N}_m(t))\leq C_m tP(\mathcal{N}_m(t)).
\end{align}}
Thus, one obtains \eqref{4.43-3}. It follows from \eqref{3.2},  \eqref{4.43},  \eqref{4.43-3} and \eqref{4.13-1}(with $|\a|=m-2$) that
{\small\begin{align}\label{4.47-1}
&\k(\v)^2\int_0^t\|\nabla\MZ^{m-2}\Delta\t\|^2d\tau
\leq C(1+P(Q(t)))\int_0^t\|\nabla\MZ^{m-2}\mbox{div}u\|^2+\|\nabla\partial_t^{m-1}\t(\tau)\|^2+P(\mathcal{N}_m(\tau))d\tau\nonumber\\
&~~~~~~~~~~~~~~~~~~~~~~~~~~~~~~~+C(1+P(Q(t)))\int_0^t\v^2\|\nabla^2u(\tau)\|^2_{\MH^{m-2}}d\tau
+\int_0^t\|\nabla(\mathcal{C}_3^\a,\mathcal{C}_4^\a,\mathcal{C}_5^\a)(\tau)\|^2d\tau\nonumber\\
&\leq  CP(\mathcal{N}_m(t))\int_0^t \|\nabla\partial_t^{m-1}(\r,\t)(\tau)\|^2d\tau
+CtP(\mathcal{N}_m(t)),
\end{align}}
which yields \eqref{4.43-1}.  Finally, it follows from \eqref{3.2} and \eqref{4.13-1} with $\MZ^\a=\partial_t^{m-1}$ that
{\small\begin{equation}\label{4.47-2}
\k(\v)^2\int_0^t\|\partial_t^{m-1}\Delta\t\|^2d\tau
\leq C(1+P(Q(t)))\int_0^tP(\mathcal{N}_m(\tau))d\tau\leq CtP(\mathcal{N}_m(t)),
\end{equation}}
which yields \eqref{4.43-2}. Therefore, the proof of this lemma is completed. $\hfill\Box$

\

Due to the difficulty on the boundary estimates, it is hard to get the uniform estimates on\\
 $\sup_{0\leq\tau\leq t}\|\partial_t^{m-1}\nabla(\r,\t)\|$. However,  the uniform estimate  on $\int_0^t\|\nabla\partial_t^{m-1}(\r,\t)\|^2d\tau$ is possible, which is crucial for us to close the a priori estimates.
\begin{lemma}\label{lem4.8}
It holds, for $m\geq3$,   that
{\small\begin{align}\label{4.94}
&\sup_{0\leq\tau\leq t}\Big(\v\|\nabla\partial_t^{m-2}\mbox{div}u(t)\|^2+\k(\v)\|\partial_t^{m-2}\Delta\t(\tau)\|^2\Big)
+\int_0^t\|(\nabla\partial_t^{m-1}\t,\nabla\partial_t^{m-1}\r)(\tau)\|^2d\tau\nonumber\\
&~~~~~~~\leq C\Big\{\mathcal{N}_m(0)+ tP(\mathcal{N}_m(t))\Big\}.
\end{align}}
\end{lemma}
\noindent\textbf{Proof}.  Multiplying $\eqref{2.26}_1$(with $\MZ^\a=\partial_t^{m-1}$) by $\nabla\partial_t^{m-2}\mbox{div}u$,  one obtains, by using \eqref{3.2},  that
{\small\begin{align}\label{4.94-1}
&\f12(2\mu+\l)\v\|\nabla\partial_t^{m-2}\mbox{div}u(t)\|^2-R\int_0^t\int\nabla\partial_t^{m-1}(\r\t)\cdot\nabla\partial_t^{m-2}\mbox{div}u dxd\tau
\nonumber\\
&\leq \f12(2\mu+\l)\v\|\nabla\partial_t^{m-2}\mbox{div}u_0\|^2+ CtP(\mathcal{N}_m(t)),
\end{align}}
where we have used the following facts
{\small\begin{align}\nonumber
\v\int_0^t\int\nabla\partial_t^{m-1}\mbox{div}u\cdot \nabla\partial_t^{m-2}\mbox{div}udxd\tau=\f12\v\|\nabla\partial_t^{m-2}\mbox{div}u(t)\|^2
-\f12\v\|\nabla\partial_t^{m-2}\mbox{div}u_0\|^2,
\end{align}}
and
{\small\begin{align}
&-\mu\v\int_0^t\int\nabla\times\partial_t^{m-1}\omega\cdot \nabla\partial_t^{m-2}\mbox{div}udxd\tau=
-\mu\v\int_0^t\int n\times\partial_t^{m-1}\omega\cdot \Pi(\nabla\partial_t^{m-2}\mbox{div}u)dxd\tau\nonumber\\
&\leq C\v\int_0^t|\partial_t^{m-1}u|_{H^{\f12}}|\partial_t^{m-2}\mbox{div}u|_{H^{\f12}} d\tau
\leq C\v\int_0^t\|\partial_t^{m-1}u\|_{H^1}\|\partial_t^{m-2}\mbox{div}u\|_{H^1} d\tau\leq  CtP(\mathcal{N}_m(t)).\nonumber
\end{align}}

On the other hand, multiplying \eqref{4.13-1}(with $\MZ^\a=\partial_t^{m-2}$) by $\f{\nabla\partial_t^{m-1}\t}{\t}$  and using \eqref{3.2}, one obtains immediately that
{\small\begin{align}\label{4.94-2}
&\int_0^t\int\f{\r}{\t}|\nabla\partial_t^{m-1}\t|^2dxd\tau+R\int_0^t\int\r\nabla\partial_t^{m-1}\t\cdot\nabla\partial_t^{m-2}\mbox{div}udxd\tau-\k(\v)\int_0^t\int\nabla\partial_t^{m-2}\Delta\t\cdot\f{\nabla\partial_t^{m-1}\t}{\t}dxd\tau
\nonumber\\
&\leq\f18\int_0^t\int\f{\r}{\t}|\nabla\partial_t^{m-1}\t|^2dxd\tau+CP(\mathcal{N}_m(t))\int_0^t\v^2\|\nabla^2u\|^2_{\MH^{m-2}}d\tau+CtP(\mathcal{N}_m(t)).
\end{align}}
Combining \eqref{4.94-1}, \eqref{4.94-2} and using \eqref{4.43-3}, one gets that
{\small\begin{align}\label{4.94-3}
&\f12(2\mu+\l)\v\|\nabla\partial_t^{m-2}\mbox{div}u(t)\|^2+\f78\int_0^t\int\f{\r}{\t}|\nabla\partial_t^{m-1}\t|^2dxd\tau
-\k(\v)\int_0^t\int\nabla\partial_t^{m-2}\Delta\t\cdot\f{\nabla\partial_t^{m-1}\t}{\t}dxd\tau\nonumber\\
&+R\int_0^t\int[\r\nabla\partial_t^{m-1}\t-\nabla\partial_t^{m-1}(\r\t)]\cdot\nabla\partial_t^{m-2}\mbox{div}udxd\tau
\leq \mathcal{N}_m(0)+ CtP(\mathcal{N}_m(t)).
\end{align}}
In order to estimate the  terms on the LHS of \eqref{4.94-3}, we first note that
{\small\begin{equation}\label{4.94-4}
\r\nabla\partial_t^{m-1}\t-\nabla\partial_t^{m-1}(\r\t)=-\t\nabla\partial_t^{m-1}\r-[\partial_t^{m-1},\r]\nabla\t
-[\partial_t^{m-1},\t]\nabla\r,
\end{equation}}
and
{\small\begin{equation}\label{4.94-5}
\nabla\partial_t^{m-2}\mbox{div}u=-\partial_t^{m-2}\nabla(\f{\r_t}{\r})
-\sum_{i=1,2}\partial_t^{m-2}\nabla(\f{u_i\cdot\partial_{y^i}\r}{\r})-\partial_t^{m-2}\nabla(\f{u\cdot N\partial_z\r}{\r}).
\end{equation}}
Then, using \eqref{4.94-4}, \eqref{4.94-5}, \eqref{3.2}, and after some  tedious calculation,  one obtains that
{\small\begin{align}\label{4.94-6}
& R\int_0^t\int[\r\nabla\partial_t^{m-1}\t-\nabla\partial_t^{m-1}(\r\t)]\cdot\nabla\partial_t^{m-2}\mbox{div}udxd\tau\nonumber\\
&\geq\f{7}{8}R\int_0^t\int\f{\t}{\r}|\nabla\partial_t^{m-1}\r|^2dxd\tau-\d_1\int_0^t\|\nabla\partial_t^{m-2}\mbox{div}u\|^2d\tau-C_{\d_1}tP(\mathcal{N}_m(t)).
\end{align}}

It follows from the trace theorem, \eqref{4.29} and  integrating by parts that
{\small\begin{align}\label{4.94-7}
&-\k(\v)\int_0^t\int\nabla\partial_t^{m-2}\Delta\t\cdot\f{\nabla\partial_t^{m-1}\t}{\t}dxd\tau
\nonumber\\
&=\k(\v)\int_0^t\int\partial_t^{m-2}\Delta\t\cdot\Big(\f{\Delta\partial_t^{m-1}\t}{\t}-\f{\nabla\t\cdot\nabla\partial_t^{m-1}\t}{\t^2}\Big)dxd\tau-\nu\k(\v)\int_0^t\int_{\partial\Om}\partial_t^{m-2}\Delta\t\cdot\f{\partial_t^{m-1}\t}{\t}d\sigma d\tau\nonumber\\
&\geq \k(\v)\int\f{|\partial_t^{m-2}\Delta\t|^2}{2\t}dx-\k(\v)\int\f{|\partial_t^{m-2}\Delta\t_0|^2}{2\t_0}dx-CtP(\mathcal{N}_m(t))-\f1{16}\int_0^t\int\f{\r}{\t}|\nabla\partial_t^{m-1}\t|^2dxd\tau\nonumber\\
&~~-C\k(\v)\int_0^t\Big(\|\nabla\partial_t^{m-2}\Delta\t\|^{\f12}\|\Delta\partial_t^{m-2}\t\|^{\f12}+\|\Delta\partial_t^{m-2}\t\|\Big)\Big(\|\nabla\partial_t^{m-1}\t\|^{\f12}\|\partial_t^{m-1}\t\|^{\f12}+\|\partial_t^{m-1}\t\|\Big)d\tau\nonumber\\
&\geq \k(\v)\int\f{|\partial_t^{m-2}\Delta\t|^2}{2\t}dx-\k(\v)\int\f{|\partial_t^{m-2}\Delta\t_0|^2}{2\t_0}dx
-\f18\int_0^t\int\f{\r}{\t}|\nabla\partial_t^{m-1}\t|^2dxd\tau\nonumber\\
&~~~~~~-\d\k(\v)^2\int_0^t\|\nabla\partial_t^{m-2}\Delta\t\|^2d\tau-C_\d tP(\mathcal{N}_m(t)).
\end{align}}
Substituting \eqref{4.94-7} and \eqref{4.94-6} into \eqref{4.94-3}, one obtains immediately that
{\small\begin{align}\label{4.94-8}
&\v\|\nabla\partial_t^{m-2}\mbox{div}u(t)\|^2+\k(\v)\|\partial_t^{m-2}\Delta\t(t)\|^2+\int_0^t\|\nabla\partial_t^{m-1}(\t,\r)\|^2d\tau\nonumber\\
&\leq C\Big\{ \mathcal{N}_m(0)+\d\k(\v)^2\int_0^t\|\nabla\partial_t^{m-2}\Delta\t\|^2d\tau+C\d_1\int_0^t\|\nabla\partial_t^{m-2}\mbox{div}u\|^2d\tau+C_{\d_1,\d} tP(\mathcal{N}_m(t))\Big\}\\
&\leq C\Big\{ \mathcal{N}_m(0)+\d_1\int_0^t\|\nabla\partial_t^{m-1}\r\|^2d\tau+\d P(\mathcal{N}_m(t))\int_0^t\|\nabla\partial_t^{m-1}(\r,\t)\|^2d\tau+C_{\d_1,\d} tP(\mathcal{N}_m(t))\Big\},\nonumber
\end{align}}
where we have used \eqref{4.43}, \eqref{4.43-1} in the last inequality. Setting  $\d=\d_1 P(\mathcal{N}_m(t))^{-1}$ and taking $\d_1$ suitably small, one  proves \eqref{4.94}.
Thus, the proof of Lemma \ref{lem4.8} is completed. $\hfill\Box$

\

Since the estimate in Lemma \ref{lem4.3} is not enough to support us to get the uniform estimate for $\nabla\partial_t^{m-1}u$, so we need some new estimate on $\|\partial_t^{m-1}\mbox{div}u\|$. Fortunately,  we have the following subtle control about $\|\partial_t^{m-1}\mbox{div}u\|$:
\begin{lemma}\label{lem4.7}
	Let us define
{\small	\begin{equation}\label{4.53}
		\Lambda_{m}(t)\triangleq\|(\r, u,\t)\|^2_{\mathcal{H}^m}+\sum_{|\b|\leq m-2}\|\MZ^{\b}\nabla(\r,\t)\|_{H^1_{co}}^2+\sum_{|\b|\leq m-2}\|\MZ^{\b}\nabla u\|_{H^1_{co}}^2.
	\end{equation}}
	Then, for every $m\geq3$,  it holds that
	{\small\begin{eqnarray}\label{4.48}
		\|\partial_t^{m-1}\mbox{div}u(t)\|^2\leq C_2 \Big\{P(\Lambda_{m}(t))
		+P(Q(t))\Big\}.
	\end{eqnarray}}
\end{lemma}
\noindent\textbf{Proof}. Since the proof is the same as the one in \cite{Wang-Xin-Yong},
we omit the details for brevity.
\begin{remark}
	We point out that it does not contain the terms $\|\nabla\partial_t^{m-1}(\r,u,\t)\|$ in the right hand side of \eqref{4.48}. The estimation of $Q(t)$ will be given in section 3.4 below. This key observation allows us  to obtain the uniform estimates for $\|\nabla\partial_t^{m-1}u\|$.
\end{remark}

\subsection{Normal Derivatives Estimates}
\setcounter{equation}{0}

Similar to the corresponding part in \cite{Wang-Xin-Yong}, in order to estimate $\|\nabla u\|_{\mathcal{H}^{m-1}}$, it remains  to estimate $\|\chi\partial_n u\|_{\mathcal{H}^{m-1}}$, where $\chi$ is supported compactly in one of the $\Om_j$ and with value one in a neighborhood of the boundary. Indeed, it follows from the definition of the norm that
$\|\chi\partial_{y^i}u\|_{\mathcal{H}^{m-1}}\leq C\| u\|_{\mathcal{H}^{m}}$ for $i=1,2$.  So it suffices to estimate $\|\chi\partial_n u\|_{\mathcal{H}^{m-1}}$.

Note that
{\small\begin{equation}\label{5.1}
\mbox{div}u=\partial_nu\cdot n+(\Pi\partial_{y^1}u)_1+(\Pi\partial_{y^2}u)_2.
\end{equation}}
and
{\small\begin{equation}\label{5.2}
\partial_nu=[\partial_nu\cdot n] n+\Pi(\partial_nu),
\end{equation}}
where $\Pi$ is  defined in \eqref{3.10-1}.
Thus it follows from \eqref{5.1} and \eqref{5.2} that
{\small\begin{equation}\label{5.3}
	\|\chi\partial_n u\|_{\mathcal{H}^{m-1}}
	\leq \|(\chi\partial_n u\cdot n,\chi\Pi(\partial_nu))\|_{\mathcal{H}^{m-1}}\leq  C_{m}\Big\{\|\chi\mbox{div}u\|_{\mathcal{H}^{m-1}}
	+\|\chi\Pi(\partial_nu)\|_{\mathcal{H}^{m-1}}+ \|u\|_{\mathcal{H}^{m}}\Big\}.
\end{equation}}
Thus it suffices to estimate $\|\chi\Pi(\partial_nu)\|_{\mathcal{H}^{m-1}}$, since  $ \|u\|_{\mathcal{H}^{m}} $ and $\|\chi\mbox{div}u\|_{\mathcal{H}^{m-1}} $ have been estimated in section 3.1 and section 3.2, respectively. We extend the smooth symmetric matrix $A$   in \eqref{1.9} to be
{\small$$A(y,z)=A(y).$$}
Define
{\small\begin{equation}\label{5.4}
\eta\triangleq \chi\Big(\omega\times n+\Pi(Bu)\Big)=\chi\Big(\Pi(\omega\times n)+\Pi(Bu)\Big).
\end{equation}}
The $\eta$ defined here, which enables one to avoid to estimate $\nabla^2p$, is different slightly from the one in \cite{Masmoudi-R}. Then in view of the Navier-slip boundary condition \eqref{1.9},  $\eta$ satisfies:
{\small\begin{equation}\label{5.5}
\eta|_{\partial\Om}=0.
\end{equation}}
Since $\omega\times n=(\nabla u-(\nabla u)^t)\cdot n$, so $\eta$ can be rewritten as
{\small\begin{equation}\label{5.6}
\eta=\chi\Big\{\Pi(\partial_nu)-\Pi(\nabla(u\cdot n))+\Pi((\nabla n)^t\cdot u)+\Pi(Bu)\Big\},
\end{equation}}
which  yields immediately that
{\small\begin{equation}\label{5.7}
\|\chi\Pi(\partial_n u)\|_{\mathcal{H}^{m-1}}\leq C_{m+1}(\|\eta\|_{\mathcal{H}^{m-1}}+\|u\|_{\mathcal{H}^{m}}).
\end{equation}}
Hence, it remains  to estimate $\|\eta\|_{\mathcal{H}^{m-1}}$. In fact, one can get the following conormal estimates for $\eta$:

\begin{lemma}\label{lem5.1}
For every $m\geq3$, it holds that
{\small\begin{equation}\label{5.8}
\sup_{0\leq \tau\leq t}\|\eta(\tau)\|^2_{\mathcal{H}^{m-1}}+\v\int_0^t\|\nabla\eta\|^2_{\mathcal{H}^{m-1}}d\tau\leq CC_{m+2}\Big\{P(\mathcal{N}_m(0))+\d\v^2\int_0^t\|\nabla^2u\|^2_{\mathcal{H}^{m-1}}d\tau+C_\d t P(\mathcal{N}_m(t))\Big\}.
\end{equation}}
\end{lemma}
\noindent\textbf{Proof}.  Notice that
{\small$$\nabla\times((u\cdot\nabla) u)=(u\cdot\nabla)\omega-(\omega\cdot\nabla)u+\mbox{div}u \cdot\omega,$$}
so $\omega$ solves the following   vorticity equation
{\small\begin{eqnarray}\label{5.9}
\r \omega_t+\r(u\cdot\nabla)\omega=\mu\v \Delta\omega+F_1,
\end{eqnarray}}
with
{\small\begin{equation}\label{5.10}
F_1\triangleq -\nabla\r\times u_t-\nabla\r\times(u\cdot\nabla)u+\r(\omega\cdot\nabla)u-\r\mbox{div}u \omega.
\end{equation}}
Consequently, we obtain that $\eta$ solves the equation
{\small\begin{align}\label{5.11}
&\r\eta_t+\r u_1\partial_{y^1}\eta+\r u_2\partial_{y^2}\eta+\r u\cdot N\partial_z\eta-\mu\v\Delta\eta\nonumber\\
&~~~=\chi[F_1\times n+\Pi(B F_2)]+\chi(F_3+F_4)+F_5+\v\Delta(\Pi B)\cdot u \triangleq F,
\end{align}}
where
{\small\begin{align}\label{5.12}
\begin{cases}
 F_2=(2\mu+\l)\v\nabla\mbox{div}u-\nabla p, \\
 F_3=-2\mu\sum_{j=1}^3\v\partial_j\omega\times\partial_jn-\mu\v\omega\times\Delta n+\sum_{i=1}^2\r u_i\omega\times\partial_{y^i}n+\r u\cdot N\omega\times\partial_zn \\
~~~~~~-[\sum_{i=1}^2\r u_i\Pi(\partial_{y^i}B\cdot u)+\r u\cdot N\Pi(\partial_zB\cdot u)]+\mu\sum_{j=1}^3\v\Pi(\partial_jB\partial_ju),\\
 F_4=-\sum_{i=1}^2\r u_i(\partial_{y^i}\Pi)(Bu)-\r u\cdot N(\partial_z\Pi)(Bu),\\
 F_5=\sum_{i=1}^2\r u_i\partial_{y^i}\chi\cdot (\omega\times n+\Pi(Bu))+\r u\cdot N\partial_z\chi\cdot (\omega\times n+\Pi(Bu))\\
~~~~~~~~~~-2\mu\sum_{j=1}^3\v\partial_j\chi\partial_j(\omega\times n+\Pi(Bu))+\v\mu\Delta\chi\cdot (\omega\times n+\Pi(Bu)).
\end{cases}
\end{align}}

Let us begin with the proof of the $L^2$-energy estimate. Multiplying \eqref{5.11} by $\eta$, yields that
{\small\begin{equation}\label{5.16}
\sup_{0\leq \tau\leq t}\int\r|\eta|^2dx+2\v\int_0^t\|\nabla\eta\|^2d\tau
\leq \int\r_0|\eta_0|^2dx+\int_0^t\int F\eta dxd\tau.
\end{equation}}
To estimate the terms on the RHS, note that
{\small\begin{align}
&\int_0^t\|\chi\Pi(F_1\times n)\|^2_{\mathcal{H}^{m-1}}d\tau\leq C_mP(\mathcal{N}_m(t))
\int_0^t\|\nabla\partial_t^{m-1}(\r,\t)(\tau)\|^2d\tau+CP(\mathcal{N}_m(t)),\label{5.17}\\
&\int_0^t\|\chi\Pi(BF_2)\|^2_{\mathcal{H}^{m-1}}d\tau\leq
C_{m+1}\Big\{\int_0^tP(\mathcal{N}_m)+\|\nabla\partial_t^{m-1}(\r,\t)\|^2d\tau
+\v^2\int_0^t\|\chi\nabla\mbox{div}u\|^2_{\mathcal{H}^{m-1}}d\tau\Big\},\label{5.18}\\
&\int_0^t\|\chi F_3\|^2_{\mathcal{H}^{m-1}}d\tau\leq
C_{m+2}\Big\{\v^2\int_0^t\|\chi\nabla^2u\|^2_{\mathcal{H}^{m-1}}d\tau
+tP(\mathcal{N}_m(t))\Big\},\label{5.19}\\
&\int_0^t\|\chi F_4\|^2_{\mathcal{H}^{m-1}}d\tau\leq C_{m+2}tP(\mathcal{N}_m(t)).\label{5.20}
\end{align}}
Since all the terms in $F_5$ are supported away from the boundary, one can estimate all the derivatives by the $\|\cdot\|_{\mathcal{H}^{m}}$ norms. Therefore, it is easy to obtain
{\small\begin{equation}\label{5.21}
\int_0^t\| F_5\|^2_{\mathcal{H}^{m-1}}d\tau\leq C_{m+1}\Big\{\v^2\int_0^t\|\chi\nabla^2 u\|^2_{\mathcal{H}^{m-1}}d\tau+tP(\mathcal{N}_m(t)) \Big\}.
\end{equation}}
Finally, by integrating by parts, it is easy to obtain, for $|\a|\leq m-1$,  that
{\small\begin{equation}\label{5.21-1}
\int_0^t\int\v\MZ^{\a}(\Delta(\Pi B)\cdot u)\cdot \MZ^{\a}\eta d\tau\leq \d\v\int_0^t \|\nabla\MZ\eta\|^2d\tau +C_{m+2}tP(\mathcal{N}_m(t)).
\end{equation}}
Consequently, substituting these estimates into \eqref{5.16} and using the Cauchy inequality,  one has that
{\small\begin{align}\label{5.22}
\sup_{0\leq \tau\leq t}\int\r|\eta|^2dx+2\v\int_0^t\|\nabla\eta\|^2d\tau &\leq \int\r_0|\eta_0|^2dx+C_{m+2}\Big\{\d\int_0^t\|\nabla\partial_t^{m-1}(\r,\t)(\tau)\|^2d\tau\nonumber\\
&~~~~~~~+\d\v^2\int_0^t\|\chi \nabla^2 u\|^2_{\mathcal{H}^{m-1}}d\tau+C_\d t P(\mathcal{N}_m(t))\Big\}.
\end{align}}
Thus,  we proved \eqref{5.8} for $k=0$ by using Lemma \ref{lem4.8}.

To prove the general case, let us assume the \eqref{5.8} is proved for $k\leq m-2$.  We apply $\MZ^\a$
to \eqref{5.11} for $|\a|=m-1$ to obtain that
{\small\begin{equation}\label{5.23}
\r \MZ^\a\eta+\r (u\cdot\nabla )\MZ^\a\eta-\mu\v\MZ^\a\Delta\eta= \MZ^\a F+\tilde{\mathcal{C}}_1^\a+\tilde{\mathcal{C}}_2^\a,
\end{equation}}
where
{\small\begin{align}\label{5.24}
\begin{cases}
\tilde{\mathcal{C}}_1^\a=-[\MZ^\a,\r]\eta_t=\sum_{|\b|\geq 1, \b+\g=\a}C_{\a,\b}\MZ^\b\r \MZ^\g\eta_t,\\[2mm]
\tilde{\mathcal{C}}_2^\a=-\sum_{|\b|\geq 1, \b+\g=\a}C_{\a,\b}\MZ^\b(\r u_i)\MZ^\g\partial_{y^i}\eta
-\sum_{|\b|\geq 1, \b+\g=\a}C_{\a,\b}\MZ^\b(\r u\cdot N)\MZ^\g\partial_{z}\eta\\
~~~~~~~~~~-\r (u\cdot N)\sum_{|\b|\leq m-2}C_{\b}(z) \partial_z\MZ^\b\eta, 
\end{cases}
\end{align}}
where $C_\b(z)$ is bounded smooth function of $z$.
Multiplying \eqref{5.23} by $\MZ^\a\eta$ and using \eqref{5.17}-\eqref{5.21}, one obtains that
{\small\begin{align}\label{5.26}
&\sup_{0\leq \tau\leq t}\int\f12\r|\MZ^\a\eta|^2dx
\leq\mu\v\int_0^t\int\MZ^\a\Delta\eta\MZ^\a\eta dxd\tau +\int\f12\r_0|\MZ^\a\eta_0|^2dx+\int_0^t\int(\tilde{\mathcal{C}}_1^\a+\tilde{\mathcal{C}}_2^\a)\MZ^\a\eta dxd\tau\nonumber\\
&~~ ~~~~~~~~~+C_{m+1}\Big\{\d\v^2\int_0^t\|\chi \nabla^2 u\|^2_{\mathcal{H}^{m-1}}d\tau+\d\int_0^t\|\nabla\partial_t^{m-1}(\r,\t)\|^2d\tau+C_\d t P(\mathcal{N}_m(t))\Big\}.
\end{align}}
By the same argument as Lemma 3.12 of \cite{Wang-Xin-Yong}, one can gets that
{\small\begin{align}\label{5.29}
\v \int_0^t\int\MZ^\a\Delta\eta\MZ^\a\eta dxd\tau&\leq -\f34\mu\v \int_0^t\|\nabla\MZ^\a\eta\|^2d\tau+C\v\int_0^t\|\nabla\eta\|^2_{\mathcal{H}^{m-2}}d\tau
+C_{m+2}t P(\mathcal{N}_m(t)),
\end{align}}
and
{\small\begin{eqnarray}\label{5.34}
\int_0^t\|(\tilde{\mathcal{C}}_1^\a,\tilde{\mathcal{C}}_2^\a)\|^2d\tau
\leq C_{m+2}[1+P(Q(t))]\int_0^t P(\mathcal{N}_m(\tau))d\tau\leq C_{m+2}t P(\mathcal{N}_m(t)).
\end{eqnarray}}
 Substituting \eqref{5.29} and \eqref{5.34} into \eqref{5.26} and using Lemma \ref{lem4.8}, one obtains that
{\small\begin{align}\label{5.35}
&\sup_{0\leq \tau\leq t}\int\r|\MZ^\a\eta|^2dx+ \mu\v \int_0^t\|\nabla\MZ^\a\eta\|^2d\tau
\leq C_{m+2}\Big\{\int\r_0|\MZ^\a\eta_0|^2dx+\v\int_0^t\|\nabla\eta\|^2_{\mathcal{H}^{m-2}}d\tau\nonumber\\
&~~~~~~~~~~~~~~~~~~~~~~~~~~~~~~~~~~~~~~~~~~~~~~~~~~~~~~~~~~~+\d\v^2\int_0^t\| \nabla^2 u\|^2_{\mathcal{H}^{m-1}}d\tau+C_\d t P(\mathcal{N}_m(t))\Big\}.
\end{align}}
By using the induction assumption, one can eliminate the term $ \v\int_0^t\|\nabla\eta\|^2_{\mathcal{H}^{m-2}}d\tau$. Therefore,  the proof Lemma \ref{lem5.1} is completed. $\hfill\Box$

\

From \eqref{5.1}, \eqref{5.2} and \eqref{5.7}, it holds that
{\small\begin{align}
&\sum_{|\b|\leq m-2}\|\MZ^\b\nabla u\|^2_{H^1_{co}}\leq C_{m+1}\Big(\|u\|^2_{\mathcal{H}^m}+\|\eta\|^2_{\mathcal{H}^{m-1}}
+\sum_{k=0}^{m-2}\|\partial_t^{k}\mbox{div}u(t)\|^2_{m-1-k}\Big),\label{5.36}\\
&\int_0^t\|\nabla^2 u\|^2_{\mathcal{H}^{m-1}}d\tau\leq C_{m+2}\int_0^t\Big(\|\nabla u\|^2_{\mathcal{H}^m}
+\|\nabla\eta\|^2_{\mathcal{H}^{m-1}}+\|\nabla\mbox{div}u\|^2_{\mathcal{H}^{m-1}}
+P(\mathcal{N}_m) \Big)d\tau,\label{5.37},
\end{align}
\begin{align}
&\sum_{k=0}^{m-2}\int_0^t\|\nabla^2\partial_t^k u\|^2_{m-1-k}d\tau\nonumber\\
&~~~~~~\leq C_{m+2}\Big\{\int_0^t \|\nabla u\|^2_{\mathcal{H}^m}
+\|\nabla\eta\|^2_{\mathcal{H}^{m-1}}d\tau+\sum_{k=0}^{m-2}\int_0^t \|\partial_t^{k}\nabla\mbox{div}u\|^2_{m-1-k}d\tau
+tP(\mathcal{N}_m) \Big\},\label{5.37-1}\\
&\v\int_0^t\|\nabla^2\MZ^{m-2}u\|^2d\tau\leq C_{m+1}\v\int_0^t\|\nabla\eta\|^2_{\mathcal{H}^{m-2}}d\tau+ C_{m+1}tP(\mathcal{N}_m),\label{5.38}
\end{align}}
where \eqref{4.43} are used  in the estimate of \eqref{5.38}. Setting $\d=\d_1 P(\mathcal{N}_m(t))^{-1}$ and taking $\d_1$ suitably small, then it follows  from \eqref{2.16-1}, \eqref{4.6}, \eqref{4.25}, \eqref{4.94},  \eqref{5.8},  \eqref{5.36}-\eqref{5.38} and Lemma \ref{lem4.5}  that
{\small\begin{align}\label{5.39}
&\sup_{0\leq \tau\leq t}\Big\{\Lambda_{m}(\tau)+\|\eta(\tau)\|^2_{\mathcal{H}^{m-1}}
+\v\|\partial_t^{m-1}(\nabla\r,\mbox{div}u)(\tau)\|^2+\v\|\nabla\partial_t^{m-2}\mbox{div}u(\tau)\|^2+\k(\v)\|\partial_t^{m-2}\Delta\t(\tau)\|^2\Big\}\nonumber\\
&+\int_0^t\|\nabla\partial_t^{m-1}(\r,\t)(\tau)\|^2d\tau+\int_0^t\v\|\nabla u(\tau)\|^2_{\mathcal{H}^{m}}+\k(\v)\|\nabla\t(\tau)\|^2_{\mathcal{H}^{m}}d\tau
+\v\int_0^t\sum_{k=0}^{m-2}\|\partial_t^{k}\nabla^2u(\tau)\|^2_{m-1-k}d\tau
\nonumber\\
&+\k(\v)\int_0^t\sum_{k=0}^{m-2}\|\partial_t^{k}\Delta\t(\tau)\|^2_{m-1-k}d\tau+\int_0^t\v^2\|\nabla^2\partial_t^{m-1}u(\tau)\|^2+\k(\v)^2\|\partial_t^{m-1}\Delta\t(\tau)\|^2d\tau\nonumber\\
&\leq CC_{m+2}\Big\{\mathcal{N}_{m}(0) +  tP(\mathcal{N}_{m}(t))\Big\},~~\mbox{for}~m\geq 3.
\end{align}}

\subsection{$L^\infty$-Estimates}
\setcounter{equation}{0}

This section is devoted to  estimate   $L^{{\infty}}$-norm parts contained in \eqref{1.11-1}. However, it is not easy to get such estimates because the equations of $\r,u,\t$ are strongly coupled and the viscosity and  heat conductivity are not at the same order. Actually, if one estimate $\|\nabla\r\|_{\MH^{1,\infty}}$ directly, then the  high order term $\int_0^t\|\nabla\t\|_{\mathcal{H}^{2,\infty}}d\tau$ must appear, and it is very hard to control such a term. To overcome the difficulty, instead, we try to estimate the $L^{{\infty}}$-norm of $\nabla(\r\t)$,$\nabla u$ and $\nabla\t$.  Firstly, we have the following useful Lemma:

\begin{lemma}\label{lem6.1}
For every $|\a|\geq 0$, it holds that
{\small\begin{align}
&\|\MZ^\a(\r, \t,u)\|^2_{L^\infty}\leq CP(\Lambda_{m}(t)),~~\mbox{for}~~m\geq 2+|\a|,\label{6.1}\\[2mm]
&\|\nabla\r\|^2_{\mathcal{H}^{1,\infty}}\leq C_3\Big(\|\nabla(\r\t)\|^2_{\mathcal{H}^{1,\infty}}+\|\nabla\t\|^2_{\mathcal{H}^{1,\infty}}\Big)\cdot P(\Lambda_m(t)),~\mbox{for}~m\geq 5,\label{6.1-1}\\[2mm]
& Q(t)\leq C_3\Big\{\|\nabla u\|^2_{\mathcal{H}^{1,\infty}}+\v\|\partial_z^2 u\|^2_{L^{\infty}}+\|\nabla(\r\t)\|^4_{\mathcal{H}^{1,\infty}}+\|\nabla\t\|^4_{\mathcal{H}^{1,\infty}}+P(\Lambda_m(t))\Big\}
~~\mbox{for}~~m\geq5.\label{6.4-3}\\[2mm]
&\|\mbox{div}u\|^2_{\mathcal{H}^{1,\infty}}\leq C_3P(\Lambda_{m}(t))[1+\|\nabla(\r\t)\|^2_{\mathcal{H}^{1,\infty}}+\|\nabla\t\|^2_{\mathcal{H}^{1,\infty}}]
,~~\mbox{for}~~m\geq 4,\label{6.4}\\[2mm]
&\|\nabla\mbox{div}u\|^2_{L^\infty}\leq
\begin{cases}
C_3 P(Q(t)),\\
C_3[P(\|\nabla\r\|^2_{\mathcal{H}^{1,\infty}})+P(\|\nabla u\|^2_{\mathcal{H}^{1,\infty}})+P(\Lambda_{m})],
\end{cases}
~~\mbox{for}~m\geq3,\label{6.4-1}\\
&\|\nabla\mbox{div}u\|^2_{\mathcal{H}^{1,\infty}}\leq C\Big\{\|\nabla(\r\t)\|^2_{\mathcal{H}^{2,\infty}}+\|\nabla\t\|^2_{\mathcal{H}^{2,\infty}}
+P(\Lambda_{m}(t))+P(Q(t))\Big\},
~~\mbox{for}~~m\geq6,\label{6.4-2}
\end{align}}
\end{lemma}

\noindent\textbf{Proof}. The proof of \eqref{6.1} is an immediately consequence of \eqref{3.3}, we omit the details here. Notice that
{\small\begin{equation}\label{6.9}
\nabla\r=\f{\nabla(\r\t)}{\t}-\f{\r}{\t}\nabla\t,
\end{equation}}
which immediately implies \eqref{6.1-1}.  Using \eqref{6.1} and \eqref{6.1-1},  we immediately obtain \eqref{6.4-3}.

Using \eqref{2.55}, \eqref{6.1}, \eqref{6.1-1}, $\eqref{1.9-1}_1$ and the facts that
{\small\begin{align}\label{2.55-1}
\mbox{div}u&=\t^{-1}[\t_t+(u\cdot\nabla)\t]-p^{-1}[p_t+(u\cdot\nabla)p],\\
\nabla\mbox{div}u&=\t^{-1}[\nabla\t_t+(u\cdot\nabla)\nabla\t]-p^{-1}[\nabla p_t+(u\cdot\nabla)\nabla p]-p^{-1}\nabla p\mbox{div}u\nonumber\\
&~~~~+Rp^{-1}\nabla\r[\t_t+(u\cdot\nabla)\t]+\t^{-1}\nabla{u}\cdot\nabla\t-p^{-1}\nabla{u}\cdot\nabla{p}.
\end{align}}
it is easy to prove \eqref{6.4}-\eqref{6.4-2}. For simplicity, we omit the details here.  Therefore, we complete the proof of lemma \ref{lem6.1}. $\hfill\Box$

\begin{remark}
Lemma \ref{lem6.1} implies that one needs only to estimate  $\|\nabla(\r\t)\|^2_{\mathcal{H}^{1,\infty}}$, $\v\|\nabla(\r\t)\|^2_{\mathcal{H}^{2,\infty}}$, $\|\nabla\t\|^2_{\mathcal{H}^{1,\infty}}$, $\|\nabla u\|^2_{\MH^{1,\infty}}$ and $\v\|\partial_z^2 u\|^2_{L^{\infty}}$. Indeed, it is crucial to estimate $\|\nabla(\r\t)\|_{\mathcal{H}^{1,\infty}}$ but not $\|\nabla\r\|_{\mathcal{H}^{1,\infty}}$.  We also point out that the condition \eqref{1.8} will be used repeatedly to  control the possible interaction between the viscous  and the thermal boundary layers in the following analysis.
\end{remark}

\noindent{\bf \underline{Uniform Estimate for $\|\nabla p\|^2_{\MH^{1,\infty}}$ and $\v\|\nabla p\|^2_{\MH^{2,\infty}}$}}:\\[2mm]
Firstly, we have the following lemma which  will be used to prove \eqref{6.115} below.
\begin{lemma}\label{lem3.17}
	Assume  \eqref{1.8} holds, one has, for  $m\geq6$,  that
{\small	\begin{align}\label{6.101}
	&\v\int_0^t\|\partial_{zz}u\|_{\MH^{1,\infty}}d\tau\leq C_{m+2}\Big\{\Lambda_{m}(0) +  tP(\mathcal{N}_{m}(t))\Big\}.
	\end{align}}
\end{lemma}
\noindent\textbf{Proof}.  Then, it follows from    the momentum equations $\eqref{1.1}_2$,   \eqref{3.3} \eqref{6.10-1}, \eqref{5.39} and \eqref{6.4-2} that
{\small\begin{align}\label{6.102}
&\v\int_0^t\|\partial_{zz}u\|_{\MH^{1,\infty}}d\tau\leq C_3\int_0^tP(\mathcal{N}_m) +\v\|\nabla\mbox{div}u\|_{\MH^{1,\infty}}+\v\|\nabla u\|_{\MH^{2,\infty}}d\tau\nonumber\\
&\leq C_3tP(\mathcal{N}_m)+C_3\v\int_0^t\|\nabla (\t,u)\|_{\MH^{2,\infty}}d\tau\leq C_3tP(\mathcal{N}_m)+C_3\v\int_0^t\|\nabla^2\MZ^2 (\t,u)\|^{\f12}_{H^1_{co}}\|\nabla\MZ^2 (\t,u)\|^{\f12}_{H^2_{co}}d\tau\nonumber\\
&\leq C_3tP(\mathcal{N}_m)+C_3\v^4\int_0^t\|(\nabla^2u,\Delta\t)\|^2_{\MH^3}d\tau\leq CC_{m+2}\Big\{\mathcal{N}_{m}(0) +  tP(\mathcal{N}_{m}(t))\Big\},
\end{align}}
where we have used the \eqref{1.8} in the last inequality.  Thus, the proof of Lemma \ref{lem3.17} is completed.  $\hfill\Box$

\

\begin{lemma}\label{lem3.18}
	Assume  \eqref{1.8} holds, one has, for  $m\geq6$,  that
{\small	\begin{align}\label{6.115}
\|\nabla p\|^2_{\MH^{1,\infty}}+\v\|\nabla p\|^2_{\MH^{2,\infty}}+\int_0^t\|\nabla p\|^2_{\MH^{2,\infty}}d\tau\leq C_{m+2}\Big\{P(\mathcal{N}_{m}(0)) +  tP(\mathcal{N}_{m}(t))\Big\}.
	\end{align}}	
\end{lemma}
\noindent\textbf{Proof}.  Note that
{\small\begin{align}\label{6.116}
&\|\nabla p\|^2_{\MH^{1,\infty}}+\v\|\nabla p\|^2_{\MH^{2,\infty}}+\int_0^t\|\nabla p\|^2_{\MH^{2,\infty}}d\tau\nonumber\\
&\leq CC_{m+2}\Big\{ P(\Lambda_{m}(t))+tP(\mathcal{N}_m(t))+\v\|\partial_n p\|^2_{\MH^{2,\infty}}+\int_0^t\|\partial_n p\|^2_{\MH^{2,\infty}}d\tau\Big\},
\end{align}}
therefore, one needs only to control $\v\|\partial_n p\|^2_{\MH^{2,\infty}}+\int_0^t\|\partial_n p\|^2_{\MH^{2,\infty}}d\tau$. Substituting \eqref{2.55} into $\eqref{2.5}_2$, one can obtain that
{\small\begin{align}\label{6.117}
(2\mu+\l)\v[\nabla\r_t+(u\cdot\nabla)\nabla\r]+\r\nabla p=-\r^2\dot{u}-\mu\v\r\nabla\times\omega-(2\mu+\l)\v\r[\nabla(\f1\r)\r_t+\nabla(\f{u}{\r})\cdot\nabla\r],
\end{align}}
where $\dot{u}=u_t+(u\cdot\nabla)u$.  It follows from \eqref{6.9} that
{\small\begin{equation}\label{6.118}
\nabla\r_t=\f1\t \nabla(\r\t)_t-\f1\t\Big(\f1\t\t_t\nabla(\r\t)+\r_t\nabla\t+\r\nabla\t_t-\f\r\t \nabla\t \t_t\Big),
\end{equation}}
and
{\small\begin{equation}\label{6.119}
(u\cdot\nabla)\nabla\r=\f1\t (u\cdot\nabla)\nabla(\r\t)-\f1\t\Big(\nabla(\r\t)(\f{u}{\t}\cdot\nabla)\t+\nabla\t(u\cdot\nabla)\r+\r(u\cdot\nabla)\nabla\t-\f\r\t \nabla\t (u\cdot\nabla)\t\Big).
\end{equation}}
Substituting \eqref{6.118} and \eqref{6.119} into \eqref{6.117}, one can obtain that
{\small\begin{align}\label{6.120}
&(2\mu+\l)\v[\nabla p_t+(u\cdot\nabla)\nabla p]+p\nabla p\nonumber\\
&=-p\r\dot{u}-\mu\v p \nabla\times\omega-(2\mu+\l)\v\Big\{ p[\nabla(\f1\r)\r_t+\nabla(\f{u}{\r})\cdot\nabla\r]-\Big(\f1\t\t_t\nabla p+R\r_t\nabla\t+R\r\nabla\t_t-R\f\r\t \nabla\t \t_t\Big)\nonumber\\
&~-\Big(\nabla p(\f{u}{\t}\cdot\nabla)\t+R\nabla\t(u\cdot\nabla)\r+R\r(u\cdot\nabla)\nabla\t-R\f\r\t \nabla\t (u\cdot\nabla)\t\Big)\Big\}\triangleq p\r\dot{u}-\mu\v p \nabla\times\omega+I.
\end{align}}
Then, it follows from \eqref{6.120} that
{\small\begin{align}\label{6.121}
&(2\mu+\l)\v[\partial_n p_t+(u\cdot\nabla)\partial_n p]+p\partial_n p\nonumber\\
&=-\r p\dot{u}\cdot n-\mu\v p  (n \cdot\nabla\times\omega)+I\cdot n+(2\mu+\l)\v \nabla p\cdot(u\cdot\nabla)n\triangleq J.
\end{align}}

Define
{\small\begin{equation}\label{6.123}
h\triangleq \MZ^\a\partial_np,
\end{equation}}
then, applying $\MZ^\a$ with $|\a|\leq2$, one can obtains that
{\small\begin{align}\label{6.122}
&(2\mu+\l)\v[h_t+(u\cdot\nabla)h]+p h=\MZ^\a J-(2\mu+\l)\v[\MZ^\a,u\cdot\nabla]\partial_np-[\MZ^\a,p]\partial_np\triangleq K.
\end{align}}
It is convenient to consider the above equation in the Lagrangian coordinates
{\small\begin{equation}\label{6.124}
\tilde{h}(t,\xi)=h(t,X(t,\xi)),~~\tilde{p}(t,\xi)=p(t,X(t,\xi)),~~\tilde{K}(t,\xi)=K(t,X(t,\xi)),
\end{equation}}
where
{\small\begin{align*}
\begin{cases}
\f{dX(t,\xi)}{dt}=u(t,X(t,\xi)),\\
X(0,\xi)=\xi\in\Om.
\end{cases}
\end{align*}}
Then \eqref{6.122} is rewritten as
{\small\begin{align}\label{6.125}
\f{d}{dt}\tilde{h}+\f{\tilde{p}}{(2\mu+\l)\v } \tilde{h}=\f{1}{(2\mu+\l)\v } \tilde{K},
\end{align}}
which yields immediately the following solution formula
{\small\begin{equation}\label{6.126}
\tilde h(t,\xi)=\tilde h(0,\xi)\exp{(-\int_0^t\f{\tilde{p}(\tau,\xi)}{(2\mu+\l)\v}d\tau)}
+\f{1}{(2\mu+\l)\v}\int_0^t\tilde{K}(\tau,\xi)\exp{(-\int_\tau^t\f{\tilde{p}(s,\xi)}{(2\mu+\l)\v}ds)}d\tau.
\end{equation}}
Notice that  $\f{\tilde{p}}{2\mu+\l}\geq c>0$ with $c$ independent of $\v$, then this together with \eqref{6.126} yield  that
{\small\begin{align}\label{6.127}
\|\tilde h(t)\|_{L^\i}\leq \|\tilde h(0)\|_{L^\i}e^{-\f{ct}{\v}}+C\int_0^t\f{1}{\v}\|\tilde{K}(\tau)\|_{L^\i}e^{-c\v^{-1}(t-\tau)}d\tau.
\end{align}}
It follows from \eqref{6.127} and   Holder inequality that
{\small\begin{align}
\|\tilde h(t)\|_{L^\i}&\leq \|\tilde h(0)\|_{L^\i}e^{-\f{ct}{\v}}+C\left(\int_0^t\|\tilde{K}(\tau)\|^2_{L^\i}d\tau\right)^{\f12}
\left( \int_0^t\f1{\v^2}\exp(-\f{2c(t-\tau)}{\v})d\tau\right)^{\f12}\nonumber\\
&\leq \|\tilde h(0)\|_{L^\i}e^{-\f{ct}{\v}}+C\f{1}{\sqrt{\v}}\left(\int_0^t\|\tilde{K}(\tau)\|^2_{L^\i}d\tau\right)^{\f12},\nonumber
\end{align}}
i.e.
{\small\begin{align}\label{6.129}
\v\|\tilde h(t)\|^2_{L^\i}&\leq \v\|\tilde h(0)\|^2_{L^\i}+C\int_0^t\|\tilde{K}(\tau)\|^2_{L^\i}d\tau.
\end{align}}

On the other hand, integrating \eqref{6.127} over $[0,t]$, one can obtains that
{\small\begin{align}\label{6.130}
\int_0^t\|\tilde{h}(\tau)\|^2_{L^\i}d\tau&\leq C\v\|\tilde h(0)\|^2_{L^\i}+C\int_0^t\left(\int_0^\tau\f{1}{\v}\|\tilde{K}(s)\|_{L^\i}e^{-c\v^{-1}(\tau-s)}ds\right)^2d\tau\nonumber\\
&\leq C\v\|\tilde h(0)\|^2_{L^\i}+C\int_0^t\left(\int_{\mathbb{R}}\f{1}{\v}\|\tilde{K}(s)\|_{L^\i}I_{(0,t)}(s)e^{-c\v^{-1}(\tau-s)}I_{(0,t)}(\tau-s)ds\right)^2d\tau\nonumber\\
&\leq C\v\|\tilde h(0)\|^2_{L^\i}+C\int_0^t\left|\Big(\|\tilde{K}(\cdot)\|_{L^\i}I_{(0,t)}(\cdot)\Big)\ast\Big(\f{1}{\v}e^{-c\v^{-1}(\cdot)}I_{(0,t)}(\cdot)\Big)\right|^2d\tau\nonumber\\
&\leq C\v\|\tilde h(0)\|^2_{L^\i}+C\left\|\Big(\|\tilde{K}(\cdot)\|_{L^\i}I_{(0,t)}(\cdot)\Big)\right\|^2_{L^2}
\cdot\left\| \Big(\f{1}{\v}e^{-c\v^{-1}(\cdot)}I_{(0,t)}(\cdot)\Big)\right\|^2_{L^1}\nonumber\\
&\leq C\v\|\tilde h(0)\|^2_{L^\i}+C\int_0^t\|\tilde{K}(\tau)\|^2_{L^\i}d\tau.
\end{align}}
Then, it follows from \eqref{6.130} and \eqref{6.129} that
{\small\begin{align}\label{6.131}
\v\|\tilde h(t)\|^2_{L^\i}+\int_0^t\|\tilde{h}(\tau)\|^2_{L^\i}d\tau&\leq C \v\|h(0)\|^2_{L^\i}+C\int_0^t\|\tilde K(\tau)\|^2_{L^\i}d\tau.
\end{align}}
In order to close the estimates of \eqref{6.131}, one still needs to control the second term on the RHS of \eqref{6.131}.
Direct calculation implies that
{\small\begin{align}\label{6.132}
(\nabla\times\omega)\cdot N=-\partial_{y^1}\omega_2+\partial_{y^2}\omega_1+\partial_1\psi\cdot\partial_{y^2}\omega_3
-\partial_2\psi\cdot\partial_{y^1}\omega_3,
\end{align}}
then, it follows from  \eqref{5.39}, \eqref{6.132}, \eqref{6.9} and \eqref{6.1} that
{\small\begin{align}\label{6.133}
&\int_0^t\|J\|^2_{\MH^{2,\i}}d\tau\leq C\int_0^tP(\mathcal{N}_m(\tau))
\cdot\Big(1+\v^2\|\nabla p\|^2_{\MH^{2,\infty}}+\v^2\|(\nabla\times\omega)\cdot N\|^2_{\MH^{2,\infty}}\nonumber\\
&~~~~~~~~~~~~~~~~~~~~~~~~~~~~~~~~~~~~~~~~~~~~~~~~+\v^2\|\nabla\t\|^2_{\MH^{3,\infty}}+\v^2\|\nabla u\|^2_{\MH^{2,\infty}}\Big)d\tau\\
&\leq C\int_0^tP(\mathcal{N}_m(\tau))
\cdot\Big(1+\v^2\|\nabla\t\|^2_{\MH^{3,\infty}}+\v^2\|\nabla u\|^2_{\MH^{3,\infty}}\Big)d\tau\nonumber\\
&\leq \v^4\int_0^t\|\Delta\t\|^2_{\MH^4}d\tau+\v^2\int_0^t\|\nabla^2u\|^2_{\MH^4}d\tau+CtP(\mathcal{N}_m(t))\leq CC_{m+2}\Big\{\mathcal{N}_m(0)+tP(\mathcal{N}_m(t))\Big\},\nonumber
\end{align}}
where we have used the condition \eqref{1.8} in the last inequality.  It follows from \eqref{5.39} that
{\small\begin{align}\label{6.134}
&\int_0^t\|(2\mu+\l)\v[\MZ^\a,u\cdot\nabla]\partial_np\|^2_{L^\i}+\|[\MZ^\a,p]\partial_np\|^2_{L^\i}d\tau\nonumber\\
&\leq C\v^2\int_0^t\|\nabla^2u\|^2_{\MH^4}d\tau+CtP(\mathcal{N}_m(t))\leq CC_{m+2}\Big\{\mathcal{N}_m(0)+tP(\mathcal{N}_m(t))\Big\}
\end{align}}
which, together with \eqref{6.133}, yields that
{\small\begin{align}\label{6.136}
\int_0^t\|\tilde K(\tau)\|^2_{L^\i}d\tau\leq CC_{m+2}\Big\{\mathcal{N}_m(0)+tP(\mathcal{N}_m(t))\Big\}.
\end{align}}
Combining \eqref{6.136}, \eqref{6.131}, \eqref{6.124}, \eqref{6.123}, \eqref{6.116} and \eqref{5.39}, one proves \eqref{6.115}. Thus, the proof of Lemma \ref{lem3.18} is completed. $\hfill\Box$

\

\noindent{\bf \underline{Estimates for  $\|\nabla u\|^2_{\mathcal{H}^{1,\infty}}$ and $\v\|\nabla u\|^2_{\mathcal{H}^{1,\infty}}$}}:
\begin{lemma}\label{lem6.2}
Assume  \eqref{1.8} holds, one has, for  $m\geq6$,  that
{\small\begin{align}\label{6.2}
&\|\nabla u\|^2_{\mathcal{H}^{1,\infty}} \leq  C_{m+2}\Big\{\mathcal{N}_{m}(0) +P(\|\nabla\t\|^2_{\MH^{1,\infty}}) +  tP(\mathcal{N}_{m}(t))\Big\}.
\end{align}}
\end{lemma}
\noindent\textbf{Proof}.
The proof of this lemma is similar to the corresponding lemma of \cite{Wang-Williams} except the pressure is a function of density and temperature.  Away from the boundary, one clearly has by the classical isotropic Sobolev embedding theorem that
{\small\begin{equation}\label{6.3}
\|\chi\MZ\nabla u\|^2_{L^\infty}+\|\chi\nabla u\|^2_{L^\infty}\leq C\|u\|^2_{\mathcal{H}^{m}}\leq \Lambda_{m}(t), ~~\mbox{for}~~m\geq 4,
\end{equation}}
where the support of  $\chi$ is away from the boundary. Therefore, by using a partition of unity subordinated to the covering \eqref{2.0}, we only need to estimate $\|\chi_j\MZ\nabla u\|_{L^\infty}+\|\chi_j\nabla u\|_{L^\infty}$ for $j\geq 1$. For notational convenience, we shall denote $\chi_j$ by $\chi$. Similar to \cite{Masmoudi-R}, we use the local parametrization in the neighborhood of the boundary given by a normal geodesic system which makes the Laplacian has a convenient form. Let us denote
{\small\begin{equation*}\label{6.5}
\Psi^n(y,z)=\left(\begin{array}{cccc} &y\\&\psi(y)\end{array}\right)
-z n(y)=x,
\end{equation*}}
where
{\small\begin{equation*}\label{6.6}
n(y)=\f{1}{\sqrt{1+|\nabla\psi(y)|^2}}\left(\begin{array}{cccc} &\partial_1\psi(y)\\&\partial_1\psi(y)\\&-1\end{array}\right),
\end{equation*}}
is the unit outward normal. As before, we can extend $n$ and $\Pi$ in the interior by setting
{\small$$n(\Psi^n(y,z))=n(y),~~\Pi(\Psi^n(y,z))=\Pi(y).$$}
Note that $n(y,z)$ and $\Pi(y,z)$ have different definitions from the ones used before. The interest of this parametrization is that in the associated local basis $(e_{y^1},e_{y^2}, e_z)$ of $\mathbb{R}^3$, it holds that $\partial_z=\partial_n$ and
{\small$$\Big(e_{y^i}\Big)\Big|_{\Psi^n(y,z)}\cdot\Big(e_{z}\Big)\Big|_{\Psi^n(y,z)}=0.$$}
The scalar product on $\mathbb{R}^3$ induces in this coordinate system the Riemannian metric $g$ under the form
{\small\begin{equation*}\label{6.7}
g(y,z)=\left(\begin{array}{cccc} &\tilde{g}(y,z)&0\\&0&1 \end{array}\right).
\end{equation*}}
Therefore, the Laplacian in this coordinate system reads
{\small\begin{equation}\label{6.8}
\Delta f=\partial_{zz}f+\f12\partial_z(\ln|g|)\partial_zf+\Delta_{\tilde{g}}f
\end{equation}}
where $|g|$ denotes the determinant of the matrix $g$, and $\Delta_{\tilde{g}}$  is given by
{\small$$\Delta_{\tilde{g}}f=\f{1}{\sqrt{|\tilde{g}|}} \sum_{i,j=1,2}\partial_{y^i}(\tilde{g}_{ij}|\tilde{g}|^\f12\partial_{y^j}f),$$}
which involves only the tangential derivatives.

It follows from \eqref{5.1}($n$ and $\Pi$  in the coordinate system we have just defined) and Lemma \ref{lem6.1},  for $m\geq 5$,  that
{\small\begin{eqnarray}\label{6.17}
&&\|\chi\nabla u\|^2_{L^\infty}+\|\chi\MZ\nabla u\|^2_{L^\infty}\leq C_2 \Big(\|\chi\Pi\partial_nu\|^2_{L^\infty}+\|\MZ(\chi\Pi\partial_nu)\|^2_{L^\infty}+\|\chi\mbox{div}u\|^2_{L^\infty}\nonumber\\
&&~~~~~~~~~~~~~~~~~~~~~~~~~~~~~+\|\MZ\mbox{div}u\|^2_{L^\infty}+\|\MZ Z_yu\|^2_{L^\infty}+\|Z_yu\|^2_{L^\infty}\Big)\nonumber\\
&&\leq C_3\Big\{\|\chi\Pi\partial_nu\|^2_{L^\infty}+\|\MZ(\chi\Pi\partial_nu)\|^2_{L^\infty}
+P(\Lambda_{m})+P(\|(\nabla(\r\t),\nabla\t)\|^2_{\MH^{1,\infty}})\Big\}.
\end{eqnarray}}
Consequently, one needs only to estimate $\|\chi\Pi\partial_nu\|^2_{L^\infty}+\|\MZ(\chi\Pi\partial_nu)\|^2_{L^\infty}$. To estimate these quantity, it is useful to use the vorticity $\omega$.  Indeed,  we have
{\small\begin{equation*}\label{6.18}
\Pi(\omega\times n)=\Pi((\nabla u-\nabla u^t)\cdot n)=\Pi(\partial_nu-\nabla(u\cdot n)-\nabla n^t\cdot u).
\end{equation*}}
Therefore, one obtains that
{\small\begin{equation}\label{6.19}
\|\chi\Pi\partial_nu\|^2_{L^\infty}+\|\MZ(\chi\Pi\partial_nu)\|^2_{L^\infty}
\leq C_3\Big\{\|\chi\Pi(\omega\times n)\|^2_{L^\infty}+\|\MZ(\chi\Pi(\omega\times n))\|^2_{L^\infty}
+\Lambda_{m}(t)\Big\},
\end{equation}}
which yields that we only need to estimate $\|\chi\Pi(\omega\times n)\|^2_{L^\infty}$ and $\|\MZ(\chi\Pi(\omega\times n))\|^2_{L^\infty}$.

\

By setting in the support of $\chi$
{\small\begin{equation*}\label{6.20}
\tilde{\omega}(y,z)=\omega(\Psi^n(y,z)),~~(\tilde{\r}, \tilde{u})(y,z)=(\r, u)(\Psi^n(y,z)).
\end{equation*}}
Then it follows from \eqref{5.9} and  \eqref{6.8} that
{\small\begin{equation} \label{6.21}
\tilde{\r} \tilde{\omega}_t+\tilde{\r} \tilde{u}^1\partial_{y^1}\tilde{\omega}+\tilde{\r} \tilde{u}^2\partial_{y^2}\tilde{\omega}+\tilde{\r} \tilde{u}\cdot n\partial_{z}\tilde{\omega}=\mu\v(\partial_{zz}\tilde{\omega}+\f12\partial_z(\ln|g|)\partial_z\tilde{\omega}
+\Delta_{\tilde{g}}\tilde{\omega})+\tilde{F}_1,
\end{equation}}
and
{\small\begin{equation} \label{6.22}
\tilde{\r} \tilde{u}_t+\tilde{\r} \tilde{u}^1\partial_{y^1}\tilde{u}+\tilde{\r} \tilde{u}^2\partial_{y^2}\tilde{u}+\tilde{\r} \tilde{u}\cdot n\partial_{z}\tilde{u}=\mu\v(\partial_{zz}\tilde{u}+\f12\partial_z(\ln|g|)\partial_z\tilde{u}
+\Delta_{\tilde{g}}\tilde{u})+\tilde{F}_2,
\end{equation}}
where
{\small\begin{equation*} \label{6.23}
\tilde{F}_1(y,z)=F_1(\Psi^n(y,z)),~~\tilde{F}_2(y,z)=F_2(\Psi^n(y,z)),
\end{equation*}}
where $F_1$ and $F_2$ are defined in \eqref{5.10} and \eqref{5.12}, respectively.  Note that we use the same convention as before for a vector $u$, and $u^j$ denotes the components of $u$ in the local basis $(e_{y^1}, e_{y^2}, e_z)$ just defined in this section, whereas $u_j$ denotes its components in the standard basis of $\mathbb{R}^3$. The vectorial equation of \eqref{6.21} and \eqref{6.22}
have to be understood component by component in the standard basis of $\mathbb{R}^3$.

Similar to \eqref{5.6}, we define that
{\small\begin{equation} \label{6.24}
\tilde{\eta}(y,z)=\chi\Big(\tilde\omega\times n+\Pi(B\tilde{u})\Big),
\end{equation}}
where $A$ is extended into the interior domain by $B(y,z)=B(y)$. Thus, from the boundary conditions \eqref{1.9}, one has
{\small\begin{equation} \label{6.25}
\tilde{\eta}(y,0)=0.
\end{equation}}
By using \eqref{6.21} and \eqref{6.22},  $\tilde\eta$ solves the equations
{\small\begin{align} \label{6.26}
&\tilde{\r} \tilde{\eta}_t+\tilde{\r} \tilde{u}^1\partial_{y^1}\tilde{\eta}+\tilde{\r} \tilde{u}^2\partial_{y^2}\tilde{\eta}+\tilde{\r} \tilde{u}\cdot n\partial_{z}\tilde{\eta}\nonumber\\
& ~~=\mu\v(\partial_{zz}\tilde{\eta}+\f12\partial_z(\ln|g|)\partial_z\tilde{\eta})
+\chi(\tilde{F}_1\times n)+\chi\Pi(B\tilde{F}_2)+F_\chi+\chi F_\kappa,
\end{align}}
where the source terms are given by
{\small\begin{align}
\begin{cases}\label{6.27}
F_\chi= \Big(\tilde{\r} \tilde{u}^1\partial_{y^1} +\tilde{\r} \tilde{u}^2\partial_{y^2} +\tilde{\r} \tilde{u}\cdot n\partial_{z}\Big)\chi\cdot(\tilde\omega\times n+\Pi(B\tilde{u}))\\
~~~~~~~~~~-\mu\v\Big(\partial_{zz}\chi+2\partial_z\chi\partial_z+\f12\partial_z(\ln|g|)\cdot\partial_z\chi\Big)
\cdot(\tilde\omega\times n+\Pi(B\tilde{u})), \\
F_\kappa= \Big(\tilde{\r} \tilde{u}^1\partial_{y^1}\Pi +\tilde{\r} \tilde{u}^2\partial_{y^2}\Pi\Big) \cdot(B\tilde{u})+ \tilde{\omega}\times\Big(\tilde{\r} \tilde{u}^1\partial_{y^1}n +\tilde{\r} \tilde{u}^2\partial_{y^2}n\Big)\\
~~~~~~~~~~+\Pi\Big((\tilde{\r} \tilde{u}^1\partial_{y^1} +\tilde{\r} \tilde{u}^2\partial_{y^2} )B\cdot\tilde{u}\Big)+\mu\v\Delta_{\tilde{g}}\tilde{\omega}\times n+\mu\v\Pi(B\Delta_{\tilde{g}}\tilde{u})
\end{cases}
\end{align}}
Note that in the calculating of the source terms, and in particular $F_\kappa$ which contains all the commutators coming from the fact that $n$ and $\Pi$ are not constant, we have used the idea that in the coordinate system that we have just defined, $B$, $n$ and $\Pi$ do not depend on the normal variable. By using that $\Delta_{\tilde{g}}$ involves only the tangential derivatives and that the derivatives of $\chi$ are compactly supported away from the boundary, one obtains the following estimates, for $m\geq6$,
{\small\begin{align}\label{6.29}
\begin{cases}
\|\chi\Pi(F_1\times n)\|^2_{\mathcal{H}^{1,\infty}}\leq C_2 P(\mathcal{N}_m(t)), \\
\|F_\chi\|^2_{\mathcal{H}^{1,\infty}}
\leq C_3\Big(\|u\|^2_{\mathcal{H}^{1,\infty}}\cdot\|u\|^2_{\mathcal{H}^{2,\infty}}
+\v^2\|u\|^2_{\mathcal{H}^{3,\infty}}\Big) \leq C_3P(\mathcal{N}_m(t)),\\ 
\|\chi F_\kappa\|^2_{\mathcal{H}^{1,\infty}}
\leq C_4\Big\{\|u\|^4_{\mathcal{H}^{1,\infty}}+\|u\|^2_{\mathcal{H}^{1,\infty}}\|\nabla u\|^2_{\mathcal{H}^{1,\infty}}
+\v^2(\|u\|^2_{\mathcal{H}^{3,\infty}}+\|\nabla u\|^2_{\mathcal{H}^{3,\infty}})\Big\}\\ 
~~~~~~~~~~~~~~~\leq C_4\Big\{ P(\mathcal{N}_m(t))
+\v^2\|\nabla^2 u\|^2_{\mathcal{H}^{4}}\Big\}, 
\end{cases}
\end{align}}
and from \eqref{6.4-2}, it holds that
{\small\begin{equation}\label{6.32}
\|\chi\Pi(B\tilde{F}_2)\|^2_{\mathcal{H}^{1,\infty}}\leq C_3\Big\{\v^2\|\nabla\mbox{div}u\|^2_{\mathcal{H}^{1,\infty}}+\|\nabla p\|^2_{\mathcal{H}^{1,\infty}}\Big\}\leq C_4\Big\{P(\mathcal{N}_m(t))+\v^2\|(\nabla\t,\nabla(\r\t))\|^2_{\MH^{2,\infty}}\Big\}.
\end{equation}}
Consequently, it follows  from \eqref{6.29}-\eqref{6.32},  for $m\geq 6$, that
{\small\begin{equation}\label{6.33}
\|\tilde{F}\|^2_{\mathcal{H}^{1,\infty}}\leq C_4\Big\{\v^2\|\nabla^2 u\|^2_{\mathcal{H}^{4}}+\v^2\|(\nabla\t,\nabla(\r\t))\|^2_{\MH^{2,\infty}}+ P(\mathcal{N}_m(t))\Big\},
\end{equation}}
where $\tilde{F}=\chi(\tilde{F}_1\times n)+\chi\Pi(B\tilde{F}_2)+F_\chi+\chi F_\kappa$.

In order to eliminate the term $\f12\partial_z(\ln|g|)\partial_z\tilde{\eta}$, one defines
{\small\begin{equation}\label{6.34}
\tilde\eta=\f{1}{|g|^\f14}\eta=\bar{\g} \bar{\eta}.
\end{equation}}
Note that
{\small\begin{equation}\label{6.35}
\|\tilde\eta\|_{\mathcal{H}^{1,\infty}} \leq C_3\|\bar{\eta}\|_{\mathcal{H}^{1,\infty}},~~\mbox{and}
~~\|\bar{\eta}\|_{\mathcal{H}^{1,\infty}} \leq C_3\|\tilde\eta\|_{\mathcal{H}^{1,\infty}},
\end{equation}}
and $\bar\eta$ solves the equations
{\small\begin{eqnarray} \label{6.36}
&&\tilde{\r} \bar{\eta}_t+\tilde{\r} \tilde{u}^1\partial_{y^1}\bar{\eta}+\tilde{\r} \tilde{u}^2\partial_{y^2}\bar{\eta}+\tilde{\r} \tilde{u}\cdot n\partial_{z}\bar{\eta}-\v\partial_{zz}\bar{\eta}\nonumber\\
&& ~~=\f1{\bar\g}\Big(\tilde{F}+\v\partial_{zz}\bar\g\cdot \bar\eta+\f12\v\partial_z(\ln|g|)\partial_z\g\cdot\bar\eta
-\tilde\r(\tilde{u}\cdot\nabla \bar\g)\bar\eta\Big)\triangleq \mathcal{S}.
\end{eqnarray}}
It is difficult to directly obtain the  explicit solution formula of \eqref{6.36}, so one rewrites it as
{\small\begin{align} \label{6.37}
&\tilde{\r}(t,y,0)\Big[\bar{\eta}_t+ \tilde{u}^1(t,y,0)\partial_{y^1}\bar{\eta}+ \tilde{u}^2(t,y,0)\partial_{y^2}\bar{\eta}+z\partial_z(\tilde{u}\cdot n)(t,y,0)\partial_{z}\bar{\eta}\Big]-\v\partial_{zz}\bar{\eta}\nonumber\\
&=\mathcal{S}+[\tilde{\r}(t,y,0)-\r(t,y,z)]\bar\eta_t
+\sum_{i=1,2}[(\tilde{\r}\tilde{u}^i)(t,y,0)-(\r\tilde{u}^i)(t,y,z)]\partial_{y^i}\bar\eta\nonumber\\
&~~~~~~-\r(t,y,z)[(\tilde{u}\cdot n)(t,y,z)-z\partial_z(\tilde{u}\cdot n)(t,y,0)]\partial_z\bar\eta\nonumber\\[2mm]
&~~~~~~+[\r(t,y,z)-\r(t,y,0)]\cdot z\partial_z(\tilde{u}\cdot n)(t,y,0)\partial_z\bar\eta\triangleq M
~~\mbox{for}~~z>0,
\end{align}}
with the boundary condition $\bar\eta(t,y,0)=0$. By using Lemma \ref{lemA.2} in Appendix, one has that
{\small\begin{align}\label{6.38}
\|\bar\eta\|_{\mathcal{H}^{1,\infty}} &\lesssim \|\bar\eta_0\|_{\mathcal{H}^{1,\infty}}
+\int_0^t\|\tilde{\r}^{-1}\|_{L^\infty} \|M\|_{\mathcal{H}^{1,\infty}}d\tau+\int_0^t(1+\|\tilde{\r}^{-1}\|_{L^\infty})(1+\|\MZ(\r,u,\nabla u)\|^2_{L^\infty}) \|\bar\eta\|_{\mathcal{H}^{1,\infty}}d\tau\nonumber\\
&\lesssim \|\bar\eta_0\|_{\mathcal{H}^{1,\infty}}
+C\int_0^t\|M\|_{\mathcal{H}^{1,\infty}}d\tau+CtP(\mathcal{N}_m(t)).
\end{align}}
It remains to estimate the right hand side of \eqref{6.38}. Firstly, by using \eqref{6.33}, one has that
{\small\begin{equation}\label{6.39}
\|\mathcal{S}\|^2_{\mathcal{H}^{1,\infty}}\leq C_4\Big\{\v^2\|\nabla^2 u\|^2_{\mathcal{H}^{4}}+\v^2\|\nabla\t\|^2_{\MH^{2,\infty}}+ CP(\mathcal{N}_m(t)) \Big\},~\mbox{for}~m\geq 6.
\end{equation}}
Next,   using the Taylor formula and the fact that $\bar\eta$ is compactly supported in $z$ and by the same argument of Lemma 3.14 in \cite{Wang-Xin-Yong},  one can obtain, for $m\geq 5$,  that
{\small\begin{align}\label{6.42}
&\|[\tilde{\r}(t,y,0)-\r(t,y,z)]\bar\eta_t\|^2_{\mathcal{H}^{1,\infty}}\|[(\tilde{\r}\tilde{u}^1)(t,y,0)-(\tilde{\r}\tilde{u}^1)(t,y,z)]\partial_{y^1}\bar\eta\|^2_{\mathcal{H}^{1,\infty}}\nonumber\\
&+\|[(\tilde{\r}\tilde{u}^2)(t,y,0)-(\tilde{\r}\tilde{u}^2)(t,y,z)]\partial_{y^2}\bar\eta\|^2_{\mathcal{H}^{1,\infty}}
\|[\r(t,y,z)-\r(t,y,0)]\cdot z\partial_z(\tilde{u}\cdot n)(t,y,0)\partial_z\bar\eta\|^2_{\mathcal{H}^{1,\infty}}\nonumber\nonumber\\
&+\|\r(t,y,z)[(\tilde{u}\cdot n)(t,y,z)-z\partial_z(\tilde{u}\cdot n)(t,y,0)]\partial_z\bar\eta\|^2_{\mathcal{H}^{1,\infty}}\leq CP(\mathcal{N}_m(t)),
\end{align}}
Then, it follows from  \eqref{6.39} and \eqref{6.42} that
{\small\begin{equation} \label{6.49}
\|M\|^2_{\mathcal{H}^{1,\infty}}\leq C_4\Big\{\v^2\|\nabla^2 u\|^2_{\mathcal{H}^{4}}+\v^2\|\nabla\t\|^2_{\MH^{2,\infty}}+P(\mathcal{N}_m(t)) \Big\}, ~~\mbox{for}~~ m\geq 6.
\end{equation}}
Substituting \eqref{6.49} into \eqref{6.38}, we have, for $m\geq 6$, that
{\small\begin{align}\label{6.50}
\|\bar\eta\|^2_{\mathcal{H}^{1,\infty}}
&\lesssim \|\bar\eta_0\|^2_{\mathcal{H}^{1,\infty}}
+C_4tP(\mathcal{N}_m(t))+C_4t\int_0^t\v^2\|\nabla^2 u\|^2_{\mathcal{H}^{4}}+\v^2\|\nabla\t\|^2_{\MH^{2,\infty}}d\tau\nonumber\\
&\leq CC_{m+2}\Big\{\mathcal{N}_m(0)+tP(\mathcal{N}_m(t))\Big\},
\end{align}}
where we have used \eqref{3.3}, \eqref{1.8} and the  Holder inequality   in the last inequality.
Therefore, combining   \eqref{6.50}, \eqref{6.35}, \eqref{6.24}, \eqref{6.19}, \eqref{6.17}, \eqref{6.3} and \eqref{5.39}, one obtains  \eqref{6.2}.  Therefore, he proof of Lemma \ref{lem6.2} is completed.
$\hfill\Box$

\begin{lemma}\label{lem3.16}
Assume  \eqref{1.8} holds, one has, for  $m\geq6$,  that
{\small\begin{align}\label{6.99}
&\v\|\partial_{zz}u\|^2_{L^\infty}\leq  CC_{m+2}\Big\{P(\mathcal{N}_{m}(0)) +P(\|\nabla\t\|^2_{\MH^{1,\infty}}) +  tP(\mathcal{N}_{m}(t))\Big\}.
\end{align}}
\end{lemma}
\noindent\textbf{Proof}.  By similar argument as the one in Lemma \ref{lem6.2}, one firstly has that
{\small\begin{align}\label{6.100}
\v\|\partial_{zz}u\|^2_{L^\infty}\leq C_2\Big\{P(\Lambda_{m})+P(\|(\nabla u,\nabla\t,\nabla(\r\t))\|^2_{\MH^{1,\infty}})+\v\|\partial_z\bar{\eta}\|^2_{L^\infty}\Big\},
\end{align}}
Thus, one needs only to estimate $\v\|\partial_z\bar{\eta}\|^2_{L^\infty}$. one rewrites \eqref{6.36} as
{\small\begin{equation} \label{6.36-2}
\bar{\eta}_t-\v\partial_{zz}\bar{\eta}=-(\tilde\r-1) \bar{\eta}_t-\tilde{\r} \tilde{u}^1\partial_{y^1}\bar{\eta}-\tilde{\r} \tilde{u}^2\partial_{y^2}\bar{\eta}-\tilde{\r} \tilde{u}\cdot n\partial_{z}\bar{\eta} +\mathcal{S}=:\Xi,
\end{equation}}
where $\bar{\eta}$  satisfying the homogenous Dirichlet boundary condition $\bar{\eta}|_{z=0}=0$. Then, $\bar{\eta}$ has the following expression:
{\small\begin{align*}
\bar{\eta}(t,y,z)=\int_0^{+\i} G(t,z,z')\eta_0(y,z')dz'+\int_0^t\int_0^{+\i}G(t-\tau,z,z')\Xi(\tau,y,z')dz'd\tau
\end{align*}}
where
{\small\begin{align*}
G(t,z,z')=\f{1}{\sqrt{4\pi\mu\v t}}\Big[\exp(-\f{|z-z'|^2}{4\m\v t})-\exp(-\f{|z+z'|^2}{4\m\v t})\Big].
\end{align*}}
Then, one can obtains  that
{\small\begin{align*}
\sqrt{\v}\partial_z\bar{\eta}(t,y,z)=\sqrt{\v}\int_0^{+\i}\partial_z G(t,z,z')\eta_0(y,z')dz'+\sqrt{\v}\int_0^t\int_0^{+\i}\partial_zG(t-\tau,z,z')\Xi(\tau,y,z')dz'd\tau.
\end{align*}}
Since $\eta_0(y,z)$ vanishes on the boundary due to the compatibility condition,  it follows from the integrating by parts to the first term that
{\small\begin{align}\label{6.103}
\sqrt{\v}\|\partial_z\bar{\eta}\|_{L^\infty}\leq \sqrt{\v}\|\partial_z\eta_0\|_{L^\i}+\int_0^t\f{1}{\sqrt{t-\tau}}\|\Xi(\tau)\|_{L^\i}d\tau.
\end{align}}
Directly calculation shows that
{\small\begin{align}\label{6.104}
\Big(\int_0^t\f{1}{\sqrt{t-\tau}}\|\Xi(\tau)\|_{L^\i}d\tau\Big)^2
\leq C\v^2\int_0^t \|\nabla^2u\|^2_{\MH^3}+CtP(\mathcal{N}_m(t)),~~\mbox{for}~~m\geq5.
\end{align}}
Substituting \eqref{6.103}, \eqref{6.104} into \eqref{6.100} and using \eqref{5.39}, \eqref{6.2}, \eqref{6.115},  one proves \eqref{6.99}.  Therefore, the proof of this lemma is completed. $\hfill\Box$

\

\noindent{\bf \underline{Estimate for $\|\nabla\t\|_{\MH^{1,\infty}}$}}:\\[2mm]
In order to estimate $\|\nabla\t\|_{\MH^{1,\infty}}$, the most difficult part is to control the term  $p\nabla\mbox{div}u$ which comes from  the term $p\mbox{div}u$ appearing on LHS the energy equation $\eqref{2.5}_3$. Actually, if $p\nabla\mbox{div}u$ is regarded as the source term,  it is very difficult to bound  the term  $\int_{0}^t\|p\nabla\mbox{div}u(\tau)\|^2_{\MH^{1,\i}}d\tau$  since the derivative is  too higher. It is noted that such difficulty does not arise in the isentropic case \cite{Wang-Xin-Yong}. So, to overcome the difficulty, new idea is needed.  Fortunately, we find that the term  $p\nabla\mbox{div}u$ can be  decomposed into two parts i.e. $\nabla(\r\t)_t$ and $\nabla\t_t$.  The most difficult term   $\nabla\t_t$ can be absorbed into the main part  of  equation,  while $\nabla(\r\t)_t$ is regarded  as the source term which has already be controlled in Lemma \ref{lem3.18}. This  observation is key to close the pointwise estimates.

\begin{lemma}\label{lem3.19}
Assume  \eqref{1.8} holds, one has, for  $m\geq6$,  that
{\small\begin{align}\label{6.108}
\|\nabla\t\|^2_{\MH^{1,\infty}}\leq C C_{m+2}\Big\{ P(\mathcal{N}_m(0))  +  tP(\mathcal{N}_{m}(t))\Big\}.
\end{align}}
\end{lemma}
\noindent\textbf{Proof}. Due to \eqref{6.1}, one needs only to estimate $\|\partial_z\t\|_{\MH^{1,\i}}$ or $\|\partial_n\t\|_{\MH^{1,\i}}$. It follows from $\eqref{2.5}_3$ that
{\small\begin{align}\label{6.105}
&\r[\nabla\t_t+(u\cdot\nabla)\nabla\t]+p\nabla\mbox{div}u-\k(\v)\Delta\nabla\t\nonumber\\
&=-\nabla p\mbox{div}u-[\nabla\r(\t_t+u\cdot\nabla\t)+\r\nabla{u}\cdot\nabla\t]+\v\nabla[2\m|Su|^2+\l|\mbox{div}u|^2].
\end{align}}
In order to deal with the term $p\nabla{div}u$, using the mass equation $\eqref{2.5}_1$,  one first notice that
{\small\begin{equation}\label{6.106}
p\mbox{div}u=R\r[\t_t+(u\cdot\nabla)\t]-[p_t+(u\cdot\nabla)p],
\end{equation}}
and
{\small\begin{align}\label{6.107}
p\nabla\mbox{div}u&=R\r[\nabla\t_t+(u\cdot\nabla)\nabla\t]-[\nabla p_t+(u\cdot\nabla)\nabla p]-\nabla p\mbox{div}u\nonumber\\
&~~~~+R\nabla\r[\t_t+(u\cdot\nabla)\t]+R\r\nabla{u}\cdot\nabla\t-\nabla{u}\cdot\nabla{p}.
\end{align}}
Then, it follows from $\eqref{2.5}_3$ and \eqref{6.105}-\eqref{6.107} that
{\small\begin{align}\label{6.109}
&(R+1)\r[\nabla\t_t+(u\cdot\nabla)\nabla\t]-\k(\v)\Delta\nabla\t\nonumber\\
&=[\nabla p_t+(u\cdot\nabla)\nabla p]+\nabla{u}\cdot\nabla p-(1+R)[\nabla\r(\t_t+u\cdot\nabla\t)+\r\nabla{u}\cdot\nabla\t]\nonumber\\
&~~~~+\v\nabla[2\m|Su|^2+\l|\mbox{div}u|^2]\triangleq B_1,
\end{align}}
and
{\small\begin{align}\label{6.110}
(R+1)\r[\t_t+(u\cdot\nabla)\t]-\k(\v)\Delta\t=[p_t+(u\cdot\nabla)p]+\v[2\m|Su|^2+\l|\mbox{div}u|^2]\triangleq B_2.
\end{align}}

We   use the local coordinates $(y,z)$ defined in Lemma \ref{lem6.2} in the neighborhood of the boundary  which makes the Laplacian has a convenient form. The functions $\tilde{\r}, \tilde{u}, \chi, n, \Pi$ are the same   ones defined in Lemma \ref{lem6.2}. By setting in the support of $\chi$
{\small\begin{align}
\tilde{\t}(y,z,t)=\t(\Psi^n(y,z),t), ~~~\widetilde{\nabla\t}(y,z,t)=(\nabla\t)(\Psi^n(y,z),t)
\end{align}}
Then it follows from \eqref{6.109} and \eqref{6.110} that
{\small\begin{equation}\label{6.109-1}
(R+1)[\tilde\r\widetilde{\nabla\t}_t+\tilde\r\tilde{u}^1\partial_{y^1}\widetilde{\nabla\t}+\tilde\r\tilde{u}^2\partial_{y^2}\widetilde{\nabla\t}+\tilde\r\tilde{u}\cdot n\partial_{z}\widetilde{\nabla\t}]
=\k(\v)(\partial_{zz}\widetilde{\nabla\t}+\f12\partial_z(\ln|g|)\partial_z\widetilde{\nabla\t}
+\Delta_{\tilde{g}}\widetilde{\nabla\t})+\tilde{B}_1,
\end{equation}}
and
{\small\begin{equation}\label{6.110-1}
(R+1)[\tilde\r\tilde\t_t+\tilde{\r}\tilde{u}^1\partial_{y^1}\tilde\t+\tilde{\r}\tilde{u}^2\partial_{y^2}\tilde\t
+\tilde{\r}\tilde{u}\cdot n\partial_{z}\tilde\t]=\k(\v)(\partial_{zz}\tilde\t+\f12\partial_z(\ln|g|)\partial_z\tilde\t
+\Delta_{\tilde{g}}\tilde\t) +\tilde{B}_2.
\end{equation}}
where
{\small\begin{align}
\tilde{B}_1 =B_1(\Psi^n(y,z),t),~~~\tilde{B}_2 =B_2(\Psi^n(y,z),t).
\end{align}}
Define
{\small\begin{equation}\label{6.112}
\zeta(y,z,t)=\chi(n\cdot\widetilde{\nabla\t}-\nu\tilde\t),
\end{equation}}
then, in view of the boundary condition $\eqref{1.9}_3$, $\zeta$ satisfies $\zeta=0~\mbox{on}~\partial\Om$.
Considering  $\eqref{6.109-1}\cdot n+\n\cdot\eqref{6.110-1}$, it is easy to know that $\zeta$ satisfies
{\small\begin{align} \label{6.111}
&\tilde{\r} \zeta_t+\tilde{\r} \tilde{u}^1\partial_{y^1}\zeta+\tilde{\r} \tilde{u}^2\partial_{y^2}\zeta+\tilde{\r} \tilde{u}\cdot n\partial_{z}\zeta\nonumber\\
& =\k(\v)\big(\partial_{zz}\zeta+\f12\partial_z(\ln|g|)\partial_z\zeta\big)
+\chi(\tilde{B}_1\cdot n)+\nu\chi\tilde{B}_2+F_\chi^{\t}+\chi F_\kappa^{\t},
\end{align}}
where the source terms are given by
{\small\begin{align}\label{6.27-1}
\begin{cases}
F_\chi^\t=(R+1) \Big(\tilde{\r} \tilde{u}^1\partial_{y^1} +\tilde{\r} \tilde{u}^2\partial_{y^2} +\tilde{\r} \tilde{u}\cdot n\partial_{z}\Big)\chi\cdot(n\cdot\widetilde{\nabla\t}-\nu\tilde\t)\\
~~~~~~~~~~-\k(\v)\Big(\partial_{zz}\chi+2\partial_z\chi\partial_z+\f12\partial_z(\ln|g|)\cdot\partial_z\chi\Big)
\cdot(n\cdot\widetilde{\nabla\t}-\nu\tilde\t), \\
F_\kappa^\t=(R+1) \widetilde{\nabla\t}\cdot\Big(\tilde{\r} \tilde{u}^1\partial_{y^1}n +\tilde{\r} \tilde{u}^2\partial_{y^2}n\Big)+\k(\v) n\cdot\Delta_{\tilde{g}}\widetilde{\nabla\t}+\nu\k(\v)\Delta_{\tilde{g}}\tilde{\t}. 
\end{cases}
\end{align}}
Then, using the similar arguments in Lemma \ref{lem6.2}, one can obtain that
{\small\begin{align}\label{6.113}
&\|\zeta(t)\|^2_{\MH^{1,\i}}\leq C\Big\{P(\mathcal{N}_m(0))+P(\Lambda_{m})+tP(\mathcal{N}_m(t))
+P(\mathcal{N}_m(t))\Big(\int_0^t\k(\v)\|\Delta\t\|^{\f12}_{\MH^4}+\|\nabla p\|_{\MH^{2,\i}}d\tau\Big)^2\nonumber\\
&~~~~~+\Big(\int_0^t\v\|\nabla u\|_{L^\infty}\|\nabla^2 u\|_{\MH^{1,\i}}d\tau\Big)^2\Big\} \leq C\Big\{P(\mathcal{N}_m(0))+tP(\mathcal{N}_m(t))\nonumber\\
& ~~~~~~~~~+\Big(\|\nabla u_0\|^2_{L^\infty}+t\mathcal{N}_m(t)\Big)\Big(\int_0^t\v\|\nabla^2 u\|_{\MH^{1,\i}}d\tau\Big)^2\Big\}\leq CC_{m+2}\Big\{ P(\mathcal{N}_m(0))+tP(\mathcal{N}_m(t))\Big\},
\end{align}}
where we have used  \eqref{6.101}, \eqref{6.115}, \eqref{5.39} and Holder inequality above. Then \eqref{6.108} follows from \eqref{6.112}, \eqref{6.113}, \eqref{6.1} and \eqref{5.39}.  Therefore, the proof of Lemma \ref{lem3.19} is completed. $\hfill\Box$

\

Combining Lemma \ref{lem3.18}-Lemma \ref{lem3.19}, one can  obtain
\begin{proposition}
Assume  \eqref{1.8} holds, one has, for  $m\geq6$,  that
{\small\begin{align}\label{6.135}
&\|\nabla(\r\t)\|^2_{\MH^{1,\infty}}+\|\nabla u\|^2_{\mathcal{H}^{1,\infty}}+\|\nabla\t\|^2_{\MH^{1,\infty}}+\v\|\nabla (\r\t)\|^2_{\MH^{2,\infty}}+\v\|\partial_{zz}u\|^2_{L^\infty}+\int_0^t\|\nabla (\r\t)\|^2_{\MH^{2,\infty}}d\tau\nonumber\\
&\leq C C_{m+2}\Big\{ P(\mathcal{N}_m(0))  +  tP(\mathcal{N}_{m}(t))\Big\}.
\end{align}}
\end{proposition}

\subsection{Proof of Theorem \ref{thm3.1}}
Firstly, it follows from \eqref{6.4-3}, \eqref{5.39} and \eqref{6.135} that
{\small\begin{equation}\label{7.30}
Q(t) \leq C C_{m+2}\Big\{ P(\mathcal{N}_m(0))  +  tP(\mathcal{N}_{m}(t))\Big\}.
\end{equation}}
In order to close the priori estimate, one still needs to get the uniform estimate for $\|\nabla\partial_t^{m-1}u\|$.  It follows from Lemma \ref{lem4.7}, \eqref{7.30} and \eqref{5.39} that
{\small\begin{align}\label{6.137}
\|\nabla\partial_t^{m-1}u\|^2&\lesssim C_{m+1}\Big\{ \|u\|^2_{\mathcal{H}^m}+\|\eta\|^2_{\mathcal{H}^m}
+ \|\partial_t^{m-1}\mbox{div}u\|^2_{L^2}\Big\}\nonumber\\
&\lesssim C_{m+2}\Big\{P(\mathcal{N}_{m}(0)) +  tP(\mathcal{N}_{m}(t))\Big\}.
\end{align}}
Combining \eqref{5.39}, \eqref{6.135} and \eqref{6.137}, one gets \eqref{3.0-1}. Finally, it follows from \eqref{1.1} that
{\small\begin{align*}\label{7.33}
|\r(x,0)|\exp(-\int_0^t\|\mbox{div}u\|_{L^\infty}d\tau)\leq\r(x,t)\leq |\r(x,0)|\exp(\int_0^t\|\mbox{div}u\|_{L^\infty}d\tau),
\end{align*}}
so we proved \eqref{3.0-3}. The Newton-Leibniz formula yields immediately that \eqref{3.0-5}.  Thus  the proof of Theorem \ref{thm3.1} is completed. $\hfill\Box$


\section{Proof of Theorem \ref{thm1.1}}
\renewcommand{\theequation}{\arabic{section}.\arabic{equation}}

\noindent\textbf{Proof of Theorem \ref{thm1.1}}: In this section, we shall show that how we can combine  our a priori estimates to obtain the uniform existence result. Let us fix $m\geq6$, we consider initial data such that
{\small\begin{equation}\label{8.1}
\mathcal{I}_m(0)=\sup_{\v\in(0,1]}\|(\r_0^\v, u_0^\v, \t_0^\v)\|^2_{X^\v_m}\leq \tilde{C}_0,~~\mbox{and}
~~0<\hat{C}_0^{-1}\leq \r_0^{\v}\leq \hat{C}_0,~~0<\hat{C}_0^{-1}\leq \t_0^{\v}\leq \hat{C}_0,
\end{equation}}
For such initial data,   we are not aware of a local existence result for \eqref{1.1} and \eqref{1.9}, so one first needs to prove the local existence result for \eqref{1.1} and \eqref{1.9} with initial data  $(\r^\v_0, u^\v_0,\t^\v_0)\in X^{\v,m}_{NS}$. For such initial data $(\r^\v_0, u^\v_0,\t^\v_0)$, it follows from the definition of $X^{\v,m}_{NS}$ there exists a sequence of smooth approximate  initial data $(\r^{\v,\d}_0, u^{\v,\d}_0,\t^{\v,\d}_0)\in X^{\v,m}_{NS,ap}$($\delta$ being a regularization parameter), which have enough space regularity so that the time derivatives at the initial time can be defined by Navier-Stokes equations and the boundary  compatibility conditions are satisfied. For fixed $\v\in(0,1]$, we construct  approximate solutions, inductively, as follows\\
(1) Define $u^0=u_0^{\v,\d}$, and \\
(2) Assuming that $u^{k-1}$ was defined for $k\geq1$, let $(\r^k,u^k,\t^k)$ be the unique solution to the following linearized  initial boundary value problem:
{\small\begin{eqnarray}\label{8.5}
\begin{cases}
\r^k_t+\mbox{div}(\r^k u^{k-1})=0 ~~\mbox{in}~(0,T)\times\Omega,\\
\r^k u^k_t+\r^k u^{k-1}\cdot\nabla{u}^k+R\nabla{(\r^k\t^k)}=\v\Delta u^k+\v\nabla\mbox{div}u^k,~~\mbox{in}~(0,T)\times\Omega,\\
\r^k \t^k_t+\r^k u^{k-1}\cdot\nabla\t^k+R\r^k\t^k \mbox{div}u^{k-1}=\k(\v)\Delta\t^k+2\mu\v|Su^{k-1}|^2+\l\v|\mbox{div}u^{k-1}|^2,
~\mbox{in}~(0,T)\times\Omega,\\
(\r^k, u^k, \t^k)|_{t=0}=(\r_0^{\v,\delta}, u_0^{\v, \delta}, \t_0^{\v, \delta}),~~\mbox{with}~\f2{3C_0}\leq \r_0^{\v,\d},~ \t_0^{\v, \delta}\leq \f32C_0,\\
\mbox{with boundary conditions } \eqref{1.9}.
\end{cases}
\end{eqnarray}}
Since $\r^k, \t^k$ and $u^k$ are decoupled,  the existence of global unique  smooth solution $(\r^k,u^k,\t^k)(t)$ of \eqref{8.5} with $0<\r^k(t),\t^k(t)<\infty$ can be obtained by using classical methods, for example,  the same argument of in Cho and Kim \cite{Kim} and the standard elliptic regularity results as in Agmon-Douglis-Nirenberg \cite{Nirinberg}.   On the other hand, since $(\r_0^{\v,\delta}, u_0^{\v, \delta},\t_0^{\v,\delta})\in H^{3m}$, one prove that there exists  a positive time $\tilde{T}_1=\tilde{T}_1(\v)$  such that
{\small\begin{equation}\label{8.6-1}
\|(\r^{k}, u^{k}, \t^k)(t)\|^2_{H^{3m}}\leq \hat{C}_1<\infty,
~~\mbox{and}~~(2\hat{C}_0)^{-1}\leq \r^k(t),~\t^k(t)\leq 2\hat{C}_0,~~\mbox{for}~~0\leq t\leq \tilde{T}_1,
\end{equation}}
where $\tilde{T}_1$, $\hat{C}_1$ depend  on $\hat{C}_0$, $\v^{-1}$ and $\|(\r_0^{\v,\delta}, u_0^{\v, \delta}, \t_0^{\v,\delta})\|_{H^{3m}}$.
Based on the above uniform estimates for $(\r^k, u^k,\t^k)$,  by the same arguments as  section 3 of \cite{Kim}, there exists a  uniform time $\hat{T}_{1}(\leq \tilde{T}_1)$(independent of $k$) such that $(\r^k, u^k,\t^k)$ converges to a limit $(\r^{\v,\d}, u^{\v,\d},\r^{\v,\d})$ as $k\rightarrow +\infty$ in the following strong sense:
{\small$$(\r^k, u^k, \t^k)\rightarrow (\r^{\v,\d}, u^{\v,\d}, \t^{\v,\d})~~\mbox{in}~~L^\infty(0,\hat{T}_{1}; L^2),~~\mbox{and}~~\nabla u^k\rightarrow \nabla u^{\v,\d}~\mbox{in}~ L^2(0,\hat{T}_{1}, L^2).$$}
It is easy to check $(\r^{\v,\d}, u^{\v,\d},\t^{\v,\d})(t)$ is a classical solution to the problem \eqref{1.1}, \eqref{1.9} with initial data $(\r_0^{\v,\delta}, u_0^{\v, \delta}, \t_0^{\v,\d})$.  Then, by virtue of the lower semi-continuity of norms, one can deduce from   \eqref{8.6-1} that
{\small\begin{equation}\label{8.6-2}
\|(\r^{\v,\d}, u^{\v,\d},\r^{\v,\d})(t)\|^2_{H^{3m}}\leq \hat{C}_1<\infty
~~\mbox{and}~~(2\hat{C}_0)^{-1}\leq \r^{\v,\d}(t),~\t^{\v,\d}(t)\leq 2\hat{C}_0,~~\mbox{for}~~0\leq t\leq \hat{T}_{1},
\end{equation}}

Applying the a priori estimates given in Theorem \ref{thm3.1} to  $(\r^{\v,\d}, u^{\v,\d},\t^{\v,\d})(t)$, we obtain a uniform time $T_0>0$ and positive constant $\tilde{C}_3$(independent of $\v$ and  $\d$),  such that it holds for $(\r^{\v,\d}, u^{\v,\d}, \t^{\v,\d})(t)$ that
{\small\begin{align}\label{8.6}
\Upsilon_m(\r^{\v,\d}, u^{\v,\d},\t^{\v,\d})(t)\leq \tilde{C}_3~~\mbox{and}~~(2\hat{C}_0)^{-1}\leq\r^{\v,\d}(t),~\t^{\v,\d}(t)\leq 2\hat{C}_0,~\forall t\in[0,\tilde{T}_2],
\end{align}}
where  $\tilde{T}_2\triangleq\min\{T_0,\hat{T}_{1}\}$ and the uniform constants $T_0$, $\tilde{C}_3$(independent of $\v,\d$) depend only on $\hat{C}_0$ and $\mathcal{I}_m(0)$.  Based on the uniform estimates \eqref{8.6}  for $(\r^{\v,\d}, u^{\v,\d},\t^{\v,\d})$, one can pass the limit $\d\rightarrow0$ to get a strong solution $(\r^\v,u^\v, \t^\v)$ of \eqref{1.1}, \eqref{1.9} with initial data $(\r_0^\v, u_0^\v,\t^{\v}_0)$ satisfying  \eqref{8.1} by using a strong compactness arguments.
It follows from \eqref{8.6} that $(\r^{\v,\d}, u^{\v,\d},\t^{\v,\d})$ is uniform bounded in $L^\infty([0,T_0];H^m_{co})$, $\nabla(\r^{\v,\d}, u^{\v,\d},\t^{\v,\d})$ is uniform bounded in $L^\infty([0,T_0];H^{m-1}_{co})$, and $\partial_t(\r^{\v,\d}, u^{\v,\d},\t^{\v,\d})$ is uniform bounded in $L^\infty([0,T_0];H^{m-1}_{co})$. Then, it follows from the compactness argument \cite{Simon}  that $(\r^{\v,\d}, u^{\v,\d},\t^{\v,\d})$ is compact in $\mathcal{C}([0,T_0];H^{m-1}_{co})$. In particular, there exists a sequence $\d_n\rightarrow 0+$ and $(\r^\v,u^\v,\t^\v)\in\mathcal{C}([0,T_0];H^{m-1}_{co})$ such that
{\small\begin{equation}\label{9.3-1}
(\r^{\v,\d_n},u^{\v,\d_n},\t^{\v,\d_n})\rightarrow (\r^\v,u^\v,\t^\v)~~\mbox{in}~~\mathcal{C}([0,T];H^{m-1}_{co})~~\mbox{as}~~\d_n\rightarrow 0+.
\end{equation}}
Moreover, applying the lower semi-continuity of norms to   the bounds \eqref{8.6}, one obtains the bounds   for $(\r^\v,u^\v,\t^\v)$ that
{\small\begin{align}\label{8.7}
	\Upsilon_m(\r^{\v}, u^{\v},\t^{\v})(t)\leq \tilde{C}_3~~\mbox{and}~~(2\hat{C}_0)^{-1}\leq\r^{\v}(t),~\t^{\v}(t)\leq 2\hat{C}_0,~\forall t\in[0,\tilde{T}_2].
	\end{align}}
It follows from \eqref{8.7} and  the anisotropic Sobolev inequality \eqref{3.3}  that
{\small\begin{align}\label{9.4-1}
&\sup_{t\in[0,T_0]}\|(\r^{\v,\d_n}-\r^\v,u^{\v,\d_n}-u^\v,\t^{\v,\d_n}-\t^\v)\|^2_{L^\infty}\leq \sup_{t\in[0,T_0]}\|(\r^{\v,\d_n}-\r^\v,u^{\v,\d_n}-u^\v,\t^{\v,\d_n}-\t^\v)\|^2_{H^2_{co}}\nonumber\\
&+ \sup_{t\in[0,T_0]}\Big( \|\nabla(\r^{\v,\d_n}-\r^\v,u^{\v,\d_n}-u^\v,\t^{\v,\d_n}-\t^\v)\|_{H^1_{co}}\cdot\|(\r^{\v,\d_n}-\r^\v,u^{\v,\d_n}-u^\v,\t^{\v,\d_n}-\t^\v)\|_{H^2_{co}}\Big)\rightarrow 0,\nonumber\\
&~~~~~~~~~~~~~~~~~~~~~~~~~~~~~~~~~~~~~~~~~~~~~~~~~~~~~~~~~~~~~~~~~~~~~~~~~~~~~~~~~~~~~~~\mbox{as}~\d_n\rightarrow 0+.
\end{align}}
Then, it is easy to check that $(\r^\v, u^\v,\t^\v)$ is a strong solution of the Navier-Stokes system.  The uniqueness of the solution $(\r^\v, u^\v,\t^\v)$ is easy since we work on functions with Lipschitz regularity. Thus the whole family of $(\r^{\v,\d}, u^{\v,\d},\t^{\v,\d})$ converges to $(\r^{\v}, u^{\v},\t^{\v})$ as $\d\rightarrow0+$. Therefore, for initial data  $(\r^\v_0, u^\v_0,\t^\v_0)\in X^{\v,m}_{NS}$, we have established the local existence result for \eqref{1.1},\eqref{1.9} such that $(\r^\v,u^\v,\t^\v)(t)\in X^{\v,m}_{NS},~~t\in [0,\tilde{T}_2]$ .

We shall use the above local existence results  to prove Theorem \ref{thm1.1}. If $T_0\leq \hat{T}_1$, then  Theorem \ref{thm1.1} follows from \eqref{8.7} with $\tilde{C}_1\triangleq \tilde{C}_3$.  On the other hand, if $ \hat{T}_1\leq T_0$, based on the uniform estimates  \eqref{8.7}, we can use the local existence results established above to extend our solution  step by step to the uniform time interval $t\in[0,T_0]$. Therefore,   the proof of Theorem \ref{thm1.1} is completed. $\hfill\Box$

\section{Proof of Theorem \ref{thm1.2}: Vanishing Dissipation Limit}
In this section, we  study  the vanishing dissipation limit of the solutions of full compressible Navier-Stokes system \eqref{1.1} to the solutions of full compressible Euler system  with a rate of convergence. It is well known that the solution $(\r,u,\t)(t)\in H^3$ of the Euler system  \eqref{1.7}, \eqref{1.7-2} and \eqref{1.20} satisfies
{\small\begin{equation}\label{15.1}
\sum_{k=0}^3\|(\r,u,\t)\|_{C^k(0,T_1;H^{3-k})}\leq \tilde{C}_4,~~~\f{1}{2\hat{C}_0}\leq \r(t),\t(t)\leq 2\hat{C}_0,
\end{equation}}
where $\tilde{C}_4$ depends only on  $\|(\r_0,u_0,\t_0)\|_{H^3}$. On the other hand,  it follows from Theorem \ref{thm1.1} that the solution $(\r^\v,u^\v,\t^\v)(t)$ of \eqref{1.1},\eqref{1.9} and \eqref{1.20} satisfies
{\small\begin{equation}\label{15.2}
\|(\r^\v,u^\v,\t^\v)(t)\|_{X_{m}^{\v}}\leq \tilde{C}_1, ~~~\f{1}{2\hat{C}_0}\leq \r^{\v}(t),\t^\v(t)\leq 2\hat{C}_0,~~\forall~t\in[0,T_0],
\end{equation}}
where $T_0,~\hat{C}_0$, and $\tilde{C}_1$  are defined in Theorem \ref{thm1.1}. In particular,  this uniform regularity implies  the following bound
{\small\begin{equation}\label{15.3}
\|(\r^\v,u^\v,\t^\v)\|_{W^{1,\infty}}+\|\partial_t(\r^\v,u^\v,\t^\v)\|_{L^\infty}\leq \tilde{C}_1,
\end{equation}}
which plays an important role in the proof of Theorem \ref{thm1.2}.

\

Define
{\small\begin{eqnarray}\label{14.1}
\phi^\v=\r^\v-\r,~~\psi^\v=u^\v-u, ~\xi^\v=\t^\v-\t.
\end{eqnarray}}
It then follows from \eqref{1.1} and \eqref{1.7} that
{\small\begin{equation}\label{14.2}
\begin{cases}
\phi^\v_t+\r\mbox{div}\psi^\v+u\cdot\nabla\phi^\v=R_1^\v,\\
\r\psi^\v_t+\r u\cdot\nabla\psi^\v+\nabla(p^\v-p)+\Phi^\v=-\mu\v\nabla\times(\nabla\times\psi^\v)
+(2\mu+\l)\v\nabla\mbox{div}\psi^\v+R^\v_2,\\
\r\xi^\v_t+\r u\cdot\nabla\xi^\v+p\mbox{div}\psi^\v+\Psi^\v=\k(\v)\Delta\xi^\v+2\mu\v|Su^\v|^2+\l\v|\mbox{div}u^\v|^2
+R^\v_3,
\end{cases}
\end{equation}}
where
{\small\begin{equation}\label{14.3}
\begin{cases}
R_1^\v=-\phi^\v\mbox{div}\psi^\v-\psi^\v\cdot\nabla\phi^\v-\phi^\v\mbox{div}u-\nabla\r\cdot\psi^\v,\\
R^\v_2=-\phi^\v\psi^\v_t-\phi^\v u_t+\mu\v\Delta{u}+(\mu+\l)\v\nabla\mbox{div}u,\\
R^\v_3=-\phi^\v\xi^\v_t-(p^\v-p)\mbox{div}\psi^\v-(p^\v-p)\mbox{div}\psi+\k(\v)\Delta\t,
\end{cases}
\end{equation}}
and
\begin{align}\label{14.4}
\begin{cases}
\Phi^\v= (\r^\v u^\v-\r u)\cdot\nabla u^\v=(\r^\v u^\v-\r u)\cdot\nabla\psi^\v
+(\r^\v u^\v-\r u)\cdot\nabla u,\\
\Psi^\v=(\r^\v u^\v-\r u)\cdot\nabla \t^\v=(\r^\v u^\v-\r u)\cdot\nabla \xi^\v+(\r^\v u^\v-\r u)\cdot\nabla\t.
\end{cases}
\end{align}
The boundary conditions to \eqref{14.2} are
\begin{equation}\label{16.2}
\psi^\v\cdot n=0,~~n\times(\nabla\times\psi^\v)=[B\psi^\v]_\tau+[Bu]_\tau-n\times\omega,
~\mbox{and}~\nabla\xi^\v\cdot n=\nu\xi^\v+\nu\t-\nabla\t\cdot n~~\mbox{on}~\partial\Om.
\end{equation}

\begin{lemma}\label{lem14.1}
	It holds that
	\begin{eqnarray}\label{14.6}
	&&\|(\phi^\v, \psi^\v,\xi^\v)(t)\|^2
	+\int_0^t\v\|\psi^\v\|^2_{H^1}+\k(\v)\|\xi^\v\|^2_{H^1}d\tau \leq C(\v^{\f32}+\k(\v)^{\f32}),~
	t\in [0,T_2],
	\end{eqnarray}
	where $T_2=\min\{T_0,T_1\}$,  $C>0$ depend only on $\hat{C}_0, \tilde{C}_1$ and $\tilde{C}_4$.
\end{lemma}
\noindent\textbf{Proof}: Multiplying $\eqref{14.2}_2$ by $\psi^\v$, one obtains that
{\small\begin{align}\label{14.7}
&\frac{d}{dt} \int_{\Omega}\f12\r|\psi^\v|^2dx+\int_{\Omega}\Phi^\v\cdot \psi^\v dx+\int_{\Omega}\nabla(p^\v-p)\cdot \psi^\v dx\nonumber\\
&=-\mu\v\int_{\Omega}\nabla\times(\nabla\times\psi^\v)\cdot \psi^\v dx+(2\mu+\l)\v\int_{\Omega}\nabla\mbox{div}\psi^\v\cdot\psi^\v dx+\int_{\Omega}R_2^\v\cdot \psi^\v dx.
\end{align}}
If follows from integrating by parts and $\eqref{14.2}_1$ that
{\small\begin{align}\label{17.4}
&\int_{\Omega}\nabla(p^\v-p)\cdot \psi^\v dx=-\int_{\Omega}(p^\v-p)\mbox{div}\psi^\v dx\geq -R\int_{\Omega}\t\phi^\v\mbox{div}\psi^\v dx-R\int_{\Omega}\r\xi^\v \mbox{div}\psi^\v dx-C\|(\phi^\v,\xi^\v)\|^2\nonumber\\
&\geq R\f{d}{dt}\int_{\Omega}\f{\t}{2\r}|\phi^\v|^2dx-R\int_{\Omega}\r\xi^\v \mbox{div}\psi^\v dx-C\|(\phi^\v,\xi^\v)\|^2.
\end{align}}
It is easy to check that
{\small\begin{align}\label{14.10}
&-\mu\v\int_{\Omega}\nabla\times(\nabla\times\psi^\v)\cdot \psi^\v dx=-\mu\v\int_{\Omega}|\nabla\times\psi^\v|^2dx
-\mu\v\int_{\partial\Omega}n\times(\nabla\times\psi^\v)\cdot \psi^\v dx\nonumber\\
&\leq -\mu\v\|\nabla\times\psi^\v\|^2
+C\v\left|\int_{\partial\Omega}[B\psi^\v+Bu-n\times\omega]\cdot \psi^\v dx\right|\nonumber\\
&\leq
-\mu\v\|\nabla\times\psi^\v\|^2+C\v\left(|\psi^\v|^2_{L^2(\partial\Omega)}
+|\psi^\v|_{L^2(\partial\Omega)}\right),
\end{align}}
{\small\begin{equation}\label{14.11}
\v\int_{\Omega}\nabla\mbox{div}\psi^\v\cdot\psi^\v dx
=-\v\|\mbox{div}\psi^\v\|^2,
\end{equation}}
{\small\begin{equation}\label{14.8}
	\left|\int_{\Omega}\Phi^\v\cdot\psi^\v dx\right|=\left|\int_{\Omega}((\r^\v u^\v-\r u)\cdot\nabla) u^\v \cdot\psi^\v dx\right|
	\leq C(1+\|(\r^\v,u^\v,\nabla u^\v)\|_{L^\infty})\|(\phi^\v,\psi^\v)\|_{L^2}^2,
	\end{equation}}
and
{\small\begin{align}\label{14.12}
\left|\int_{\Omega}R_2^\v\cdot \psi^\v dx \right|\leq C\|(\phi^\v,\psi^\v)\|_{L^2}^2+C\v^2.
\end{align}}
Collecting all the above estimates, one gets that
{\small\begin{align}\label{14.6-1}
&\frac{d}{dt}\int_{\Omega}\f{\r}2|\psi^\v|^2+R\f{\t}{2\r}|\phi^\v|^2dx-R\int_{\Omega}\r\xi^\v \mbox{div}\psi^\v dx
+\mu\v\|\nabla\times\psi^\v\|^2+(2\mu+\l)\v\|\mbox{div}\psi^\v\|^2\nonumber\\
&\leq
C\|(\phi^\v,\psi^\v)\|_{L^2}^2+C\v^2+C\v\left(|\psi^\v|^2_{L^2}+|\psi^\v|_{L^2}\right),
\end{align}}

On the other hand, multiplying $\f{\xi^\v}{\t}$, one can obtain that
{\small\begin{align}\label{17.1}
&\f{d}{dt}\int_{\Omega}\f{\r}{2\t}|\xi^\v|^2dx+R\int_{\Omega}\r\xi^\v\mbox{div}\psi^\v dx+\f{3\k(\v)}4\int_{\Omega}\f{|\nabla\xi^\v|^2}{\t}dx\nonumber\\
&\leq C\|(\phi^\v,\psi^\v,\xi^\v)\|^2+C\k(\v)(|\xi^\v|^2_{L^2}+|\xi^\v|_{L^2})+C(\v^2+\k(\v)^2).
\end{align}}
where we have used the facts that
{\small\begin{align}\label{17.3}
	&\k(\v)\int_{\Omega}\Delta\xi^\v\f{\xi^\v}{\t}dx
	 =-\k(\v)\int_{\Omega}\f{|\nabla\xi^\v|^2}{\t}dx+\k(\v)\int_{\Omega}\f{\xi^\v}{\t}\nabla\t\cdot\nabla\xi^\v dx+
	\k(\v)\int_{\partial\Omega}\f{\xi^\v}{\t}n\cdot\nabla\xi^\v d\sigma\nonumber\\
	&\leq -\f{3\k(\v)}4\int_{\Omega}\f{|\nabla\xi^\v|^2}{\t}dx+C\k(\v)\|\xi^\v\|^2
	 +C\k(\v)\Big|\int_{\partial\Omega}\f{\xi^\v}{\t}[\nu\xi^\v+\v\t-n\cdot\nabla\t]d\sigma\Big|\nonumber\\
	&\leq -\f{3\k(\v)}4\int_{\Omega}\f{|\nabla\xi^\v|^2}{\t}dx+C\|\xi^\v\|^2+C\k(\v)(|\xi^\v|^2_{L^2}+|\xi^\v|_{L^2}),
	\end{align}}
and
\begin{align}\label{17.2}
&|\int_{\Omega}\f{\xi^\v}{\t} \Psi^\v dx|+|\int_{\Omega}\f{\xi^\v}{\t} R^\v_3 dx|
+|\int_{\Omega}\f{\xi^\v}{\t} (2\mu\v|Su^\v|^2+\l\v|\mbox{div}u^\v|^2) dx|\nonumber\\
&\leq C \|(\phi^\v,\psi^\v,\xi^\v)\|^2+C(\v+\k(\v))\|\xi^\v\|
\leq C \|(\phi^\v,\psi^\v,\xi^\v)\|^2+C(\v^2+\k(\v)^2).
\end{align}

It follows from  \eqref{3.1}  that
{\small\begin{equation}\label{14.6-2}
	\|\psi^\v\|^2_{H^1}\leq C_1 \left(\|\nabla\times\psi^\v\|^2+\|\mbox{div}\psi^\v\|^2+\|\psi^\v\|^2\right).
	\end{equation}}
The trace theorem yields that
{\small\begin{equation}\label{14.6-3}
	|\psi^\v|^2_{L^2}\leq \d\|\nabla\psi^\v\|^2+C_\d\|\psi^\v\|^2,~ \mbox{and}~
	|\xi^\v|^2_{L^2}\leq \d\|\nabla\xi^\v\|^2+C_\d\|\xi^\v\|^2
	\end{equation}}
and
{\small\begin{align}\label{14.6-4}
	\begin{cases}
	\v|\psi^\v|_{L^2}\leq \d\v\|\nabla\psi^\v\|^2+C_\d\v\|\psi^\v\|^{\f23}\leq\d\v\|\nabla\psi^\v\|^2+ \|\psi^\v\|^2+C_\d\v^{\f32},\\
	\k(\v)|\xi^\v|_{L^2}\leq \d\k(\v)\|\nabla\xi^\v\|^2+C_\d\k(\v)\|\xi^\v\|^{\f23}\leq\d\k(\v)\|\nabla\xi^\v\|^2+ \|\xi^\v\|^2+C_\d\k(\v)^{\f32}.
	\end{cases}
	\end{align}}
Adding   \eqref{14.6-1} and \eqref{17.1} together, using \eqref{14.6-2}-\eqref{14.6-4} and choosing $\d$ suitably small, one obtains that
{\small\begin{align}\label{17.5}
	 &\frac{d}{dt}\Big(\int_{\Omega}\f{\r}2|\psi^\v|^2+R\f{\t}{2\r}|\phi^\v|^2+\f{\r}{2\t}|\xi^\v|^2dx\Big)
	+c_1\Big(\v\|\psi^\v\|^2_{H^1}+\k(\v)\|\xi^\v\|^2_{H^1}\Big)\nonumber\\
	&\leq
	C\|(\phi^\v,\psi^\v,\xi^\v)\|^2_{L^2} +C\left(\v^{\f32}+\k(\v)^{\f32}\right),
	\end{align}}
where $c_1>0$ is a positive constant independent of $\v$.
Then Gronwall's inequality yields immediately that \eqref{14.6}. Therefore, the proof of Lemma \ref{lem14.1} is completed.
$\hfill\Box$

\begin{lemma}\label{lem14.2}
	It holds that
{\small\begin{align}\label{14.13}
	&\|(\mbox{div}\psi^\v,\nabla\phi^\v,\nabla\xi^\v)(t)\|^2
	+\v\int_0^t \|\nabla\mbox{div}\psi^\v(\tau)\|^2 d\tau+\k(\v)\int_0^t\|\Delta\xi^\v\|^2d\tau\nonumber\\
	&\leq C
	 \d\int_0^t\|(\psi^\v_t,\xi^\v_t)\|^2+\v\|\psi^\v\|^2_{H^2}d\tau+C_\d\int_0^t \|(\phi^\v,\psi^\v,\xi^\v)\|^2_{H^1} d\tau+C_\d[\v^{\f12}+\k(\v)^{\f13}],
	~t\in[0,T_2],
	\end{align}}
	where $\d>0$ will be chosen  later.	
\end{lemma}
\noindent\text{\bf Proof:} Multiplying $\eqref{14.2}_2$ by $\nabla\mbox{div}\psi^\v$ leads to
{\small\begin{align}\label{14.14}
&\int_{\Omega} \left(\r\psi^\v_t+\r u\cdot\nabla\psi^\v\right)\cdot\nabla\mbox{div}\psi^\v dx+\int_{\Omega}\nabla(p^\v-p)\cdot\nabla\mbox{div}\psi^\v dx\\
&=-\mu\v\int_{\Omega}\nabla\times(\nabla\times\psi^\v)\cdot \nabla\mbox{div}\psi^\v dx
+(2\mu+\l)\v\|\nabla\mbox{div}\psi^\v\|^2+\int_{\Omega} R^\v_2\cdot\nabla\mbox{div}\psi^\v dx-\int_{\Omega}\Phi^\v\cdot\nabla\mbox{div}\psi^\v dx.\nonumber
\end{align}}

It follows from $\eqref{14.2}_1$ that
{\small\begin{align}\label{14.17-1}
&\nabla\mbox{div}\psi^\v=-\f1{\r^\v}\Big[\nabla\phi^\v_t+(u^\v\cdot\nabla)\nabla\phi^\v\Big]-\f1{\r^\v}\Big[\nabla \r^\v\mbox{div}\psi^\v-\nabla u^\v\nabla\phi^\v+\nabla(\phi^\v\mbox{div}u+\psi^\v\cdot\nabla\r)\Big].
\end{align}}
which, together with integrating by parts, yields that
{\small\begin{align}\label{14.18}
&\int_{\Omega}\nabla(p^\v-p)\cdot\nabla\mbox{div}\psi^\v dx = -R\int_{\Omega}\nabla(\r\xi^\v+\t\phi^\v+\phi^\v\xi^\v)\cdot\nabla\mbox{div}\psi^\v dx\nonumber\\
&\leq  R\int_{\Omega}\r\nabla\xi^\v  \cdot\nabla\mbox{div}\psi^\v dx+R\int_{\Omega} \t\nabla\phi^\v  \cdot\nabla\mbox{div}\psi^\v dx+R\Big|\int_{\Omega}\mbox{div}(\xi^\v\nabla\r+\phi^\v\nabla\t) \mbox{div}\psi^\v dx\Big|\nonumber\\
&~~~+R\Big|\int_{\partial\Omega}(\r\xi^\v+\t\phi^\v) n\cdot\mbox{div}\psi^\v d\sigma\Big|+C\|(\phi^\v,\xi^\v)\|^2_{H^1}\nonumber\\
&\leq -R\f{d}{dt}\int_{\Omega}\f{\t}{2\r^\v}|\nabla\phi^\v|^2dx+R\int_{\Omega} \r\nabla\xi^\v  \cdot\nabla\mbox{div}\psi^\v dx+C\|(\phi^\v,\psi^\v,\xi^\v)\|^2_{H^1}+C|(\phi^\v,\xi^\v)|_{L^2}.
\end{align}}

It follows from \eqref{6.4-1} and \eqref{1.18} that
{\small\begin{eqnarray}\label{14.21-2}
\|\nabla\mbox{div}u^\v\|_{L^\infty}+\|\nabla\mbox{div}u^\v\|_{L^2}\leq C< \infty,
\end{eqnarray}}
where $C>0$ depends only on $\tilde{C}_1$.  Integrating by parts and using the Holder inequality, one has that
{\small\begin{align}\label{14.15}
&\int_{\Omega} \left(\r\psi^\v_t+\r u\cdot\nabla\psi^\v\right)\cdot\nabla\mbox{div}\psi^\v dx\leq -\int_{\Omega} \left(\r\mbox{div}\psi^\v_t+\r u\cdot\nabla\mbox{div}\psi^\v\right)\mbox{div}\psi^\v dx
\nonumber\\
&~~~~~~~+\Big|\int_{\Omega} \left(\nabla\r \psi^\v_t+\nabla(\r u)^t\nabla\psi^\v\right)\mbox{div}\psi^\v dx\Big|+\Big|\int_{\partial\Omega}\r (u\cdot\nabla)\psi^\v\cdot n \mbox{div}\psi^\v d\sigma\Big|\nonumber\\
&\leq -\f{d}{dt}\int_{\Omega}\f\r2|\mbox{div}\psi^\v|^2dx+\d\|\psi^\v_t\|^2+C_\d \|\nabla\psi^\v\|^2+\left|\int_{\partial\Omega}\r (u\cdot\nabla)n \psi^\v \mbox{div}\psi^\v dx\right|\nonumber\\
&\leq
-\f{d}{dt}\int_{\Omega}\f\r2|\mbox{div}\psi^\v|^2dx+\d\|\psi^\v_t\|^2+C_\d \|\nabla\psi^\v\|^2+C|\psi^\v|_{L^2},
\end{align}}
{\small\begin{align}\label{14.16}
&\left|\int_{\Omega}\Phi^\v\cdot\nabla\mbox{div}\psi^\v dx\right|=\left|\int_{\Omega}[(\r^\v u^\v-\r u)\cdot\nabla\psi^\v+(\r^\v u^\v-\r u)\cdot\nabla u]\nabla\mbox{div}\psi^\v dx\right|\nonumber\\
&\leq C\|(\phi^\v,\psi^\v)\|^2_{H^1}+\left|\int_{\partial\Om}((\r^\v u^\v-\r u)\cdot\nabla u)\cdot n  \mbox{div}\psi^\v d\sigma\right|\nonumber\\
&\leq
C\left(\|(\phi^\v,\psi^\v)\|^2_{H^1}
+|(\phi^\v,\psi^\v)|_{L^2}\right),
\end{align}}
and
{\small\begin{align}\label{14.20}
&\v\left|\int\nabla\times(\nabla\times\psi^\v)\cdot \nabla\mbox{div}\psi^\v dx\right|=\v\left|\int_{\partial\Omega}n\times(\nabla\times\psi^\v)\cdot \nabla\mbox{div}\psi^\v d\sigma \right|\nonumber\\
&=\v\left|\int_{\partial\Omega}\left(B\psi^\v+Bu-n\times\omega\right)\cdot \Pi(\nabla\mbox{div}\psi^\v) d\sigma \right|\leq C\v\left(1+|\psi^\v|_{H^{\f12}}\right)|\mbox{div}\psi^\v|_{H^{\f12}}\nonumber\\
&\leq C\v \|\mbox{div}\psi^\v\|_{H^1} (1+\|\psi^\v\|_{H^1})\leq  C (\v+\|\psi^\v\|_{H^1}^2).
\end{align}}
For the term involving $R^\v_2$.  It follows from \eqref{14.21-2} and  integrating by parts  that
{\small\begin{align}\label{14.21}
&\left|\int R^\v_2\cdot\nabla\mbox{div}\psi^\v dx\right|
\leq C(1+\|\nabla\mbox{div}u^\v\|_{L^\infty})[\|\phi^\v\|\|\psi^\v_t\|
+\|(\phi^\v,\psi^\v)\|^2_{H^1}]+\v\|u\|_{H^2}\|\nabla\mbox{div}\psi^\v\|\nonumber\\
&\leq \d \|\psi^\v_t\|^2+ C_\d[\|(\phi^\v,\psi^\v)\|^2_{H^1} +\v^2]+C\v.
\end{align}}
Then the trace theorem implies that
{\small\begin{equation}\label{16.4}
|(\phi^\v,\psi^\v,\xi)|_{L^2}\leq C\Big\{\|(\phi^\v,\psi^\v,\xi^\v)\|^2_{H^1}+ \|(\phi^\v,\psi^\v,\xi^\v)\|^{\f23}_{L^2}\Big\}\leq C\Big\{ \|(\phi^\v,\psi^\v,\xi)\|^2_{H^1}+ \k(\v)^{\f12}+\v^{\f12}\Big\}.
\end{equation}}
Substituting \eqref{14.18} and \eqref{14.15}-\eqref{14.21} into \eqref{14.14} and using the \eqref{16.4}, one has that
{\small\begin{align}\label{17.6}
&\f{d}{dt}\Big(\int_{\Omega}\f\r2|\mbox{div}\psi^\v|^2+\f{R\t}{2\r^\v}|\nabla\phi^\v|^2dx\Big)-R\int_{\Omega} \r\nabla\xi^\v  \cdot\nabla\mbox{div}\psi^\v dx+\f34(2\mu+\l)\v\|\nabla\mbox{div}\psi^\v\|^2\nonumber\\
&\leq \d \|\psi^\v_t\|^2+C_\d\Big[\|(\phi^\v,\psi^\v,\xi^\v)\|^2_{H^1}+\k(\v)^{\f12}+\v^{\f12}\Big].
\end{align}}

Applying $\nabla$ to $\eqref{14.2}_3$, one can obtain that
{\small\begin{align}\label{17.7}
&\r^\v\nabla\xi^\v_t+\r^\v(u^\v\cdot\nabla)\nabla\xi^\v+p\nabla\mbox{div}\psi^\v-\k(\v)\Delta\nabla\xi^\v\nonumber\\
&=\v\nabla(2\mu|Su^\v|^2+\l|\mbox{div}u^\v|^2)+\nabla\tilde{R}^\v_3-\nabla\r^\v\xi^\v_t-\nabla(\r^\v u^\v)\nabla\xi^\v-\nabla p\mbox{div}\psi^\v,
\end{align}}
where
{\small\begin{equation}\label{17.8}
\tilde{R}^\v_3=-\xi^\v\t_t-(\r^\v u^\v-\r u)\cdot\nabla\t-(p^\v-p)\mbox{div}u^\v+\k(\v)\Delta\t.
\end{equation}}
Multiplying \eqref{17.7} by $\f{\nabla\xi^\v}{\t}$, one has that
{\small\begin{align}\label{17.9}
&\f{d}{dt}\int_{\Omega}\f{\r^\v}{2\t}|\nabla\xi^\v|^2dx+R\int_{\Omega} \r\nabla\xi^\v  \cdot\nabla\mbox{div}\psi^\v dx-\k(\v)\int_{\Omega}\Delta\nabla\xi^\v\cdot\f{\nabla\xi^\v}{\t}dx\nonumber\\
&\leq \d\v\|\psi^\v\|^2_{H^2}+\d\|\xi^\v_t\|^2+C_\d\|(\phi^\v,\psi^\v,\xi^\v)\|^2_{H^1}+C_\d(\v^2+\k(\v)^2).
\end{align}}
It follows from integrating by parts and the boundary condition \eqref{16.2} that
{\small\begin{align}\label{17.10}
-\k(\v)\int_{\Omega}\Delta\nabla\xi^\v\cdot\f{\nabla\xi^\v}{\t}dx
&\geq\k(\v)\int_{\Omega}\f{1}{\t}|\Delta\xi^\v|^2 dx-C\k(\v)\|\nabla\xi^\v\|\|\Delta\xi^\v\|
-C\k(\v)\Big|\int_{\partial\Omega}\Delta\xi^\v\f{n\cdot\nabla\xi^\v}{\t}d\sigma\Big|\nonumber\\
&\geq\f34\k(\v)\int_{\Omega}\f{1}{\t}|\Delta\xi^\v|^2 dx-C\k(\v)\|\nabla\xi^\v\|^2-C\k(\v)|\Delta\xi^\v|_{L^2}\nonumber\\
&\geq\f34\k(\v)\int_{\Omega}\f{1}{\t}|\Delta\xi^\v|^2 dx-C\k(\v)\|\nabla\xi^\v\|^2-C\k(\v)\|\Delta\xi^\v\|_{L^2}^{\f12}\|\Delta\xi^\v\|_{H^1}^{\f12}\nonumber\\
&\geq\f12\k(\v)\int_{\Omega}\f{1}{\t}|\Delta\xi^\v|^2 dx-C\|\xi^\v\|^2_{H^1}-C\k(\v)^{\f13}\Big[1+\k(\v)^2\|\nabla\Delta\xi^\v\|^2\Big].
\end{align}}
Substituting \eqref{17.10} into \eqref{17.9}, one obtains that
{\small\begin{align}\label{17.11}
	&\f{d}{dt}\int_{\Omega}\f{\r^\v}{2\t}|\nabla\xi^\v|^2dx+R\int_{\Omega} \r\nabla\xi^\v  \cdot\nabla\mbox{div}\psi^\v dx+\f12\k(\v)\int_{\Omega}\f{1}{\t}|\Delta\xi^\v|^2 dx\nonumber\\
	&\leq \d\v\|\psi^\v\|^2_{H^2}+\d\|\xi^\v_t\|^2+C_\d\|(\phi^\v,\psi^\v,\xi^\v)\|^2_{H^1}+C\k(\v)^{\f13}\Big[1+\k(\v)^2\|\nabla\Delta\xi^\v\|^2\Big]+C_\d(\v^2+\k(\v)^2).
	\end{align}}
Combining \eqref{17.6} and  \eqref{17.11}, one has that
{\small\begin{align}\label{17.12}
	 &\f{d}{dt}\Big(\int_{\Omega}\f\r2|\mbox{div}\psi^\v|^2+\f{R\t}{2\r^\v}|\nabla\phi^\v|^2+\f{\r^\v}{2\t}|\nabla\xi^\v|^2dx\Big)+\f34(2\mu+\l)\v\|\nabla\mbox{div}\psi^\v\|^2+\f12\k(\v)\int_{\Omega}\f{1}{\t}|\Delta\xi^\v|^2 dx\nonumber\\
	&\leq \d\v\|\psi^\v\|^2_{H^2}+\d\|\xi^\v_t\|^2+C_\d\Big[\|(\phi^\v,\psi^\v,\xi^\v)\|^2_{H^1}+\k(\v)^{\f12}+\v^{\f12}\Big]
	+C\k(\v)^{\f13}\Big[1+\k(\v)^2\|\nabla\Delta\xi^\v\|^2\Big].
	\end{align}}
It follows from Theorem \ref{thm1.1} and \eqref{15.1} that
{\small\begin{equation}\label{17.13}
\k(\v)^2\int_0^t\|\nabla\Delta\xi^\v\|^2 d\tau \leq C<\infty.
\end{equation}}
Then, integrating \eqref{17.12} over $[0,T_2]$ and using \eqref{17.13}, one gets \eqref{14.13}.  Thus,   the proof of Lemma \ref{lem14.2} is completed.    $\hfill\Box$

\begin{lemma}\label{lem14.3}
	It holds that
{\small	\begin{align}\label{14.22}
	 &\|\nabla\times\psi^\v\|^2+\v\int_0^t\|(\nabla\times\psi^\v)(\tau)\|^2_{H^1}d\tau\leq \d\|\nabla(\phi^\v,\psi^\v,\xi^\v)\|^2_{L^2}\nonumber \\
	& ~~~~~+C\d\int_0^t\|\psi^\v_t\|^2+\v\|\nabla^2\psi^\v\|^2d\tau+C_\d\int_0^t \|(\phi^\v,\psi^\v,\xi^\v)\|^2_{H^1}d\tau+C_\d[\v^\f{1}{2}+\k(\v)^{\f12}],	
	\end{align}}
	where $\d>0$ will be chosen later.
\end{lemma}
\noindent\text{\bf Proof:} Multiplying $\eqref{14.2}_2$ by $\nabla\times(\nabla\times\psi^\v)$ gives that
{\small\begin{align}\label{14.23}
&\int_{\Omega} \r^\v\psi^\v_t\cdot\nabla\times(\nabla\times\psi^\v) dx+\int_{\Omega}\nabla(p^\v-p)\cdot\nabla\times(\nabla\times\psi^\v) dx+\mu\v\|\nabla\times(\nabla\times\psi^\v)\|^2\nonumber\\
&=(2\mu+\l)\v\int_{\Omega}\nabla\times(\nabla\times\psi^\v)\cdot \nabla\mbox{div}\psi^\v dx+\int_{\Omega} \tilde\Phi^\v \cdot\nabla\times(\nabla\times\psi^\v) dx+\int_{\Omega} \tilde{R}^\v_2\cdot\nabla\times(\nabla\times\psi^\v) dx,
\end{align}}
where one has rewritten $\eqref{14.2}_2$ and
{\small\begin{align}
\tilde{R}_2^\v=-\phi^\v u_t+\mu\v\Delta u+(\mu+\l)\v\nabla\mbox{div}u,~~ \mbox{and}~~ \tilde\Phi^\v =\r^\v u^\v\cdot\nabla\psi^\v
+(\r^\v u^\v-\r u)\cdot\nabla u.\nonumber
\end{align}}
Integrating along the boundary, one has that
{\small\begin{align}\label{14.25}
	&\left|\int_{\Omega}\nabla(p^\v-p)\cdot\nabla\times(\nabla\times\psi^\v) dx\right|
	=\left|\int_{\partial\Omega}\nabla(p^\v-p)\cdot (n\times(\nabla\times\psi^\v)) d\sigma\right|\nonumber\\
	 &=\left|\int_{\partial\Omega}\Pi(\nabla(p^\v-p))\cdot[B\psi^\v+Bu-n\times\omega] dx\right|\leq C\left[|p^\v-p|_{H^{\f12}}|\psi^\v|_{H^{\f12}}+|p^\v-p|_{L^2}\right]\nonumber\\
	&\leq C\left[\|(\phi^\v,\psi^\v,\xi^\v)\|^2_{H^1} +|(\phi^\v,\xi^\v)|_{L^2}\right].
\end{align}}
Note that the first term on the right hand side of \eqref{14.23} has been estimated in \eqref{14.20}. It remains to estimate the other terms of \eqref{14.23}.  By the same argument as the Lemma 6.3 of \cite{Wang-Xin-Yong}, one can obtains that
{\small\begin{align}\label{14.24}
\int_{\Omega}\r^\v\psi^\v_t \cdot\nabla\times(\nabla\times\psi^\v) dx
&\geq \f{d}{dt}\left(\int_{\Omega}\f12\r^\v|\nabla\times\psi^\v|^2dx+\int_{\partial\Omega}\f12\r^\v\psi^\v B\psi^\v +\r^\v\psi^\v\cdot(Bu-n\times\omega)d\sigma  \right)\nonumber\\[1mm]
&~~~~~~~~~~~~-\d\|\psi^\v_t\|^2-C_\d\left(\|\psi^\v\|^2_{H^1}
+|\psi^\v|_{L^2}\right),
\end{align}}
{\small\begin{align}\label{14.26}
&\Big|\int_{\Omega} \tilde{R}^\v_2\cdot\nabla\times(\nabla\times\psi^\v) dx\Big| \nonumber\\
&\leq C\|(\phi^\v,\psi^\v)\|^2_{H^1}+\Big|\int_{\partial\Omega}\phi^\v u_t\cdot (n\times(\nabla\times\psi^\v)) \Big|+C\v\|u\|_{H^3}\|\psi^\v\|_{H^1}+C\v\|u\|_{H^3}|\nabla\times\psi^\v|_{L^2}
\nonumber\\
&\leq \d\v\|\nabla\times(\nabla\times\psi^\v)\|^2
+C_\d\left(\|(\phi^\v,\psi^\v)\|^2_{H^1}+|(\phi^\v,\psi^\v)|_{L^2}
+\v^{\f32} \right).
\end{align}}
and
\begin{align}\label{14.44}
\Big|\int_{\Omega} \tilde\Phi^\v\cdot\nabla\times(\nabla\times\psi^\v) dx\Big|\leq C[\|(\phi^\v,\psi^\v)\|^2_{H^1}+|(\phi^\v,\psi^\v)|_{L^2}].
\end{align}
Then, combining \eqref{14.23}, \eqref{14.25}-\eqref{14.44}, one obtains that
{\small\begin{align}\label{14.45-1}
&\f{d}{dt}\left(\int_{\Omega}\f12\r^\v|\nabla\times\psi^\v|^2dx+\int_{\partial\Omega}\f12\r^\v\psi^\v B\psi^\v +\r^\v\psi^\v\cdot(Bu-n\times\omega)d\sigma  \right)
+\f12\mu\v\|\nabla\times(\nabla\times\psi^\v)\|^2\nonumber\\
&\leq
C\d\|\psi^\v_t\|^2+C\d\v\|\nabla^2\psi^\v\|^2+C_\d\left(\|(\phi^\v,\psi^\v,\xi^\v)\|^2_{H^1}+\v^{\f12}+\k(\v)^{\f12}\right),
\end{align}}
where we have used
{\small\begin{eqnarray}\label{14.47}
\begin{cases}
|(\phi^\v,\psi^\v,\xi^\v)|_{L^2}\leq C\|(\phi^\v,\psi^\v,\xi^\v)\|^{\f12}_{H^1}\cdot\|(\phi^\v,\psi^\v,\xi^\v)\|^{\f12}
\leq  \d\|\nabla(\phi^\v,\psi^\v,\xi^\v)\|^2+C_\d(\v^{\f12}+\k(\v)^{\f12}),\\[2mm]
|(\phi^\v,\psi^\v,\xi^\v)|^2_{L^2}\leq C\|(\phi^\v,\psi^\v,\xi^\v)\|_{H^1}\cdot\|(\phi^\v,\psi^\v,\xi^\v)\|
\leq  \d\|\nabla(\phi^\v,\psi^\v,\xi^\v)\|^2+C_\d(\v^{\f32}+\k(\v)^{\f32}),
\end{cases}
\end{eqnarray}}
which are consequences of the trace theorem and  \eqref{14.6}. It follows  from \eqref{3.1-1} that
{\small\begin{align}\label{14.46}
&\|\nabla\times\psi^\v\|^2_{H^1}\leq C_1 \left(\|\nabla\times(\nabla\times\psi^\v)\|^2+\|\mbox{div}(\nabla\times\psi^\v)\|^2+\|\nabla\times\psi^\v\|^2+|n\times(\nabla\times\psi^\v)|^2_{H^{\f12}}\right)\nonumber\\
&\leq
C_1 (\|\nabla\times(\nabla\times\psi^\v)\|^2+\|\nabla\times\psi^\v\|^2+|B\psi^\v|^2_{H^{\f12}}+|(Bu)_\tau-n\times\omega|^2_{H^{\f12}}),\nonumber\\
&\leq
C_1 (\|\nabla\times(\nabla\times\psi^\v)\|^2+\|\psi^\v\|^2_{H^1}+C),
\end{align}}
Substituting \eqref{14.46} into  \eqref{14.45-1} yields that
{\small\begin{align}\label{14.45}
&\f{d}{dt}\left(\int_{\Omega}\f12\r^\v|\nabla\times\psi^\v|^2dx+\int_{\partial\Omega}\f12\r^\v\psi^\v B\psi^\v +\r^\v\psi^\v\cdot(Bu-n\times\omega)d\sigma  \right)+c_1\v\|\nabla\times\psi^\v\|^2_{H^1}
\nonumber\\
&\leq
C\d\|\psi^\v_t\|^2+C\d\v\|\nabla^2\psi^\v\|^2+C_\d\left(\|(\phi^\v,\psi^\v,\xi^\v)\|^2_{H^1}+\v^{\f12}+\k(\v)^{\f12}\right),
\end{align}}
Integrating \eqref{14.45}  over $[0,t]$ and using \eqref{14.47}, one  gets \eqref{14.22}, respectively.
Therefore, the proof of Lemma \ref{lem14.3} is completed. $\hfill\Box$

\

\noindent\textbf{Proof of Theorem \ref{thm1.2}:} It follows from  \eqref{3.1} that
{\small\begin{eqnarray}\label{14.49}
\|\psi^\v\|^2_{H^1}&\leq& C\left(\|\nabla\times\psi^\v\|^2+\|\mbox{div}\psi^\v\|^2+\|\psi^\v\|^2+
|\psi^\v\cdot n|_{H^{\f12}}\right)\nonumber\\
&\leq&C\left(\|\nabla\times\psi^\v\|^2
+\|\mbox{div}\psi^\v\|^2+\|\psi^\v\|^2\right),
\end{eqnarray}}
and
{\small\begin{eqnarray}\label{14.50}
\|\psi^\v\|^2_{H^2}&\leq& C\left(\|\nabla\times\psi^\v\|^2_{H^1}+\|\mbox{div}\psi^\v\|^2_{H^1}+\|\psi^\v\|^2_{H^1}+
|\psi^\v\cdot n|_{H^{\f32}}\right)\nonumber\\
&\leq&C\left(\|\nabla\times\psi^\v\|^2_{H^1}+\|\mbox{div}\psi^\v\|^2_{H^1}
+\|\psi^\v\|^2_{H^1}\right).
\end{eqnarray}}
While  $\eqref{14.2}_2$ and $\eqref{14.2}_3$ imply  that
{\small\begin{equation}\label{14.5}
\begin{cases}
\|\psi^\v_t\|_{L^2}^2\leq C\left(\|(\phi^\v,\psi^\v)\|_{H^1}^2
+\v^2\|\nabla^2\psi^\v\|_{L^2}^2+\v^2 \right),\\
\|\xi^\v_t\|_{L^2}^2\leq C\left(\|(\phi^\v,\psi^\v)\|_{H^1}^2
+\k(\v)^2\|\Delta\xi^\v\|_{L^2}^2+\v^2+\k(\v)^2 \right).
\end{cases}
\end{equation}}
Then, collecting  \eqref{14.22}, \eqref{14.5}, \eqref{14.49}-\eqref{14.50}, \eqref{14.13}, \eqref{14.6} and choosing $\d$ suitably small, one obtains that
{\small\begin{align}
	 &\|\nabla(\psi^\v,\phi^\v,\xi^\v)\|^2+\v\int_0^t\|\psi^\v(\tau)\|^2_{H^2}d\tau+\k(\v)\int_0^t\|\Delta\xi^\v(\tau)\|^2d\tau
	\nonumber\\
& \leq
	C\int_0^t \|\nabla(\phi^\v,\psi^\v,\xi^\v)(\tau)\|^2d\tau+C[\v^\f{1}{2}+\k(\v)^{\f13}],\nonumber
\end{align}}
which, together with the  Gronwall's inequality yields immediately that
{\small\begin{equation}\label{14.52-1}
	 \|\nabla(\psi^\v,\phi^\v,\xi^\v)\|^2+\v\int_0^t\|\psi^\v(\tau)\|^2_{H^2}d\tau+\k(\v)\int_0^t\|\Delta\xi^\v(\tau)\|^2d\tau\leq C[\v^\f{1}{2}+\k(\v)^{\f13}].
\end{equation}}
Then, \eqref{14.6} and \eqref{14.52-1} imply \eqref{1.8-0}-\eqref{1.9-0}.  On the other hand, \eqref{1.10-2} is an immediately consequences of \eqref{1.8-0}, \eqref{15.1} and \eqref{15.3}.  Thus,  the proof of Theorem \ref{thm1.2} is completed.  $\hfill\Box$


\section{Appendix}

We have the following Lemma whose proof can be found in Appendix of \cite{Wang-Xin-Yong}:
\begin{lemma}\label{lemA.2}
Consider $h$ a smooth solution of
{\small\begin{align}\label{A.13}
\begin{cases}
a(t,y)[\partial_th+b_1(t,y)\partial_{y^1}h+b_2(t,y)\partial_{y^2}h+z b_3(t,y)\partial_zh]-\v\partial_{zz}h=G,~z>0,\\
h(0,y,z)=h_0(y,z),~~
~h(t,y,0)=0,
\end{cases}
\end{align}}
for some smooth function $a(t,y)$ satisfies $c_1\leq a(t,y)\leq \f1{c_1}$ and vector fields $b=(b_1,b_2,b_3)^t(t,y)$. Assume that $h$ and $G$ are compactly supported in $z$. Then, one has the estimate:
{\small\begin{equation}\label{A.14}
\|h\|_{\mathcal{H}^{1,\infty}} \lesssim \|h_0\|_{\mathcal{H}^{1,\infty}}
+\int_0^t\|\f{1}{a}\|_{L^\infty} \|G\|_{\mathcal{H}^{1,\infty}}d\tau
+\int_0^t(1+\|(\f{1}{a},b)\|_{L^\infty})(1+\sum_{i=0}^2\|Z_i(a,b)\|^2_{L^\infty}) \|h\|_{\mathcal{H}^{1,\infty}}d\tau.
\end{equation}}

\end{lemma}

\noindent {\bf Acknowledgments:} The author would like to thank Prof. Dehua Wang for valuable  discussions during his stay at AMSS, CAS. Yong Wang is partially supported by National Natural
Sciences Foundation of China No. 11371064 and 11401565.

\end{document}